\documentclass[10pt]{article}

\usepackage{amstext}
\usepackage{amssymb}
\usepackage{amsmath}
\usepackage{amsthm}

\newtheorem{theorem}{Theorem}[section]
\newtheorem{corol}[theorem]{Corollary}

\newcounter{notezz}[section]

\newtheorem{lemma}[theorem]{Lemma}

\theoremstyle{definition}
\newtheorem{definition}[theorem]{Definition}

\theoremstyle{remark}\newtheorem{remark}[theorem]{Remark}

\renewcommand{\Im}{\operatorname{Im}}
\renewcommand{\Re}{\operatorname{Re}}

\def\sgn{\mathop{\rm sign}\nolimits}

\newcommand{\SCF}{\operatorname{{\rm V}}}

\def\sgn{\mathop{\rm sign}\nolimits}
\def\eqbd{\mathop{{:}{=}}}

\begin{document}
\vspace{\baselineskip}

\title{Generalized Hurwitz polynomials}

\author{Mikhail Tyaglov\thanks{The work was supported by the
Sofja Kovalevskaja Research Prize of Alexander von Humboldt
Foundation.
Email: {\tt tyaglov@math.tu-berlin.de}}  \\
\small Institut f\"ur Mathematik,  Technische Universit\"at
Berlin }

\date{\small \today}

\maketitle

\begin{abstract}
We describe a wide class of polynomials, which is a natural
generalization of Hurwitz stable polynomials. We also give a
detailed account of so-called self-interlacing polynomials, which
are dual to Hurwitz stable polynomials but have only real and
simple zeroes. All proofs are given using properties of rational
functions mapping the upper half-plane of the complex plane to the
lower half-plane. Matrices with self-interlacing spectra and other
applications of generalized Hurwitz polynomials are discussed.
\end{abstract}

\null \thispagestyle{empty} \tableofcontents \vfill
\eject


\section*{Introduction}

\hspace{4mm} The polynomials with zeroes in the open left
half-plane of the complex plane are called Hurwitz stable, or just
stable. They play an important role in several areas of
mathematics and engineering such as stability and bifurcation
theory, control theory etc. The Hurwitz stable polynomials are
well-studied and closely related to the moment problem, spectral
theory, operator theory and orthogonal polynomials. The Hurwitz
stable polynomials were studied by C.\,Hermite, E.\,Routh,
A.\,Stodola,
A.\,Hurwitz~\cite{Hermite,Routh,Hurwitz,KreinNaimark,Gantmakher}.
The present paper was motivated by some problems of the
bifurcation theory and concerns the most natural (from our point
of view) generalization of Hurwitz stable polynomials in the
following way.

Let
\begin{equation*}\label{intro.poly}
p(z)=p_0(z^2)+zp_1(z^2)
\end{equation*}
be a real polynomial, where $p_0$ and $p_1$ are its even and odd
parts, respectively. Consider its associated rational function
\begin{equation*}\label{intro.Phi}
\Phi(u)=\dfrac{p_1(u)}{p_0(u)}~.
\end{equation*}
In the present paper we study polynomials whose associated
function $\Phi$ is a so-called \textit{R}-function, that is, the
function mapping the upper half-plane of the complex plane to the
lower half-plane. These functions have real, simple, and
interlacing zeroes and poles~\cite{KreinKatz,Holtz_Tyaglov}. Such
functions play a very important role in spectral theory of
self-adjoint operators~\cite{Simon_1998}, in the moment
problem~\cite{AkhiezerKrein,Akhiezer_moment,Nudelman}, in the
theory of the Stieltjes string~\cite{KreinKatz}, in the theory of
oscillation operators~\cite{KreinGantmaher}, in the theory of
Jacobi and Stieltjes continued fractions~\cite{Wall} etc. In
particular, it is very well
known~\cite{KreinNaimark,ChebotarevMeiman,Gantmakher,Barkovsky.2}
that the polynomial $p$ is Hurwitz stable if and only if its
associated function $\Phi$ is an \textit{R}-function with negative
zeroes and poles, and $\deg p_0\geqslant\deg p_1$. In
Section~\ref{section:Hurwitz.polys} of the present paper, we prove
this fact using methods of complex analysis, and then from this
fact we obtain all main properties of Hurwitz stable polynomials
and the main criteria of the polynomial stability including
connections with total positivity and Stieltjes continued
fractions. Section~\ref{subsection:quasi-stable.poly} is devoted
to polynomials with zeroes in the \textit{closed} half-plane of
the complex plane. Note that in~\cite{Asner,Kemperman} there was
proved that the \textit{infinite} Hurwitz matrix of a quasi-stable
polynomial is totally nonnegative. In
Section~\ref{subsection:quasi-stable.poly}, we proved the converse
statement, which is new.

In Section~\ref{section:self-interlacings}, we investigate
polynomials whose associated function $\Phi$ maps the upper
half-plane to the lower half-plane (with $\deg p_0\geqslant\deg
p_1$) but has only positive zeroes and poles in contrast with the
Hurwitz stable case . We call those polynomials
\textit{self-interlacing}. If $p(z)$ is a self-interlacing
polynomial, then it has only simple real zeroes which interlace
with zeroes of the polynomial $p(-z)$. More exactly, if
$\lambda_i$, $i=1,\ldots,n$, are the zeroes of the
self-interlacing polynomial $p$ (of type~I), then
\begin{equation*}\label{intro.self-interlacing.zero.distribution.I.type}
0<\lambda_1<-\lambda_2<\lambda_3<\ldots<(-1)^{n-1}\lambda_n.
\end{equation*}
This fact is very curious, because the deep structure of
self-interlacing polynomials is very close to the structure of the
Hurwitz polynomials but the former ones have only \text{real and
simple} zeroes. The main result that reveals this connection is
Theorem~\ref{Theorem.connection.Hurwitz.self-interlacing}, which
says that a polynomial $p(z)=p_0(z^2)+zp_1(z^2)$ is
self-interlacing if and only if the polynomial
$q(z)=p(-z^2)-zp_1(z^2)$ is Hurwitz stable.
Section~\ref{section:self-interlacings} is devoted to the
comprehensive description of  self-interlacing polynomials and
includes analogues of the Hurwitz criterion of stability and of
the Li\'enard and Chipart criterion of stability, and some other
criteria of self-interlacing. We also give some examples of
matrices whose characteristic polynomials are self-interlacing and
describe some numerical observations related to the
self-interlacing polynomials. Note that there are three
\textit{R}-functions connected to the self-interlacing
polynomials. Namely, let $p$ be self-interlacing. Then it has only
real and simple zeroes, so its logarithmic derivative
$\dfrac{p'}{p}$ is an \textit{R}-function (see
e.g.~\cite{Holtz_Tyaglov}). Moreover, $p$ is self-interlacing if
and only if its associated function $\Phi=\dfrac{p_1}{p_0}$ is an
\textit{R}-function with positive zeroes and poles
(Theorem~\ref{Theorem.main.self-interlacing}). Finally, $p$ is
self-interlacing if and only if the function
$-\dfrac{p(z)}{p(-z)}$ is an \textit{R}-function (see the proof of
Theorem~\ref{Theorem.main.self-interlacing}).

The full generalization of the Hurwitz stable polynomials (in
terms of the function $\Phi$) is given in
Section~\ref{section:Generalized.Hurwitz}. We describe the
polynomials whose associated function $\Phi$ is an
\textit{R}-function without any restriction on the signs of its
zeroes and poles. We call such polynomials \textit{generalized
Hurwitz} polynomials. More exactly, the polynomial $p$ of degree
$n$ is called generalized Hurwitz of order $k$ (of type~I),
$0\leqslant k\leqslant\left[\dfrac{n+1}2\right]$, if it has
exactly $k$ zeroes in the closed right half-plane, all of which
are nonnegative and simple:
\begin{equation*}\label{intro.Hurwitz.poly.def}
0\leqslant\mu_1<\mu_2<\cdots<\mu_k,
\end{equation*}
such that $p(-\mu_i)\neq0$, $i=1,\ldots,k$, and $p$ has an odd
number of zeroes, counting multiplicities, on each interval
$(-\mu_{k},-\mu_{k-1}),\ldots,(-\mu_3,-\mu_2),(-\mu_2,-\mu_1)$.
Moreover, the number of zeroes of $p$ on the interval~$(-\mu_1,0)$
(if any) is even, counting multiplicities. The other real zeroes
lie on the interval $(-\infty,-\mu_{k})$: an odd number of zeroes,
counting multiplicities, when $n=2l$, and an even number of
zeroes, counting multiplicities, when $n=2l+1$. All nonreal zeroes
of $p$ (if any) are located in the open left half-plane of the
complex plane.

We prove a generalization of Hurwitz theorem,
Theorem~\ref{Theorem.general.Hurwitz.criterion}, which, in fact,
is an analogue of Hurwitz criterion of stability for the
generalized Hurwitz polynomials. We also establish a
generalization of the Li\'enard and Chipart criterion of
stability, Theorem~\ref{Theorem.Lienard.Chipart.Generalized}. Note
that in~\cite{Korsakov} G.\,F.\,Korsakov made an attempt to prove
this generalization of the Li\'enard and Chipart criterion.
However, his methods did not allow him to prove the whole theorem
(for details, see
Section~\ref{subsection:General.Hurwitz.general.theory}) .

We also should note that V.\,Pivovarchik
in~\cite{Pivovarchik_gen_Hurw} describes some aspects of the
generalized Hurwitz polynomials, however, in different form,
replacing the lower half-plane by the upper half-plane. In
particular, he indirectly establishes that the polynomial $p$ is
generalized Hurwitz (shifted Hermite-Biehler with symmetry, in
Pivovarchik's terminology) if and only if its associated function
$\Phi$ is an \textit{R}-function. Precisely, the cases considered
in~\cite{Pivovarchik_gen_Hurw} are $\deg p_0>\deg p_1$ and $\deg
p_0=\deg p_1$ with the leading coefficients of $p_0$ and $p_1$ of
the same sign. Neither generalizations of Hurwitz theorem nor
generalizations of the Li\'enard and Chipart criterion, nor the
theory of self-interlacing polynomials, nor connections of the
generalized Hurwitz polynomials with continued fractions were
established in~\cite{Pivovarchik_gen_Hurw}. However,
V.\,Pivovarchik with H.\,Woracek
described~\cite{Pivovarchik_Worachek_1,Pivovarchik_Worachek_2} a
generalization of his shifted Hermite-Biehler polynomials with
symmetry (which are generalized Hurwitz polynomials up to rotation
of the upper half-plane to the left half-plane) to entire
functions. The main result
of~\cite{Pivovarchik_Worachek_1,Pivovarchik_Worachek_2} is an
analogue of our Theorem~\ref{Theorem.main.general.Hurwitz} for
entire functions, provided by two different methods. The method
in~\cite{Pivovarchik_Worachek_2} coincides with that
of~\cite{Pivovarchik_gen_Hurw}, where
as~\cite{Pivovarchik_Worachek_2} uses the theory of de Branges and
Pontryagin spaces. Our proof of
Theorem~\ref{Theorem.main.general.Hurwitz} differs from the
methods used by V.\,Pivovarchik and H.\,Woracek and is based on
using the properties of \textit{R}-functions.

It is curious that V.\,Pivovarchik came to his shifted
Hermite-Biehler polynomials and entire functions with symmetry due
to his investigation of spectra of quadratic operator
pencils~\cite{Pivovarchik_2007} and Sturm-Liouville operators
whose boundary conditions contain spectral
parameter~\cite{Mee_Piv_2001,Mee_Piv_2005}. The last two articles
deal with a generalization of the Regge problem, which has a
direct relation to the scattering theory. So the distribution of
zeroes of the generalized Hurwitz polynomials has a physical
interpretation.

Finally, in Section~\ref{subsection:application.bifurcation} we
explain how the generalized Hurwitz polynomials can be used in
bifurcation theory. Also, we describe how specific zeroes of a
real polynomial $p(z,\alpha)$ move with the movement of the
parameter $\alpha$ along $\mathbb{R}$ if $p(z,\alpha)$ is
generalized Hurwitz for all real $\alpha$.

\setcounter{equation}{0}

\section{Auxiliary definitions and theorems}\label{section:auxiliary.theorems}

\hspace{4mm} Consider a rational function
\begin{equation}\label{basic.rational.function}
R(z)=\dfrac{q(z)}{p(z)},
\end{equation}
where $p$ and $q$ are polynomials with complex coefficients
 \begin{eqnarray}\label{polynomial.1}
& p(z)= a_0z^n+a_1z^{n-1}+\dots+a_n,\qquad &
a_0, a_1, \dots, a_n\in\mathbb{R},\ a_0>0, \\
\label{polynomial.2} & q(z)=
b_0z^{n}+b_1z^{n-1}+\dots+b_{n},\qquad & b_0, b_1, \dots,
b_{n}\in\mathbb{R},
\end{eqnarray}
so $\deg p=n$ and $\deg q\leqslant n$. If the greatest common
divisor of $p$ and $q$ has degree $l$, $0 \leq l\leq n$, then the
rational function $R$ has exactly $r=n-l$ poles.

\subsection{Hankel matrices and Hankel minors}\label{subsection:Hankel}

\hspace{4mm} Expand the function~\eqref{basic.rational.function}
into its Laurent series at $\infty$:
\begin{equation}\label{basic.rat.func.expansion}
R(z)=s_{-1}+\dfrac{s_{0}}{z}+\dfrac{s_{1}}{z^{2}}+\dfrac{s_{2}}{z^{3}}+\cdots.
\end{equation}
Here $s_j = 0$ for $j<n-1-m$ and $s_{n-1-m}=\dfrac{b_0}{a_0}$,
where $m=\deg q$.

The sequence of coefficients of negative powers of $z$
\begin{equation*}\label{main.sequence}
s_0,s_1,s_2,\ldots
\end{equation*}
defines the infinite \textit{Hankel} matrix $S \eqbd S(R)
\eqbd\|s_{i+j}\|_{i,j=0}^{\infty}$.

\begin{definition}\label{def.Hankel.determinants}
For a given infinite sequence $(s_j)_{j=0}^\infty$, consider the
determinants
\begin{equation}\label{Hankel.determinants.1}
D_j(S)=
\begin{vmatrix}
    s_0 &s_1 &s_2 &\dots &s_{j-1}\\
    s_1 &s_2 &s_3 &\dots &s_j\\
    \vdots&\vdots&\vdots&\ddots&\vdots\\
    s_{j-1} &s_j &s_{j+1} &\dots &s_{2j-2}
\end{vmatrix},\quad j=1,2,3,\dots,
\end{equation}
i.e., the \textit{leading principal minors} of the infinite Hankel
matrix~$S$. These determinants are referred to as the
\textit{Hankel minors} or  \textit{Hankel determinants}.
\end{definition}

An infinite matrix is said to have finite rank $r$ if all its
minors of order greater than $r$ equal zero whereas there exists
at least one nonzero minor of order $r$.
Kronecker~\cite{Kronecker} proved that, for any infinite Hankel
matrix, any minor of order $r$ where $r$ is the rank  of the
matrix, is a multiple of its leading principal minor of order $r$.
This implies the following result.
\begin{theorem}[Kronecker~\cite{Kronecker}]\label{Th.Hankel.matrix.rank.2}
An infinite Hankel matrix $S=\|s_{i+j}\|_{0}^{\infty}$ has finite
rank $r$ if and~only~if
\begin{eqnarray}
D_r(S) & \neq & 0, \label{finite.rank.condition.1}  \\
D_j(S)& = & 0, \quad \hbox{\rm for all } j>r.
\label{finite.rank.condition.2}
\end{eqnarray}
\end{theorem}
The following theorem was established by Gantmacher
in~\cite{Gantmakher}.
\begin{theorem}\label{Th.Hankel.matrix.rank.1}
An infinite Hankel matrix $S=\|s_{i+j}\|_{0}^{\infty}$ has finite
rank if and only if the sum of the series
\begin{equation*}
R(z)=\dfrac{s_0}{z}+\dfrac{s_1}{z^2}+\dfrac{s_2}{z^3}+\cdots
\end{equation*}
is a rational function of $z$. In this case the rank of the matrix
$S$ is equal to the number of poles of the~func\-tion~$R$.
\end{theorem}

Theorems~\ref{Th.Hankel.matrix.rank.2}
and~\ref{Th.Hankel.matrix.rank.1} imply the following:  if the
greatest common  divisor of the polynomials $p$ and~$q$ defined
in~\eqref{polynomial.1}--\eqref{polynomial.2} has degree $l$,
$0\leqslant l\leqslant m$, then the
formul\ae~\eqref{finite.rank.condition.1}--\eqref{finite.rank.condition.2}
hold for $r=n-l$ for the rational
function~\eqref{basic.rational.function}.

Let $\widehat{D}_j(S)$ denote the following determinants
\begin{equation}\label{Hankel.determinants.2}
\widehat{D}_j(S)=
\begin{vmatrix}
    s_1 &s_2 &s_3 &\dots &s_{j}\\
    s_2 &s_3 &s_4 &\dots &s_{j+1}\\
    \vdots&\vdots&\vdots&\ddots&\vdots\\
    s_{j} & s_{j+1} & s_{j+2} & \dots &s_{2j-1}
\end{vmatrix},\quad j=1,2,3,\dots
\end{equation}
With a slight abuse of notation, we will also write $D_j(R) \eqbd
D_j(S(R))$ and $\widehat{D}_j(R) \eqbd\widehat{D}_j(S(R))$
 if the matrix $S=S(R)$ is made of the
coefficients~\eqref{basic.rat.func.expansion}
 of the function $R$.

Theorems~\ref{Th.Hankel.matrix.rank.2}
and~\ref{Th.Hankel.matrix.rank.1} have the following simple
corollaries, which will be useful later.
\begin{corol}[\cite{Holtz_Tyaglov}]\label{corol.Hankel.matrix.rank.2}
A rational function $R$ represented by the
series~\eqref{basic.rat.func.expansion} has at most $r$ poles if
and~only~if
\begin{equation*}\label{zero.pole.condition}
\widehat{D}_{j}(R)=0\quad {\rm for } \;\; j=r+1, r+2, \ldots .
\end{equation*}
\end{corol}%
\begin{corol}[\cite{Holtz_Tyaglov}]\label{corol.zero.pole}
A rational function $R$ with exactly $r$ poles represented by the
series~\eqref{basic.rat.func.expansion}  has a pole at the point
$0$ if and only if
\begin{equation*}\label{zero.pole.condition.2}
\widehat{D}_{r-1}(R)\neq0\quad\text{and}\quad
\widehat{D}_{r}(R)=0.
\end{equation*}
Otherwise, $\widehat{D}_r(R)\neq0$.
\end{corol}%

\vspace{2mm}

In the sequel, we are also interested in the number of sign
changes in the sequences of Hankel minors $D_j(R)$ and
$\widehat{D}_j(R)$, which we denote by\footnote{Here we set
$D_0(R)=\widehat{D}_0(R)\equiv1$.}
$\SCF(D_0(R),D_1(R),\ldots,D_r(R))$ and
$\SCF(\widehat{D}_0(R),\widehat{D}_1(R),\ldots,\widehat{D}_r(R))$,
respectively. These numbers exist, since $R$ is a real function.

In his remarkable work~\cite{Frobenius}, G.\,Frobenius proved the
following fact that allow us to calculate the number of sign
changes in the sequence $D_0(R),D_1(R),\ldots,D_r(R)$ in the case
when some minors $D_j(R)$ are zero.
\begin{theorem}[Frobenius \cite{Frobenius,{Gantmakher.1}}]\label{Th.Frobenius}
If, for some integers $i$ and $j$ $(0\leq i<j)$,
\begin{equation}\label{Th.Frobenius.condition}
D_i(R)\neq 0, \quad
D_{i+1}(R)=D_{i+2}(R)=\cdots=D_{i+j}(R)=0,\quad D_{i+j+1}(R)\neq
0,
\end{equation}
then the number $\SCF(D_0(R),D_1(R),D_2(R),\ldots,D_r(R))$ of
Frobenius sign changes should be calculated by assigning  signs as
follows:
\begin{equation}\label{Th.Frobenius.rule}
\sgn D_{i+\nu}(R)=(-1)^{\tfrac{\nu(\nu-1)}2}\sgn D_{i}(R),\quad
\nu=1,2,\ldots,j.
\end{equation}
\end{theorem}

\vspace{2mm}

\begin{definition}\label{def.tot.noneg.infinite.matrix}
A matrix (finite or infinite) is called \textit{totally
nonnegative} (strictly totally positive) if all its minors are
nonnegative (positive).
\end{definition}

\begin{definition}\label{def.tot.noneg.infinite.matrix.with.rank}
An infinite matrix of finite rank~$r$ is called
$m$-\textit{totally nonnegative} ($m$-\textit{strictly totally
positive}) of rank~$r$, $m\leqslant r$, if all its minors are
nonnegative (positive) up to order~$m$ inclusively.

An $r$-totally nonnegative ($r$-strictly totally positive) matrix
of rank $r$ is called \textit{totally nonnegative}
(\textit{strictly totally positive}) of rank $r$.
\end{definition}

\begin{definition}\label{def.pos.def.infinite.matrix}
An infinite matrix is called \textit{positive
\emph{(}nonnegative\emph{)} definite} of rank $r$ if all its
leading principal minors are positive (nonnegative) up to order
$r$ inclusive.
\end{definition}

For the matrix $A$ (finite or infinite), we denote its minor of
order $j$ constructed with rows $i_1,i_2,\ldots,i_j$ and columns
$l_1,l_2,\ldots,l_j$ by
\begin{equation*}\label{minor.denotation}
A\begin{pmatrix}
    i_1 &i_2 &\dots &i_j\\
    l_1 &l_2 &\dots &l_j\\
\end{pmatrix}.
\end{equation*}

\begin{definition}[\cite{KreinGantmaher}]\label{def.stricly.sign.reg.infinite.matrix}
An infinite matrix $A$ of finite rank $r$ is called
$m$-\textit{sign regular}, $m\leqslant r$, if all its minors up to
order~$m$ (inclusively) satisfy the following inequalities:
\begin{equation}\label{Sign.regular.matrix.minors.ineq}
(-1)^{\sum\limits_{k=1}^{j}i_k+\sum\limits_{k=1}^{j}l_k}
A\begin{pmatrix}
    i_1 &i_2 &\dots &i_j\\
    l_1 &l_2 &\dots &l_j
\end{pmatrix}>0,\qquad j=1,\ldots,m.
\end{equation}
An $r$-sign regular matrix of rank $r$ is called \textit{sign
regular} of rank $r$.
\end{definition}

In~\cite{Gantmakher} there was proved the following criterion of
total positivity of infinite Hankel matrices of finite order.

\begin{theorem}\label{Theorem.total.nonnetativity.Hankel.matrix}
An infinite Hankel matrix $S=\|s_{i+j}\|_{0}^{\infty}$ of finite
rank $r$ is strictly totally positive if and only if the following
inequalities hold:
\begin{equation}\label{Theorem.total.nonnetativity.Hankel.matrix.conditions}
\begin{split}
&D_j(S)>0,\\
&\widehat{D}_j(S)>0,
\end{split}\qquad j=1,\ldots,r.
\end{equation}
\end{theorem}

In~\cite{Holtz_Tyaglov}, the following fact was proved.

\begin{theorem}\label{Theorem.sign.regularity.Hankel.matrix}
An infinite Hankel matrix $S=\|s_{i+j}\|_{0}^{\infty}$ of finite
rank $r$ is sign regular if and only if the~following inequalities
hold:
\begin{equation}\label{Theorem.sign.regularity.Hankel.matrix.conditions}
\begin{split}
&D_j(S)>0,\\
&(-1)^j\widehat{D}_j(S)>0,
\end{split}\qquad j=1,\ldots,r.
\end{equation}
\end{theorem}

\subsection{Matrices and minors of Hurwitz type}\label{subsection:Hurwitz}

\hspace{4mm} In this section we introduce infinite matrices
(matrices of Hurwitz type) associated with the polynomials $p$ and
$q$ defined in~\eqref{polynomial.1}--\eqref{polynomial.2} and
discuss the Hurwitz formula\ae that connect Hankel matrices with
matrices of Hurwitz type.

\begin{definition}\label{def.Hurwitz.matrix.infinite}
Given polynomials $p$ and $q$
from~\eqref{polynomial.1}--\eqref{polynomial.2}, define the
infinite matrix $H(p,q)$ as follows:
if $\deg q < \deg p $, that is, if $b_0=0$, then
\begin{equation}\label{Hurwitz.matrix.infinite.case.1}
H(p,q)=
\begin{pmatrix}
a_0&a_1&a_2&a_3&a_4&a_5&\dots\\
0  &b_1&b_2&b_3&b_4&b_5&\dots\\
0  &a_0&a_1&a_2&a_3&a_4&\dots\\
0  &0  &b_1&b_2&b_3&b_4&\dots\\
\vdots &\vdots  &\vdots  &\vdots &\vdots  &\vdots &\ddots
\end{pmatrix},
\end{equation}
if $\deg q =\deg p$, that is, $b_0\neq0$, then
\begin{equation}\label{Hurwitz.matrix.infinite.case.2}
H(p,q)=
\begin{pmatrix}
b_0&b_1&b_2&b_3&b_4&b_5&\dots\\
0  &a_0&a_1&a_2&a_3&a_4&\dots\\
0  &b_0&b_1&b_2&b_3&b_4&\dots\\
0  &0  &a_0&a_1&a_2&a_3&\dots\\
\vdots &\vdots  &\vdots  &\vdots &\vdots  &\vdots &\ddots
\end{pmatrix}.
\end{equation}
The matrix $H(p,q)$ is called an \textit{infinite matrix of
Hurwitz type}. We denote the leading principal minor of~$H(p,q)$
of order $j$, $j=1,2,\ldots$, by $\eta_{j}(p,q)$.
\end{definition}

\begin{remark}\label{remark.1.9}
The matrix $H(p,q)$ is of infinite rank since its submatrix
obtained by deleting the even or odd rows of the original matrix
is a triangular infinite matrix with $a_0 \neq 0$ on the main
diagonal.
\end{remark}

Together with the infinite matrix $H(p,q)$, we consider its
specific finite submatrices:
\begin{definition}\label{def.Hurwitz.matrix.finite}
Let the polynomials $p$ and $q$ are given
by~\eqref{polynomial.1}--\eqref{polynomial.2}.

\noindent If $\deg q < \deg p = n$, then we construct the ${2n}
\times {2n}$ matrix
\begin{equation}\label{Hurwitz.matrix.finite.case.1}
\mathcal{H}_{2n}(p,q)=
\begin{pmatrix}
    b_1 &b_2 &b_3 &\dots &b_{n}   &   0    &0      &\dots &0&0\\
    a_0 &a_1 &a_2 &\dots &a_{n-1} & a_{n}  &0      &\dots &0&0\\
     0  &b_1 &b_2 &\dots &b_{n-1} & b_{n}  &0      &\dots &0&0\\
     0  &a_0 &a_1 &\dots &a_{n-2} & a_{n-1}&a_{n}  &\dots &0&0\\
    \vdots&\vdots&\vdots&\ddots&\vdots&\vdots&\vdots&\ddots&\vdots&\vdots\\
     0  &  0 &  0 &\dots &a_{0} & a_{1}&a_{2}&\dots &a_{n}&0\\
     0  &  0 &  0 &\dots &   0  & b_{1}&b_{2}&\dots &b_{n}&0\\
     0  &  0 &  0 &\dots &0     & a_{0}&a_{1}&\dots &a_{n-1}&a_{n}\\
\end{pmatrix},
\end{equation}

\noindent if $\deg q = \deg p =n$, then we construct the
$(2n{+}1)\times (2n{+}1)$ matrix
\begin{equation}\label{Hurwitz.matrix.finite.case.2}
\mathcal{H}_{2n+1}(p,q)=
\begin{pmatrix}
    a_0 &a_1 &a_2 &\dots &a_{n-1} & a_{n}  &0&\dots &0&0\\
    b_0 &b_1 &b_2 &\dots &b_{n-1} & b_{n}  &0&\dots &0&0\\
     0  &a_0 &a_1 &\dots &a_{n-2} & a_{n-1}&a_{n}&\dots &0&0\\
     0  &b_0 &b_1 &\dots &b_{n-2} & b_{n-1}&b_{n}&\dots &0&0\\
    \vdots&\vdots&\vdots&\ddots&\vdots&\vdots&\vdots&\ddots&\vdots&\vdots\\
     0  &  0 &  0 &\dots &a_{0} & a_{1}&a_{2}&\dots &a_{n}&0\\
     0  &  0 &  0 &\dots &b_{0} & b_{1}&b_{2}&\dots &b_{n}&0\\
     0  &  0 &  0 &\dots &0     & a_{0}&a_{1}&\dots &a_{n-1}&a_{n}\\
\end{pmatrix}   .
\end{equation}
Both kinds of matrices $\mathcal{H}_{2n}(p,q)$ and
$\mathcal{H}_{2n+1}(p,q)$  are called \textit{finite matrices of
Hurwitz type}. The leading principal minors of these matrices will
be denoted by\footnote{That is, $\nabla_j(p,q)$ is the leading
principal minor  of the matrix $\mathcal{H}_{2n}(p,q)$ of order
$j$ if $\deg q < \deg p $. Otherwise (when $\deg q =\deg p$),
$\nabla_j(p,q)$ denotes the  leading principal minor of the matrix
$\mathcal{H}_{2n+1}(p,q)$ of order $j$.} $\nabla_j(p,q)$.
\end{definition}

In his celebrated work \cite{Hurwitz}, A.~Hurwitz found
relationships between the minors $D_j(R)$ and $\widehat{D}_j(R)$
defined in~\eqref{Hankel.determinants.1}
and~\eqref{Hankel.determinants.2} and the leading principal minors
$\eta_{i}(p,q)$ of the matrix $H(p,g)$ (see
Definition~\ref{def.Hurwitz.matrix.infinite}).
\begin{lemma}[\cite{Hurwitz,{KreinNaimark},{Gantmakher},{Barkovsky.2},{Holtz_Tyaglov}}]\label{lemma.relations.eta.nabla.D}
Let the polynomials $p$ and $q$ be defined
in~\eqref{polynomial.1}--\eqref{polynomial.2} and let
\begin{equation*}
R(z)=\dfrac{q(z)}{p(z)}=s_{-1}+\dfrac{s_0}{z}+\dfrac{s_1}{z^2}+\cdots
\end{equation*}
The following relations hold between the determinants
$\eta_j(p,q)$ and $D_j(R),\widehat{D}_j(R)$ defined
in~\eqref{Hankel.determinants.1}
and~\eqref{Hankel.determinants.2}, respectively.
\begin{itemize}
\item[] If $\deg q < \deg p $, then
\begin{equation}\label{Hurwitz.determinants.relations.infinite.2.even}
\eta_{2j}(p,q)=a_0^{2j}D_j(R),\quad j=1,2,\ldots;
\end{equation}
\begin{equation}\label{Hurwitz.determinants.relations.infinite.2.odd}
\eta_{2j+1}(p,q)=(-1)^{j}a_0^{2j+1}\widehat{D}_j(R),\quad
j=0,1,2,\ldots
\end{equation}
\item[] If $\deg q = \deg p $, then
\begin{equation}\label{Hurwitz.determinants.relations.infinite.1.odd}
\eta_{2j+1}(p,q)=b_0a_0^{2j}D_j(R),\quad j=0,1,2,\ldots;
\end{equation}
\begin{equation}\label{Hurwitz.determinants.relations.infinite.1.even}
\eta_{2j}(p,q)=(-1)^{j-1}b_0a_0^{2j-1}\widehat{D}_{j-1}(R),\quad
j=1,2,\ldots
\end{equation}
\end{itemize}
\end{lemma}
From~\eqref{Hurwitz.matrix.infinite.case.1}--\eqref{Hurwitz.matrix.finite.case.2}
we have
\begin{itemize}
\item[] If $\deg q < \deg p$, then
\begin{equation}\label{Hurwitz.finite.infnite.determinants.relation.1}
\nabla_{i}(p,q)=a_0^{-1}\eta_{i+1}(p,q),\quad i=1,2,\ldots,2n;
\end{equation}
\item[] If $\deg q = \deg p $, then
\begin{equation}\label{Hurwitz.finite.infnite.determinants.relation.2}
\nabla_{i}(p,q)=b_0^{-1}\eta_{i+1}(p,q),\quad i=1,2,\ldots,2n+1.
\end{equation}%
\end{itemize}

\subsection{Rational functions mapping the upper half-plane to
the lower half-plane}\label{subsection:rat.func.map.UHPtoLHP}

\hspace{4mm} All main theorems related to the generalized Hurwitz
polynomials are based on properties of the following very
important class of rational functions:

\begin{definition}\label{def.S-function}
A rational function $F$ is called \textit{R}-function if it maps
the upper half-plane of the complex plane to the lower
half-plane\footnote{In~\cite{Holtz_Tyaglov} these functions are
called \textit{R}-functions of \textit{negative type}.}:
\begin{equation*}\label{R-function.negative.type.condition.doblicate}
\Im z>0\Rightarrow\Im F(z)<0.
\end{equation*}
\end{definition}
By now, these functions, as well as their meromorphic analogues,
have been considered by many authors and have acquired various
names. For instance, these functions are called \emph{strongly
real functions} in the monograph~\cite{Sheil-Small} due to their
property to take real values \textit{only} for real values of the
argument (more general and detailed consideration can be found
in~\cite{ChebotarevMeiman}, see also~\cite{Holtz_Tyaglov}).

The following theorem provides the most important (for the present
study) properties of \textit{R}-functions. Some parts of this
theorem was considered
in~\cite{Pick,{KreinNaimark},{ChebotarevMeiman},{Gantmakher},{Atkinson},{Barkovsky.2}}.
An extended version of this theorem can be found (with proof)
in~\cite[Theorem~3.4]{Holtz_Tyaglov}.

\begin{theorem}\label{Th.R-function.general.properties}
Let $h$ and $g$ be real polynomials such that $\deg
h-1\leqslant\deg g=n$. For the real rational function
\begin{equation*}\label{Rational.function.for.R-functions}
F=\dfrac{h}{g}
\end{equation*}
with exactly $r\,(\leqslant n)$ poles, the following conditions
are equivalent:
\begin{itemize}
\item[$1)$] $F$ is an \textit{R}-function:
\begin{equation}\label{R-function.negative.type.condition}
\Im z>0\Rightarrow\Im F(z)<0;
\end{equation}
\item[$2)$] The function $F$ can
be represented in the form
\begin{equation}\label{Mittag.Leffler.1}
\displaystyle F(z)=-\alpha z+\beta+\sum^{r}_{j=1}\frac
{\gamma_j}{z+\omega_j},\quad\alpha\geqslant0,\,\,\,
\beta,\omega_j\in\mathbb{R},
\end{equation}
where
\begin{equation}\label{Mittag.Leffler.1.poles.formula}
\gamma_j=\dfrac{h(\omega_j)}{g'(\omega_j)}>0,\quad j=1,\ldots,r;
\end{equation}
\item[$3)$] The zeroes of the polynomials $\widetilde{g}(z)$ and
$\widetilde{h}(z)$, where $\widetilde{g}=g/f$,
$\widetilde{h}=h/f$, $f=\gcd(g,h)$, are real, simple and
interlacing, that is, between any two consecutive zeroes of one of
polynomials there is exactly one zero of the other polynomial,
counting multiplicity, and
\begin{equation}\label{R-functions.normalization}
\exists\; \omega\in\mathbb{R}\,:\,\,\,
\widetilde{g}(\omega)\widetilde{h}'(\omega)-\widetilde{g}'(\omega)\widetilde{h}(\omega)<0;
\end{equation}
\item[$4)$] The polynomial
\begin{equation}\label{R-functions.linear.combination}
g(z)=\lambda p(z)+\mu q(z)
\end{equation}
has only real zeroes for any real $\lambda$ and $\mu$,
$\lambda^2+\mu^2\neq0$, and the
condition~\eqref{R-functions.normalization} is satisfied;
\item[$5)$]
Let the function $F$ be represented by the series
\begin{equation}\label{R-function.series}
F(z)=s_{-2}z+s_{-1}+\frac{s_0}z+\frac{s_1}{z^2}+\frac{s_2}{z^3}+\dots,
\end{equation}
then $s_{-2}\leqslant0$ and $s_{-1}\in\mathbb{R}$. The following
inequalities hold
\begin{equation}\label{Grommer.condition.1}
D_j(F)>0,\quad j=1,2,\dots,r,
\end{equation}
where determinants $D_j(F)$ are defined
in~\eqref{Hankel.determinants.1}.
\end{itemize}
\end{theorem}

\begin{remark}\label{remark.R-func.pos.type}
If a function $F$ maps the upper half-plane to the upper
half-plane\footnote{In~\cite{Holtz_Tyaglov} such functions are
called \textit{R}-functions of \textit{positive type}.}, then the
function~$-F$ is an \textit{R}-function.
\end{remark}

Note that from~\eqref{Mittag.Leffler.1} it follows that if $F$ is
an \textit{R}-function, then the function $F$ is decreasing
between its poles, so $-F$ is obviously increasing. Therefore,
from Theorem~\ref{Th.R-function.general.properties} we obtain the
following well-known fact.
\begin{corol}\label{Corol.monotonicity.R.functions}
Let $g$ and $h$ be real polynomials such that $|\deg h-\deg
g|\leqslant1$ and let the zeroes of the polynomials $g$ and $h$ be
real and simple. If zeroes of $g$ and $h$ interlace, then the
function $\dfrac{h}{g}$ is decreasing or increasing between its
poles.
\end{corol}

\begin{remark}\label{remark.3.2}
If $R$ defined in~\eqref{basic.rational.function} is an
\textit{R}-function, then by
Theorem~\ref{Th.R-function.general.properties}, we have $\deg
q\geqslant\deg p-1$, that is, the equality $b_0=0$ implies
$b_1\neq0$ or, in other words, $b_0^2+b_1^2\neq0$.
\end{remark}

In the sequel, we also use the following
fact~\cite{Holtz_Tyaglov}, which is a simple consequence of a
theorem of V.\,Markov~\cite{ChebotarevMeiman}.

\begin{theorem}\label{Corol.differentiation.of.R-functions}
Let $p$ and $q$ be real coprime polynomials such that $|\deg
p-\deg q|\leqslant1$. If the function $R=q/p$ is an $R$-function,
then the functions $R_j=q^{(j)}/p^{(j)}$, $j=1,\ldots,\deg p-1$,
are also $R$-functions.
\end{theorem}

The following simple properties of \textit{R}-functions can be
established using Theorem~\ref{Th.R-function.general.properties}.

\begin{theorem}\label{Theorem.properties.R-functions}
Let $F$ and $G$ be \textit{R}-functions. Then
\begin{itemize}
\item[1)] the function  $R=F+G$ is an
\textit{R}-function;
\item[2)] the function $-F(-z)$ is an \textit{R}-function;
\item[3)] the function $-\dfrac1{F(z)}$ is an \textit{R}-function;
\item[4)] the function $zF(z)$ is an \textit{R}-function if all
poles of $F$ are positive and $F(z)\to0$ as $z\to\infty$.
\end{itemize}
\end{theorem}

Now consider again the function $R$ defined
in~\eqref{basic.rational.function}.
By~\eqref{polynomial.1}--\eqref{polynomial.2}, the degree of its
numerator is not greater than the degree of its denominator. If
$R$ is an \textit{R}-function, then one can find the numbers of
its negative and positive poles (see e.g.~\cite{Holtz_Tyaglov}).
\begin{theorem}\label{Th.number.of.negative.poles.of.R-frunction}
Let a rational function $R$ with exactly $r$ poles be an
\textit{R}-function of negative type and let $R$ have a series
expansion~\eqref{R-function.series}. Then the number $r_{-}$ of
negative poles of $R$ equals\footnote{Recall that the number
$\SCF(1,\widehat{D}_1(R),\widehat{D}_2(R),\ldots,\widehat{D}_k(R))$
of Frobenius sign changes must be calculated according to
Frobenius rule provided by Theorem~\ref{Th.Frobenius}.}
\begin{equation}\label{number.negative.poles.R-functions}
r_{-}=\SCF(1,\widehat{D}_1(R),\widehat{D}_2(R),\ldots,\widehat{D}_k(R)),
\end{equation}
where $k=r-1$, if $R(0)=\infty$, and $k=r$, if $|R(0)|<\infty$.
The~determinants $\widehat{D}_j(R)$ are defined
in~\eqref{Hankel.determinants.2}.
\end{theorem}

As was shown in~\cite{Gantmakher} (see also~\cite{Holtz_Tyaglov}),
if $R$ is an \textit{R}-function has exactly $r$ poles, all of
which are of the same sign, then all minors $\widehat{D}_j(R)$ are
non-zero up to order $r$.

\begin{corol}[\cite{Gantmakher}]\label{Corol.R-functions.pos.poles}
Let a rational function $R$ with exactly~$r$ poles (counting
multiplicities) be an \textit{R}-function of negative type. All
poles of $R$ are positive if and only if
\begin{equation*}\label{Corol.R-functions.pos.poles.condition}
\widehat{D}_j(R)>0,\qquad j=1,\ldots,r,
\end{equation*}
where the determinants $\widehat{D}_j(R)$ are defined
in~\eqref{Hankel.determinants.2}.
\end{corol}

\begin{corol}[\cite{Holtz_Tyaglov}]\label{Corol.R-functions.neg.poles}
Let a rational function $R$ with exactly~$r$ poles be an
\textit{R}-function of negative type. All poles of $R$ are
negative if and only if
\begin{equation*}\label{Corol.R-functions.neg.poles.condition}
(-1)^j\widehat{D}_j(R)>0,\qquad j=1,\ldots,r,
\end{equation*}
where the determinants $\widehat{D}_j(R)$ are defined
in~\eqref{Hankel.determinants.2}.
\end{corol}

One can
use~\eqref{Hurwitz.determinants.relations.infinite.2.even}--\eqref{Hurwitz.finite.infnite.determinants.relation.2}
to obtain criteria for the function $R$ defined
in~\eqref{basic.rational.function} to be an \textit{R}-function in
terms of Hurwitz minors.

\begin{definition}\label{Def.strong.sign.changes}
For a sequence of real numbers $a_0,a_1,\ldots,a_n$, we denote the
number of sign changes in the sequence of its \textit{nonzero}
entries  by $v(a_0,a_1,\ldots,a_n)$. The number
$v(a_0,a_1,\ldots,a_n)$ is usually called the \textit{number of
strong sign changes} of the sequence $a_0,a_1,\ldots,a_n$.
\end{definition}

In~\cite{Holtz_Tyaglov}, the following theorem was established.

\begin{theorem}\label{Th.generalized.Lienard.Chipart}
Let $R$ be a real rational function as
in~\eqref{basic.rational.function}. If $R$ is an
\textit{R}-function with exactly $n$ poles, that is,
$\gcd(p,q)\equiv1$, then the number of its positive poles equals
$v(a_0,a_1,\ldots,a_n)$. In particular, $R$ has only negative
poles if and only if~\footnote{In fact, the coefficients must be
of the same signs, but we assume that $a_0>0$
(see~\eqref{polynomial.1}).} $a_j>0$ for $j=1,2,\ldots,n$, and it
has only positive poles if and only if $a_{j-1}a_j<0$ for
$j=1,2,\ldots,n$.
\end{theorem}

\begin{remark}\label{Th.generalized.Lienard.Chipart.remark}
It is easy to see that
Theorem~\ref{Th.generalized.Lienard.Chipart} is also true for an
\textit{R}-function with degree of numerator greater than degree
of denominator.
\end{remark}

For \textit{R}-functions with only negative poles, there are a few
more criteria~\cite{Holtz_Tyaglov}.
\begin{theorem}\label{Th.general.Lienard.Chipart.negative.case.1}
The function~\eqref{basic.rational.function}, where $\deg q<\deg
p$, is an \textit{R}-function and has exactly $n$ negative poles
if and only if one of the following equivalent conditions holds
\begin{itemize}
\item[$1)$]
$a_n>0,a_{n-1}>0,\ldots,a_0>0,\quad\nabla_1(p,q)>0,\nabla_3(p,q)>0,\ldots,\nabla_{2n-1}(p,q)>0$;
\item[$2)$]
$a_n>0,b_n>0,b_{n-1}>0,\ldots,b_1>0,\quad\nabla_1(p,q)>0,\nabla_3(p,q)>0,\ldots,\nabla_{2n-1}(p,q)>0$;
\item[$3)$]
$a_n>0,a_{n-1}>0,\ldots,a_0>0,\quad\nabla_2(p,q)>0,\nabla_4(p,q)>0,\ldots,\nabla_{2n}(p,q)>0$;
\item[$4)$]
$a_n>0,b_n>0,b_{n-1}>0,\ldots,b_1>0,\quad\nabla_2(p,q)>0,\nabla_4(p,q)>0,\ldots,\nabla_{2n}(p,q)>0$;
\end{itemize}
where $\nabla_i(p,q)$ are defined in
Definition~\ref{def.Hurwitz.matrix.finite}.
\end{theorem}

\noindent In the case $\deg q=\deg p$ we have similar criteria.
\begin{theorem}\label{Th.general.Lienard.Chipart.negative.case.2}
The function~\eqref{basic.rational.function}, where $\deg
q(z)=\deg p(z)$, is an \textit{R}-function and has exactly $n$
negative poles if and only if one of the following equivalent
conditions holds
\begin{itemize}
\item[$1)$]
$a_n>0,a_{n-1}>0,\ldots,a_0>0,\quad\nabla_2(p,q)>0,\nabla_4(p,q)>0,\ldots,\nabla_{2n}(p,q)>0$;
\item[$2)$]
$a_n>0,b_n>0,b_{n-1}>0,\ldots,b_0>0,\quad\nabla_2(p,q)>0,\nabla_4(p,q)>0,\ldots,\nabla_{2n}(p,q)>0$;
\item[$3)$]
$a_n>0,a_{n-1}>0,\ldots,a_0>0,\quad\nabla_1(p,q)>0,\nabla_3(p,q)>0,\ldots,\nabla_{2n+1}(p,q)>0$;
\item[$4)$]
$a_n>0,b_n>0,b_{n-1}>0,\ldots,b_0>0,\quad\nabla_1(p,q)>0,\nabla_3(p,q)>0,\ldots,\nabla_{2n+1}(p,q)>0$;
\end{itemize}
where $\nabla_i(p,q)$ are defined in
Definition~\ref{def.Hurwitz.matrix.finite}.
\end{theorem}

At last, we recall a relation between \textit{R}-functions with
nonpositive poles and some properties of Hurwitz
matrices~\cite{Holtz_Tyaglov}.

\begin{theorem}[\textbf{Total Nonnegativity of the Hurwitz Matrix}]\label{Th.Hurwitz.Matrix.Total.Nonnegativity}
The following are equivalent:
\begin{itemize}
\item[$1)$] The polynomials $p$ and $q$ defined in~\eqref{polynomial.1}--\eqref{polynomial.2} have only nonpositive
zeroes \emph{(}or $q(z)\equiv0$\emph{)}, and the function~$R$
defined in~\eqref{basic.rational.function} is either an
\textit{R}-function \emph{(}or $R(z)\equiv0$\emph{)}.
\item[$2)$] The infinite matrix of Hurwitz type $H(p,q)$ defined in~\eqref{Hurwitz.matrix.infinite.case.1}--\eqref{Hurwitz.matrix.infinite.case.2} is totally nonnegative.
\end{itemize}
\end{theorem}

\vspace{2mm}

Theorem~\ref{Th.Hurwitz.Matrix.Total.Nonnegativity} implies the
following corollary~\cite{Holtz_Tyaglov}.

\begin{corol}\label{Corol.Hurwitz.Matrix.Total.Nonnegativity.1and2}
Let the polynomials $p$ and $q$ be defined
in~\eqref{polynomial.1}--\eqref{polynomial.2}. The following are
equivalent:
\begin{itemize}
\item[$1)$] The function~$R=\dfrac{p}{q}$ is an \textit{R}-function with
exactly $n$ poles, all of which are negative.
\item[$2)$] The infinite matrix of Hurwitz type $H(p,q)$ defined in~\eqref{Hurwitz.matrix.infinite.case.1}--\eqref{Hurwitz.matrix.infinite.case.2} is totally nonnegative and $\eta_{k+1}(p,q)\neq0$.
\item[$3)$]  The finite matrix of Hurwitz type $\mathcal{H}_{k}(p,q)$ defined in~\eqref{Hurwitz.matrix.finite.case.1}--\eqref{Hurwitz.matrix.finite.case.2} is nonsingular and totally
nonnegative. Here $k=2n$ if $\deg q<\deg p$, and $k=2n+1$ whenever
$\deg q=\deg p$.
\end{itemize}
\end{corol}

The total nonnegativity of the matrices of Hurwitz matrix also
implies the following curious result~\cite{Holtz_Tyaglov}.
\begin{theorem}\label{Th.interlacity.preserving}
Let the polynomials $p$ and $q$ defined
in~\eqref{polynomial.1}--\eqref{polynomial.2} have only
nonpositive zeroes, and let the function~$R=q/p$ be an
\textit{R}-function of negative type. Given any two positive
integers\footnote{The number $r$ is choosen such that
$b_0^2+b_r^2\neq0$.} $r$ and $l$ such that $rl\leqslant n<(l+1)r$,
the polynomials
\begin{eqnarray*}\label{Th.interlacity.preserving.poly.1}
p_{r,l}(z)&=&a_0z^l+a_{r}z^{l-1}+a_{2r}z^{l-2}+\ldots+a_{rl},\\
\label{Th.interlacity.preserving.poly.2}%
q_{r,l}(z)&=&b_0z^l+b_{r}z^{l-1}+b_{2r}z^{l-2}+\ldots+b_{rl}
\end{eqnarray*}
have only negative zeroes, and the
function~$R_{r,l}=q_{r,l}/p_{r,l}$ is an \textit{R}-function.
\end{theorem}

\subsection{Stieltjes continued fractions}\label{subsection:Stieltjes.cont.frac}

\hspace{4mm} Let a real rational function $F$ finite at infinity
be expanded into its Laurant series at infinity:
\begin{equation}\label{working.function}
F(z)=s_{-1}+\dfrac{s_{0}}{z}+\dfrac{s_{1}}{z^{2}}+\dfrac{s_{2}}{z^{3}}+\cdots
\end{equation}
\begin{definition}\label{def.Stieltjes.fraction}
The function $F$ is said to have Stieltjes continued fraction
expansion if $F$ can be represented in the form
\begin{equation}\label{Stieltjes.fraction.1}
F(z)=c_0+\dfrac1{c_1z+\cfrac1{c_2+\cfrac1{c_{3}z+\cfrac1{\ddots+\cfrac1{T}}}}},\quad
c_j\neq0,\quad\text{where}\quad
T=\begin{cases}
         &c_{2r}\qquad\ \,\text{if}\ |F(0)|<\infty,\\
         &c_{2r-1}z\quad\text{if}\ F(0)=\infty.
       \end{cases}
\end{equation}
Here $r$ is the number of poles of $F$, counting multiplicity.
\end{definition}

There exists the following criterion for a rational function to
have a Stieltjes continued fraction expansion~\cite{Wall} (see
also~\cite{Holtz_Tyaglov}).
\begin{theorem}\label{Th.Stieltjes.fraction.criterion}
Suppose that a rational function $F$ defined
in~\eqref{working.function} has exactly $r$ poles. The function
$F$ has a Stieltjes continued fraction
expansion~\eqref{Stieltjes.fraction.1} if and only if it satisfies
the conditions
\begin{equation}\label{minors.inequalities.for.Stieltjes.fraction.1}
D_{j}(F)\neq0,\quad j=1,2,\ldots,r,
\end{equation}
\begin{equation}\label{minors.inequalities.for.Stieltjes.fraction.2}
\widehat{D}_{j}(F)\neq0,\quad j=1,2,\ldots,r-1,
\end{equation}
\begin{equation}\label{minors.inequalities.for.Stieltjes.fraction.3}
D_{j}(F)=\widehat{D}_{j}(F)=0,\quad j=r+1,r+2,\ldots,
\end{equation}
where $D_j(F)$ and $\widehat{D}_j(F)$ are defined
in~\eqref{Hankel.determinants.1}
and~\eqref{Hankel.determinants.2}, respectively.

The coefficients of the continued
fraction~\eqref{Stieltjes.fraction.1} can be found by the
formul\ae
\begin{equation}\label{even.coeff.Stieltjes.fraction.main.formula}
c_{2j}=-\dfrac{D_{j}^2(F)}{\widehat{D}_{j-1}(F)\cdot\widehat{D}_{j}(F)},\quad
j=1,2,\ldots,r,
\end{equation}
\begin{equation}\label{odd.coeff.Stieltjes.fraction.main.formula}
c_{2j-1}=\dfrac{\widehat{D}_{j-1}^2(F)}{D_{j-1}(F)\cdot
D_{j}(F)},\quad j=1,2,\ldots,r,
\end{equation}
where $D_0(F)=\widehat{D}_0(F)\equiv1$.
\end{theorem}
We note that the determinant $\widehat{D}_r(F)$ in
Theorem~\ref{Th.Stieltjes.fraction.criterion} can be equal zero.
But in this case (and only in this case), the function~$F$ has a
pole at $0$ according to Corollary~\ref{corol.zero.pole}, and
therefore $T=c_{2r-1}z$ in~\eqref{Stieltjes.fraction.1}.

\vspace{2mm}

Let again the polynomials $p$ and $q$ be defined
in~\eqref{polynomial.1}--\eqref{polynomial.2}. Consider the
function $R=\dfrac{q}{p}$ and suppose that $R$ has a Stieltjes
continued fraction expansion~\eqref{Stieltjes.fraction.1}, where
$r\leqslant n=\deg p$. Then
from~\eqref{Hurwitz.determinants.relations.infinite.2.even}--\eqref{Hurwitz.determinants.relations.infinite.1.even}
and~\eqref{even.coeff.Stieltjes.fraction.main.formula}--\eqref{odd.coeff.Stieltjes.fraction.main.formula}
we obtain:
\begin{itemize}
\item [] if $\deg q(z)<\deg p(z)$, then
\begin{eqnarray}\label{odd.coeff.Stieltjes.fraction.main.formula.2.1}
c_{2j-1}&=&\dfrac{\eta_{2j-1}^2(p,q)}{\eta_{2j-2}(p,q)\cdot
\eta_{2j}(p,q)},\qquad j=1,2,\ldots,r;  \\
\label{even.coeff.Stieltjes.fraction.main.formula.2.1}
c_{2j}&=&\dfrac{\eta_{2j}^2(p,q)}{\eta_{2j-1}(p,q)\cdot\eta_{2j}(p,q)},\qquad
j=1,2,\ldots,\left[\dfrac{k}2\right];
\end{eqnarray}
\item [] if $\deg q(z)=\deg p(z)$, then
\begin{eqnarray}\label{odd.coeff.Stieltjes.fraction.main.formula.2.2}
c_{2j-1}&=&\dfrac{\eta_{2j}^2(p,q)}{\eta_{2j-1}(p,q)\cdot
\eta_{2j+1}(p,q)},\qquad j=1,2,\ldots,r; \\
\label{even.coeff.Stieltjes.fraction.main.formula.2.2}
c_{2j}&=&\dfrac{\eta_{2j+1}^2(p,q)}{\eta_{2j}(p,q)\cdot\eta_{2j+2}(p,q)},\qquad
\quad j=0,1,2,\ldots,\left[\dfrac{k}2\right];
\end{eqnarray}
\end{itemize}
where $k=2r-1$ if $R(0)=\infty$,  $k=2r$ if $|R(0)|<\infty$, and
$\eta_j(p,q)$ are the leading principal minors of the infinite
Hurwitz matrix $H(p,q)$ (see
Definition~\ref{def.Hurwitz.matrix.infinite}). Here we set
$\eta_0(p,q)\equiv1$, and $[\rho]$ denotes the largest integer not
exceeding $\rho$.

\vspace{3mm}

From
Theorems~\ref{Th.R-function.general.properties},~\ref{Th.number.of.negative.poles.of.R-frunction}
and~\ref{Th.Stieltjes.fraction.criterion} and from the
formul\ae~\eqref{even.coeff.Stieltjes.fraction.main.formula}--\eqref{odd.coeff.Stieltjes.fraction.main.formula}
one obtain the following theorem~\cite{Holtz_Tyaglov}.
\begin{theorem}\label{Th.R-functions.Stiljes.fractions}
Let the function $R$ with exactly $r$ poles, counting
multiplicities, be defined~\eqref{basic.rational.function}. If $R$
has a Stieltjes continued fraction
expansion~\eqref{Stieltjes.fraction.1}, then $R$ is an
\textit{R}-function if and only if
\begin{equation}\label{R-function.Stieltjes.fraction.condition}
c_{2j-1}>0,\quad j=1,2,\ldots,r.
\end{equation}
Moreover, the number of negative poles of the function $R$ equals
the number of positive coefficients $c_{2j}$, $j=1,2,\ldots,k$,
where $k=r$, if $|R(0)|<\infty$, and $k=r-1$, if $R(0)=\infty$.
\end{theorem}

Note that every \textit{R}-function with poles of the same sign
always has a Stieltjes continued fraction expansion.
\begin{corol}[A.\,Markov, Stieltjes, \cite{Markov,{Stieltjes1},{Stieltjes2}}]\label{corol.R-function.Stieltjes.fractions.nonpositive.poles}
A real rational function $R$ with exactly $r$ poles, counting
multiplicities, is an~\textit{R}-function with all nonpositive
poles if and only if $R$ has a Stieltjes continued fraction
expansion~\eqref{Stieltjes.fraction.1}, where
\begin{equation*}\label{R-function.Stieltjes.fraction.nonpositive.poles.condition}
c_{i}>0,\quad i=1,2,\ldots,2r-1.
\end{equation*}
Moreover, if $|R(0)|<\infty$, then $c_{2r}>0$.
\end{corol}
\begin{corol}[\cite{Holtz_Tyaglov}]\label{corol.R-function.Stieltjes.fractions.nonnegative.poles}
A real rational function $R$ with exactly $r$ poles, counting
multiplicities, is an~\textit{R}-function with all nonnegative
poles if and only if $R$ has a Stieltjes continued fraction
expansion~\eqref{Stieltjes.fraction.1}, where
\begin{equation*}\label{R-function.Stieltjes.fraction.nonnegative.poles.condition}
(-1)^{i-1}c_{i}>0,\quad i=1,2,\ldots,2r-1,
\end{equation*}
and if $|R(0)|<\infty$, then $c_{2r}<0$.
\end{corol}

\setcounter{equation}{0}

\section{Even and odd parts of polynomials. Associated
function}\label{section:assoc.func}

\hspace{4mm} Consider a real polynomial
\begin{equation}\label{main.polynomial}
p(z)=a_0z^n+a_1z^{n-1}+\dots+a_n,\qquad a_1,\dots,a_n\in\mathbb
R,\ a_0>0.
\end{equation}
In the rest of the paper we use the following notation
\begin{equation}\label{floor.poly.degree}
l=\left[\dfrac n2\right],
\end{equation}
where $n=\deg p$, and $[\rho]$ denotes the largest integer not
exceeding $\rho$.

The polynomial $p$ can always be represented as follows
\begin{equation}\label{app.poly.odd.even}
p(z)=p_0(z^2)+zp_1(z^2),
\end{equation}
where

for $n=2l$,
\begin{equation}\label{poly1.13}
\begin{split}
&p_0(u)=a_0u^l+a_2u^{l-1}+\ldots+a_{n},\\
&p_1(u)=a_1u^{l-1}+a_3u^{l-2}+\ldots+a_{n-1},
\end{split}
\end{equation}

and for $n=2l+1$,
\begin{equation}\label{poly1.12}
\begin{split}
&p_0(u)=a_1u^l+a_3u^{l-1}+\ldots+a_{n},\\
&p_1(u)=a_0u^l+a_2u^{l-1}+\ldots+a_{n-1}.
\end{split}
\end{equation}
The polynomials $p_0(z^2)$ and $p_1(z^2)$ satisfy the following
equalities:
\begin{equation}\label{poly1.2}
\begin{split}
&\displaystyle p_0(z^2)=\frac{p(z)+p(-z)}2,\\
&\displaystyle p_1(z^2)=\frac {p(z)-p(-z)}{2z}.
\end{split}
\end{equation}
Introduce the following function\footnote{In the
book~\cite[Chapter XV]{Gantmakher}, F.\,Gantmacher used the
function $-\dfrac{p_1(-u)}{p_0(-u)}$.}:
\begin{equation}\label{assoc.function}
\Phi(u)=\displaystyle\frac {p_1(u)}{p_0(u)}.
\end{equation}
\begin{definition}\label{def.associated.function}
We call $\Phi$ the \textit{function  associated with the
polynomial}~$p$.
\end{definition}
From~\eqref{poly1.2} and~\eqref{assoc.function} one can derive the
following relations:
\begin{equation}\label{poly1.1}
z\Phi(z^2)=\displaystyle\frac{p(z)-p(-z)}{p(z)+p(-z)}
=\displaystyle\frac{1-\displaystyle\frac{p(-z)}{p(z)}}{1+\displaystyle\frac{p(-z)}{p(z)}},
\end{equation}
\begin{equation*}
\frac{p(-z)}{p(z)}=\displaystyle\frac{1-z\Phi(z^2)}{1+z\Phi(z^2)}\,.
\end{equation*}

Now let us introduce two Hurwitz matrices associated with the
polynomial $p$.
\begin{definition}\label{def.Hurwitz.matrix.infinite.for.poly} The
matrix
\begin{equation}\label{Hurwitz.matrix.infinite.for.poly}
H_{\infty}(p)=
\begin{pmatrix}
a_0&a_2&a_4&a_6&a_8&a_{10}&\dots\\
0  &a_1&a_3&a_5&a_7&a_9&\dots\\
0  &a_0&a_2&a_4&a_6&a_8&\dots\\
0  &0  &a_1&a_3&a_5&a_7&\dots\\
\vdots &\vdots  &\vdots  &\vdots &\vdots  &\vdots &\ddots
\end{pmatrix}
\end{equation}
is called the \textit{infinite Hurwitz matrix} of the polynomial
$p$. The leading principal minors of the matrix $H_{\infty}(p)$
will be denoted by $\eta_j(p)$, $j=1,2,\ldots$
\end{definition}

\begin{remark}\label{remark.4.1}
According to Definitions~\ref{def.Hurwitz.matrix.infinite}
and~\ref{def.Hurwitz.matrix.infinite.for.poly}, we have
$H_{\infty}(p)=H(p_0,p_1)$, where $p_0$ and $p_1$ are the even and
odd parts of the polynomial $p$, respectively.
\end{remark}

Together with the infinite matrix $H_{\infty}(p)$, we consider its
specific finite submatrix.
\begin{definition}\label{def.Hurwitz.matrix.finite.for.poly} The
$n\times n$ matrix
\begin{equation}\label{HurwitzMatrix}
\mathcal{H}_n(p)=
\begin{pmatrix}
a_1&a_3&a_5&a_7&\dots&0&0\\
a_0&a_2&a_4&a_6&\dots&0&0\\
0  &a_1&a_3&a_5&\dots&0&0\\
0  &a_0&a_2&a_4&\dots&0&0\\
\vdots&\vdots&\vdots&\vdots&\ddots&\vdots&\vdots\\
0  &0  &0  &0  &\dots&a_{n-1} &0\\
0  &0  &0  &0  &\dots&a_{n-2} &a_n
\end{pmatrix}
\end{equation}
is called the \textit{finite Hurwitz matrix} or the
\textit{Hurwitz matrix}.
\end{definition}

\noindent The leading principal minors of this matrix we denote by
$\Delta_j(p)$:
\begin{equation}\label{delta}
\Delta_{j}(p)=
\begin{vmatrix}
a_1&a_3&a_5&a_7&\dots&a_{2j-1}\\
a_0&a_2&a_4&a_6&\dots&a_{2j-2}\\
0  &a_1&a_3&a_5&\dots&a_{2j-3}\\
0  &a_0&a_2&a_4&\dots&a_{2j-4}\\
\vdots&\vdots&\vdots&\vdots&\ddots&\vdots\\
0  &0  &0  &0  &\dots&a_{j}
\end{vmatrix},\quad j=1,\ldots,n,
\end{equation}
where we set $a_i\equiv0$ for $i>n$.
\begin{definition}\label{def.Hurwitz.dets}
The determinants $\Delta_{j}(p)$, $j=1,\ldots,n$, are called the
\textit{Hurwitz determinants} or the \textit{Hurwitz minors} of
the polynomial~$p$.
\end{definition}

\begin{remark}\label{remark.4.2}
From Definitions~\ref{def.Hurwitz.matrix.finite}
and~\ref{def.Hurwitz.matrix.finite.for.poly} it follows that
$H_{n}(p)=H_{2l}(p_0,p_1)$ if $n=2l$, and
$H_{n}(p)=H_{2l+1}(p_0,p_1)$ if $n=2l+1$, where $p_0$ and $p_1$
are defined in~\eqref{app.poly.odd.even}--\eqref{poly1.13}.
\end{remark}

\noindent Obviously,
\begin{equation}\label{Formulae.Gurwitz.3333}
\eta_j(p)=a_0\Delta_{j-1}(p),\quad j=1,\dots,n,
\end{equation}
where $\Delta_0(p)\equiv1$.

\vspace{2mm}

Suppose that $\deg p_0\geqslant\deg p_1$ and expand the function
$\Phi$ into its Laurent series at infinity:
\begin{equation}\label{app.assoc.function.series}
\Phi(u)=\dfrac{p_1(u)}{p_0(u)}=s_{-1}+\frac{s_0}u+\frac{s_1}{u^2}+\frac{s_2}{u^3}+\frac{s_3}{u^4}+\dots,
\end{equation}
where $s_{-1}\neq0$ if $\deg p_0=\deg p_1$, and $s_{-1}=0$ if
$\deg p_0>\deg p_1$.

From~\eqref{Hurwitz.determinants.relations.infinite.2.even}--\eqref{Hurwitz.finite.infnite.determinants.relation.2}
we obtain the following relations~\cite{Hurwitz,Gantmakher}
between the determinants $D_j(\Phi)$, $\widehat{D}_j(\Phi)$
defined in~\eqref{Hankel.determinants.1}
and~\eqref{Hankel.determinants.2}, the Hurwitz
minors~$\Delta_j(p)$, and the determinants $\eta_j(p)$
$(j=1,2,\dots,n)$.
\begin{itemize}
\item[1)]
If $n=2l$, then
\begin{equation}\label{Formulae.Gurwitz.1}
\begin{split}
&\Delta_{2j-1}(p)=a_0^{-1}\eta_{2j}(p)=a_0^{2j-1}D_j(\Phi),\\
&\Delta_{2j}(p)=a_0^{-1}\eta_{2j+1}(p)=(-1)^ja_0^{2j}\widehat{D}_j(\Phi),
\end{split}\qquad
j=1,2,\dots,l;
\end{equation}
\item[2)]
If $n=2l+1$, then
\begin{equation}\label{Formulae.Gurwitz.2}
\begin{split}
&\Delta_{2j}(p)=a_0^{-1}\eta_{2j+1}(p)={\left(\frac{a_0}{s_{-1}}\right)}^{2j}D_j(\Phi),\qquad j=1,2,\dots,l;\\
&\Delta_{2j+1}(p)=a_0^{-1}\eta_{2j+2}(p)=(-1)^j{\left(\frac{a_0}{s_{-1}}\right)}^{2j+1}\widehat{D}_j(\Phi),\qquad
j=0,1,\dots,l,
\end{split}
\end{equation}
where $\widehat{D}_0(\Phi)\equiv1$.
\end{itemize}

By~\eqref{app.poly.odd.even}--\eqref{poly1.12}
and~\eqref{assoc.function} and by
Theorem~\ref{Th.Hankel.matrix.rank.2} and
Corollary~\ref{corol.Hankel.matrix.rank.2}, we have
\begin{equation}\label{Theorem.Hurwitz.stable.poly.and.sign.regularity.condition.3}
D_j(\Phi)=\widehat{D}_j(\Phi)=0,\quad j>l,
\end{equation}
where $l$ is defined in~\eqref{floor.poly.degree}. So in the
sequel, we deal only with the determinants $D_j(\Phi)$,
$\widehat{D}_j(\Phi)$ of order at most $l$.

Also in the sequel, we deal only with the determinants $\eta_j(p)$
of order at most $n+1$, since by the
formul\ae~\eqref{Hurwitz.determinants.relations.infinite.2.even}--\eqref{Hurwitz.determinants.relations.infinite.1.even}
and~\eqref{Theorem.Hurwitz.stable.poly.and.sign.regularity.condition.3}
and by Remark~\ref{remark.4.1}, we have

\begin{equation}\label{Theorem.Hurwitz.stable.poly.and.sign.regularity.condition.44}
\eta_j(p)=0,\quad\text{for}\quad j>n+1.
\end{equation}

In Section~\ref{section:Generalized.Hurwitz} we also consider the
case when $n=2l+1$ with $a_1=0$ and $a_3\neq0$. In this case, we
have $\deg p_0=\deg p_1-1=l-1$, so the function $\Phi$ has the
form
\begin{equation}\label{app.assoc.function.series.Phi.1}
\Phi(u)=\dfrac{p_1(u)}{p_0(u)}=s_{-2}u+s_{-1}+\frac{s_0}u+\frac{s_1}{u^2}+\frac{s_2}{u^3}+\frac{s_3}{u^4}+\dots=\dfrac{a_0u^l+a_2u^{l-1}+\ldots+a_{2l}}{a_3u^{l-1}+a_5u^{l-2}+\ldots+a_{2l+1}}.
\end{equation}
It is easy to see that $t_{-2}=\dfrac{a_0}{a_3}$, so we can
represent $\Phi$ as a sum
$\Phi(u)=\dfrac{a_0}{a_3}\,u+\dfrac{p_2(u)}{p_0(u)}$, where
\begin{equation*}\label{poly.p2}
p_2(u)=\left(a_2-\dfrac{a_0}{a_3}\,a_5\right)u^{l-1}+\left(a_4-\dfrac{a_0}{a_3}\,a_7\right)u^{l-2}+\ldots+\left(a_{2l-2}-\dfrac{a_0}{a_3}\,a_{2l+1}\right)u+a_{2l}.
\end{equation*}
Note that for the function
\begin{equation*}\label{app.assoc.function.series.Phi.2}
\Theta(u):=\dfrac{p_2(u)}{p_0(u)}=\Phi(u)-\dfrac{a_0}{a_3}\,u=s_{-1}+\frac{s_0}u+\frac{s_1}{u^2}+\frac{s_2}{u^3}+\frac{s_3}{u^4}+\dots,
\end{equation*}
the following equalities hold
\begin{equation}\label{app.assoc.function.series.Phi.3}
D_j(\Theta)=D_j(\Phi),\quad\text{and}\quad\widehat{D}_j(\Theta)=\widehat{D}_j(\Phi)\quad\text{for}\quad
j=1,2,\ldots
\end{equation}
By Definition~\ref{def.Hurwitz.matrix.finite} and by the
formul\ae~\eqref{Hurwitz.determinants.relations.infinite.1.odd},
and~\eqref{Hurwitz.finite.infnite.determinants.relation.2} applied
to the polynomials $p_0$ and $p_2$, we have, for $j=1,2,\ldots$,
\begin{equation*}
\begin{array}{c}
-a_0a_3^{2j+1}D_j(\Theta)=-a_0a_3\nabla_{2j}(p_0,p_2)=-a_0a_3
\begin{vmatrix}
a_3&a_5&\dots&a_{4j+1}\\
a_2-\dfrac{a_0}{a_3}\,a_5& a_4-\dfrac{a_0}{a_3}\,a_7&\dots& a_{4j}-\dfrac{a_0}{a_3}\,a_{4j+3}\\
0   &a_3&\dots&a_{4j-1}\\
0 & a_2-\dfrac{a_0}{a_3}\,a_5&\dots& a_{4j-2}-\dfrac{a_0}{a_3}\,a_{4j+1}\\
\vdots&\vdots&\ddots&\vdots\\
0  &0  &\dots&a_{2j+2}-\dfrac{a_0}{a_3}\,a_{2j+5}
\end{vmatrix}=\\
=
\begin{vmatrix}
0&a_3&a_5&a_7&\dots&a_{4j+3}\\
a_0&a_2&a_4&a_6&\dots&a_{4j+2}\\
0&0&a_3&a_5&\dots&a_{4j+1}\\
0&0&a_2-\dfrac{a_0}{a_3}\,a_5& a_4-\dfrac{a_0}{a_3}\,a_7&\dots& a_{4j}-\dfrac{a_0}{a_3}\,a_{4j+3}\\
0&0&0   &a_3&\dots&a_{4j-1}\\
0&0&0 & a_2-\dfrac{a_0}{a_3}\,a_5&\dots& a_{4j-2}-\dfrac{a_0}{a_3}\,a_{4j+1}\\
\vdots&\vdots&\vdots&\vdots&\ddots&\vdots\\
0&0&0  &0  &\dots&a_{2j+2}-\dfrac{a_0}{a_3}\,a_{2j+5}
\end{vmatrix}=\\
=
\begin{vmatrix}
0&a_3&a_5&a_7&\dots&a_{4j+3}\\
a_0&a_2&a_4&a_6&\dots&a_{4j+2}\\
0&0&a_3&a_5&\dots&a_{4j+1}\\
0&a_0&a_2& a_4&\dots& a_{4j}\\
0&0&0   &a_3&\dots&a_{4j-1}\\
0&0&a_0 & a_2&\dots& a_{4j-2}\\
\vdots&\vdots&\vdots&\vdots&\ddots&\vdots\\
0&0&0  &0  &\dots&a_{2j+2}
\end{vmatrix}=\Delta_{2j+2}(p)=-a_0a_3^{2j+1}D_j(\Phi),
\end{array}
\end{equation*}
where we set $a_i\equiv0$ for $i>2l+1$.

Denote the coefficients of the polynomial $p_2$ by $b_j$:
\begin{equation}\label{Formulae.Hurwitz.1.Psi}
b_j=a_{2j+2}-\dfrac{a_0}{a_3}a_{2j+5},\qquad j=0,1,\ldots,l-1,
\end{equation}
where $a_{2l+3}\equiv0$. Then the function $u\Phi(u)$ has the form
\begin{equation*}
u\Phi(u)=\dfrac{a_0}{a_3}\,u^2+\dfrac{b_0}{a_3}\,u+\dfrac{p_3(u)}{p_0(u)},
\end{equation*}
where
\begin{equation*}
p_3(u):=up_2(u)-\dfrac{b_0}{a_3}\,up_0(u)=\left(b_1-\dfrac{b_0}{a_3}\,a_5\right)u^{l-1}+\left(b_2-\dfrac{b_0}{a_3}\,a_7\right)u^{l-2}+\ldots+\left(b_{l-1}-\dfrac{b_0}{a_3}\,a_{2l+1}\right)u.
\end{equation*}

\noindent By Definition~\ref{def.Hurwitz.matrix.finite} and by the
formul\ae~\eqref{app.assoc.function.series.Phi.3},~\eqref{Hurwitz.determinants.relations.infinite.1.odd}
and~\eqref{Hurwitz.finite.infnite.determinants.relation.2} applied
to the polynomials $p_0$ and $p_3$ we get (putting $b_i\equiv0$
for $i>l-1$):
\begin{equation*}
\begin{array}{c}
a_0a_3^{2j+2}\widehat{D}_j(\Theta)=a_0a_3^2\nabla(p_0,p_3)=a_0a_3^2
\begin{vmatrix}
a_3&a_5&\dots&a_{4j+1}\\
b_1-\dfrac{b_0}{a_3}\,a_5& b_2-\dfrac{b_0}{a_3}\,a_7&\dots& b_{2j}-\dfrac{b_0}{a_3}\,a_{4j+3}\\
0   &a_3&\dots&a_{4j-1}\\
0 & b_1-\dfrac{b_0}{a_3}\,a_5&\dots& b_{2j-1}-\dfrac{b_0}{a_3}\,a_{4j+1}\\
\vdots&\vdots&\ddots&\vdots\\
0  &0  &\dots&b_{j+1}-\dfrac{b_0}{a_3}\,a_{2j+5}
\end{vmatrix}=\\
=a_0a_3
\begin{vmatrix}
a_3&a_5&a_7&\dots&a_{4j+3}\\
0&a_3&a_5&\dots&a_{4j+1}\\
0&b_1-\dfrac{b_0}{a_3}\,a_5& b_2-\dfrac{b_0}{a_3}\,a_7&\dots& b_{2j}-\dfrac{b_0}{a_3}\,a_{4j+3}\\
0&0   &a_3&\dots&a_{4j-1}\\
0&0 & b_1-\dfrac{b_0}{a_3}\,a_5&\dots& b_{2j-1}-\dfrac{b_0}{a_3}\,a_{4j+1}\\
\vdots&\vdots&\vdots&\ddots&\vdots\\
0&0  &0  &\dots&b_{j+1}-\dfrac{b_0}{a_3}\,a_{2j+5}
\end{vmatrix}
=(-1)^j a_0a_3
\begin{vmatrix}
a_3&a_5&a_7&\dots&a_{4j+3}\\
b_0&b_1&b_2&\dots& b_{2j}\\
0&a_3&a_5&\dots&a_{4j+1}\\
0&b_0&b_1&\dots&b_{2j-1}\\
0&0   &a_3&\dots&a_{4j-1}\\
0&0&b_0&\dots&b_{2j-2}\\
\vdots&\vdots&\vdots&\ddots&\vdots\\
0&0  &0  &\dots&a_{2j+3}
\end{vmatrix}=\\
=(-1)^{j+1}
\begin{vmatrix}
0&a_3&a_5&a_7&a_9&\dots&a_{4j+5}\\
a_0&a_2&a_4&a_6&a_8&\dots&a_{4j+4}\\
0&0&a_3&a_5&a_7&\dots&a_{4j+3}\\
0&0&a_2-\dfrac{a_0}{a_3}\,a_5& a_4-\dfrac{a_0}{a_3}\,a_7& a_6-\dfrac{a_0}{a_3}\,a_9&\dots& a_{4j+2}-\dfrac{a_0}{a_3}\,a_{4j+5}\\
0&0&0&a_3&a_5&\dots&a_{4j+1}\\
0&0&0&a_2-\dfrac{a_0}{a_3}\,a_5& a_4-\dfrac{a_0}{a_3}\,a_7&\dots& a_{4j}-\dfrac{a_0}{a_3}\,a_{4j+3}\\
0&0&0&0   &a_3&\dots&a_{4j-1}\\
0&0&0&0&a_2-\dfrac{a_0}{a_3}\,a_5&\dots& a_{4j-2}-\dfrac{a_0}{a_3}\,a_{4j+1}\\
\vdots&\vdots&\vdots&\vdots&\vdots&\ddots&\vdots\\
0&0&0&0  &0  &\dots&a_{2j+3}
\end{vmatrix}=\\
 \\
=(-1)^{j+1}\Delta_{2j+3}(p)=a_0a_3^{2j+2}\widehat{D}_j(\Phi),
\end{array}
\end{equation*}
where we set $a_i\equiv0$ for $i>2l+1$. Here we used the
equalities $\widehat{D}_j(\Theta)=D_j(u\Theta)$, $j=1,2,\ldots$

Thus, we obtain the following formul\ae
\begin{equation}\label{Formulae.Gurwitz.3}
\begin{split}
&\Delta_{2j}(p)=a_0^{-1}\eta_{2j+1}(p)=-a_0a_3^{2j-1}D_{j-1}(\Phi)=a_3^{2j-2}\Delta_2(p)D_{j-1}(\Phi),\quad\qquad\qquad\qquad j=2,3,\dots,l,\\
&\Delta_{2j+1}(p)=a_0^{-1}\eta_{2j+2}(p)=(-1)^{j}a_0a_3^{2j}\widehat{D}_{j-1}(\Phi)=(-1)^{j-1}a_3^{2j-2}\Delta_3(p)\widehat{D}_{j-1}(\Phi),\qquad
j=2,3,\dots,l,
\end{split}
\end{equation}
where $\eta_{i}(p)$ are the leading principal minors of the matrix
$H_{\infty}(p)$ defined
in~\eqref{Hurwitz.matrix.infinite.for.poly}. Here we used the
formul\ae
\begin{equation}\label{Formulae.Gurwitz.333}
\Delta_2(p)=a_0^{-1}\eta_{3}(p)=-a_0a_3,\quad
\Delta_3(p)=a_0^{-1}\eta_{4}(p)=-a_0a_3^2,
\end{equation}
that follow from the equality $a_1=0$.

Note that, as above, the following holds:
\begin{equation*}
D_j(\Phi)=\widehat{D}_j(\Phi)=0,\quad j>l-1.
\end{equation*}

\setcounter{equation}{0}

\section{Stable polynomials}\label{section:Hurwitz.polys}

\hspace{4mm} This section is devoted to some basic facts of the
theory of Hurwitz stable and quasi-stable polynomials.~We expound
those facts from the viewpoint of the theory of
\textit{R}-functions.
%

\subsection{Hurwitz polynomials}\label{subsection:Hurwitz.poly}

\begin{definition}[\cite{Hurwitz,{Gantmakher}}]\label{def.Hurwitz.stable.poly}
The polynomial $p$ defined in~\eqref{main.polynomial} is called
\textit{Hurwitz} or \textit{Hurwitz stable} if all its zeroes lie
in the \textit{open} left half-plane of the complex plane.
\end{definition}

At first, we recall the following simple necessary condition for
polynomials to be Hurwitz stable. This condition is called usually
\textit{Stodola theorem}~\cite{Gantmakher,{Barkovsky.2}}.

\begin{theorem}[Stodola]\label{Th.Stodola.necessary.condition.Hurwitz}
If the polynomial $p$ is Hurwitz stable, then all its coefficients
are positive\footnote{More exactly, the coefficients must be of
the same sign, but $a_0>0$ by~\eqref{main.polynomial}.}.
\end{theorem}

Now we show a connection between Hurwitz polynomials and
\textit{R}-functions. This connection will allow us to apply all
statements from the Section~\ref{section:auxiliary.theorems} to
Hurwitz polynomials to obtain such basic criteria of Hurwitz
stability as Hurwitz criterion and Li\'enard and Chipart criterion
(Theorems~\eqref{Theorem.Hurwitz.stable.Hurwitz.matrix.criteria}
and~\ref{Th.Lienard.Chipart}, respectively). We prove the
following theorem using a method suggested by
Yu.S.\,Barkovsky~\cite{Barkovsky.pr}.

\begin{theorem}[\cite{Gantmakher,{Barkovsky.2}}]\label{Theorem.main.Hurwitz.stability}
A real polynomial $p$ defined in~\eqref{main.polynomial} is
Hurwitz stable if and only if its associated function~$\Phi$
defined in~\eqref{assoc.function} is an \textit{R}-function with
exactly $l$ poles, all of which  are negative, and the limit
$\displaystyle\lim_{u\to\pm\infty}\Phi(u)$ is positive whenever
$n=2l+1$. The number $l$ is defined in~\eqref{floor.poly.degree}.
\end{theorem}
\begin{proof}
Let the polynomial $p$ be Hurwitz stable, that is,
\begin{equation}\label{Theorem.main.Hurwitz.stability.proof.1}
p(\lambda)=0\ \ \Longrightarrow\ \ \Re\lambda<0.
\end{equation}

First, we show that
\begin{equation}\label{Theorem.main.Hurwitz.stability.proof.2.2}
\displaystyle\left|\frac {p(-z)}{p(z)}\right|<1,\quad\forall z:\,
\Re{z}>0.
\end{equation}
Note that the polynomials $p(z)$ and $p(-z)$ have no common zeroes
if $p$ is Hurwitz stable, so the function $\dfrac {p(-z)}{p(z)}$
has exactly $n$ poles. The Hurwitz stable polynomial $p$ can be
represented in the form
\begin{equation*}\label{Theorem.main.Hurwitz.stability.proof.2.5}
p(z)=a_0\prod_k(z-\lambda_k)\prod_j\left(z-\xi_j\right)\left(z-\overline{\xi}_j\right),
\end{equation*}
where $\lambda_k<0,\Re\xi_j<0$ and $\Im\xi_j\neq0$. Then we have
\begin{equation}\label{poly1.5}
\displaystyle\left|\frac {p(-z)}{p(z)}\right|=\prod_{k}\frac
{\left|z+\lambda_k\right|}{\left|z-\lambda_k\right|}\prod_{j}\frac
{\left|z+\xi_j\right|\left|z+\overline{\xi}_j\right|}{\left|z-\overline{\xi}_j\right|\left|z-\xi_j\right|}.
\end{equation}
It is easy to see that the function of type
$\dfrac{z+a}{z-\overline{a}}$~, where $\Re a<0$, maps the right
half-plane of the complex plane to the~unit disk. In fact,
\begin{equation*}\label{Theorem.main.Hurwitz.stability.proof.3}
\left|\dfrac{z+a}{z-\overline{a}}\right|^2=\dfrac{(\Re z+\Re
a)^2+(\Im z+\Im a)^2}{(\Re z-\Re a)^2+(\Im z+\Im
a)^2}<1\qquad\text{whenever}\quad\Re z>0\quad\text{and}\quad\Re
a<0.
\end{equation*}
Now from~\eqref{poly1.5} it follows that the function
$\displaystyle\frac {p(-z)}{p(z)}$ also maps the right half-plane
to the~unit disk as a~product of functions of such a type. Thus,
the inequality~\eqref{Theorem.main.Hurwitz.stability.proof.2.2} is
valid.

At the same time, the fractional linear transformation
$\displaystyle z\mapsto\frac{1-z}{1+z}$ conformally maps the unit
disk to the right half-plane:
\begin{equation}\label{Theorem.main.Hurwitz.stability.proof.5}
|z|<1\implies\Re\left(\frac{1-z}{1+z}\right)=\dfrac{1-|z|^2}{|1+z|^2}>0.
\end{equation}
Consequently, from the
relations~\eqref{poly1.1},~\eqref{Theorem.main.Hurwitz.stability.proof.2.2}
and~\eqref{Theorem.main.Hurwitz.stability.proof.5} we obtain that
the function $z\Phi(z^2)$ maps the~right half-plane to itself, so
the function $-z\Phi(-z^2)$ maps the upper half-plane of the
complex plane to the lower half-plane:
\begin{equation*}\label{Theorem.main.Hurwitz.stability.proof.6}
\Im z>0\implies\Re(-iz)>0\implies
\Re\left[-iz\Phi(-z^2)\right]=\Im\left[z\Phi(-z^2)\right]>0\implies\Im\left[-z\Phi(-z^2)\right]<0.
\end{equation*}
Since $p$ is Hurwitz stable by assumption, the polynomials $p(z)$
and $p(-z)$ have no common zeroes, therefore, $p_0$ and $p_1$ also
have no common zeroes, and $p_0(0)\neq0$
by~\eqref{app.poly.odd.even}. Moreover, by
Theorem~\ref{Th.Stodola.necessary.condition.Hurwitz} we have
$a_0>0$ and $a_1>0$, so $\deg p_0=l$
(see~\eqref{poly1.13}--\eqref{poly1.12}). Thus, the number of
poles of the function $-z\Phi(-z^2)$ equals the number of zeroes
of the polynomial $p_0(-z^2)$, i.e. exactly $2l$.

So according to Theorem~\ref{Th.R-function.general.properties},
the function $-z\Phi(-z^2)$ can be represented in the
form~\eqref{Mittag.Leffler.1}, where all poles are located
symmetrically with respect to~$0$ and $\beta=0$, since
$-z\Phi(-z^2)$ is an odd function. Denote the poles of
$-z\Phi(-z^2)$ by $\pm\nu_1,\ldots,\pm\nu_l$ such that
\begin{equation*}\label{Th.stable.poly.and.R-function.proof.6.5}
0<\nu_1<\nu_2<\ldots<\nu_l.
\end{equation*}
Note that $\nu_1\neq0$, since $p_0(0)\neq0$ as we mentioned above.

Thus, the function $-z\Phi(-z^2)$ can be represented in the
following form
\begin{equation*}\label{Th.stable.poly.and.R-function.proof.7}
-z\Phi(-z^2)=-\alpha z+\sum_{j=1}^l\frac{\gamma_j}{z-\nu_j}+
\sum_{j=1}^l\frac{\gamma_j}{z+\nu_j}= -\alpha z
+\sum_{j=1}^l\frac{2\gamma_jz}{z^2-\nu_j^2},
\quad\alpha\geqslant0,\,\gamma_j,\nu_j>0.
\end{equation*}
Dividing this equality by $-z$ and changing variables as follows
$-z^2\to u$, $2\gamma_j\to\beta_j$, $\nu_j^2\to\omega_j$, we
obtain the following representation of the function $\Phi$:
\begin{equation}\label{Mittag.Leffler}
\Phi(u)=\dfrac{p_1(u)}{p_0(u)}=\alpha+\sum_{j=1}^l\frac{\beta_j}{u+\omega_j},
\end{equation}
where $\alpha\geqslant0,\,\beta_j>0$ and
\begin{equation*}\label{Th.stable.poly.and.R-function.proof.8}
0<\omega_1<\omega_2<\ldots<\omega_l.
\end{equation*}
Here $\alpha=0$ whenever $n=2l$, and
$\alpha=\displaystyle\frac{a_0}{a_1}>0$ whenever $n=2l+1$. Since
$\Phi$ can be represented in the form~\eqref{Mittag.Leffler}, we
have that by Theorem~\ref{Th.R-function.general.properties},
$\Phi$ is an \textit{R}-function with exactly $l$ poles, all which
are negative, and
$\displaystyle\lim_{u\to\pm\infty}\Phi(u)=\alpha>0$ as $n=2l+1$.

\vspace{3mm}

Conversely, let the polynomial $p$ be defined
in~\eqref{main.polynomial} and let its associated function $\Phi$
be an \textit{R}-function with exactly $l$ poles, all of which are
negative, and $\displaystyle\lim_{u\to\pm\infty}\Phi(u)>0$ as
$n=2l+1$. We will show that $p$ is Hurwitz stable.

By Theorem~\ref{Th.R-function.general.properties}, $\Phi$ can be
represented in the form~\eqref{Mittag.Leffler}, where
$\alpha=\displaystyle\lim_{u\to\pm\infty}\Phi(u)\geqslant0$ such
that $\alpha>0$ if $n=2l+1$, and $\alpha=0$ if $n=2l$. Thus, the
polynomial $p_0$ has only negative zeroes, and the polynomials
$p_0$ and $p_1$ have no common zeroes. Together
with~\eqref{app.poly.odd.even} and~\eqref{assoc.function}, this
implies that the set of zeroes of the polynomial $p$ coincides
with the set of roots of the equation
\begin{equation}\label{poly1.9}
z\Phi(z^2)=-1.
\end{equation}
Let $\lambda$ be a zero of the polynomial $p$ and therefore, of
the equation~\eqref{poly1.9}. Then from~\eqref{Mittag.Leffler}
and~\eqref{poly1.9} we obtain
\begin{equation*}
-1=\Re\left[\lambda\Phi(\lambda^2)\right]=\left[\alpha+
\sum_{j=1}^l\beta_j\frac{|\lambda|^2+\omega_j}{|\lambda^2+\omega_j|^2}
\right]\Re\lambda,
\end{equation*}
where $\alpha\geqslant0$, and $\beta_j,\,\omega_j>0$ for
$j=1,\ldots,l$. Thus, if $\lambda$ is a zero of $p$, then
$\Re\lambda<0$, so $p$ is Hurwitz stable.
\end{proof}

Consider the Hankel matrix $S(\Phi)=\|s_{i+j}\|^{\infty}_{i,j=0}$
constructed with the coefficients of the
series~\eqref{app.assoc.function.series} and its determinants
$D_j(\Phi)$ and $\widehat{D}_j(\Phi)$ defined
in~\eqref{Hankel.determinants.1}
and~\eqref{Hankel.determinants.2}, respectively. From
Theorems~\ref{Th.R-function.general.properties},~\ref{Theorem.main.Hurwitz.stability}
and~\ref{Theorem.sign.regularity.Hankel.matrix} and from
Corollary~\ref{Corol.R-functions.neg.poles} we obtain the
following Hurwitz stability criteria:
\begin{theorem}\label{Theorem.Hurwitz.stable.poly.and.sign.regularity}
Let a real polynomial $p$ be defined by~\eqref{main.polynomial}.
The following conditions are equivalent:
\begin{itemize}
\item[1)] the polynomial $p$ is Hurwitz stable;
\item[2)] the following hold
\begin{equation}\label{Theorem.Hurwitz.stable.poly.and.sign.regularity.condition.0}
\begin{split}
&s_{-1}>0\quad\text{for}\quad n=2l+1,\\
&D_j(\Phi)>0, \qquad j=1,\ldots,l,\\
&(-1)^{j}\widehat{D}_{j}(\Phi)>0, \qquad j=1,\ldots,l,
\end{split}
\end{equation}
where $l=\left[\dfrac n2\right]$;
\item[3)] the matrix $S(\Phi)$ is sign regular of rank $l$, and $s_{-1}>0$ for $n=2l+1$.
\end{itemize}
\end{theorem}
\begin{proof}
In fact, by Theorem~\ref{Theorem.main.Hurwitz.stability}, the
polynomial $p$ is Hurwitz stable if and only if its associated
function~$\Phi$ is an \textit{R}-function with exactly $l$ poles,
all of which are negative and
$\lim\limits_{u\to\pm\infty}\Phi(u)=s_{-1}\geqslant0$. According
to Theorem~\ref{Th.R-function.general.properties} and
Corollary~\ref{Corol.R-functions.neg.poles}, this is equivalent to
the
inequalities~\eqref{Theorem.Hurwitz.stable.poly.and.sign.regularity.condition.0}.
But these inequalities also are equivalent to the sign regularity
of the matrix $S(\Phi)$ by
Theorem~\ref{Theorem.sign.regularity.Hankel.matrix}.
\end{proof}

\noindent Our next theorem provides stability criteria in terms of
coefficients of the polynomial $p$.

\begin{theorem}\label{Theorem.Hurwitz.stable.Hurwitz.matrix.criteria}
Given a real polynomial $p$ of degree $n$ as
in~\eqref{main.polynomial}, the following conditions are
equivalent:
\begin{itemize}
\item[1)] the polynomial $p$ is Hurwitz stable;
\item[2)] all Hurwitz determinants $\Delta_j(p)$ are positive:
\begin{equation}\label{Hurvitz.det.noneq}
\Delta_1(p)>0,\ \Delta_2(p)>0,\dots,\ \Delta_n(p)>0;
\end{equation}
\item[3)] the determinants $\eta_j(p)$ are positive up to order $n+1$:
\begin{equation}\label{Hurvitz.det.noneq.infinite}
\eta_1(p)>0,\ \eta_2(p)>0,\dots,\ \eta_{n+1}(p)>0;
\end{equation}
\item[4)] the matrix $\mathcal{H}_n(p)$ defined in~\eqref{HurwitzMatrix} is nonsingular and totally
nonnegative;
\item[5)] the matrix $H_{\infty}(p)$ defined in~\eqref{Hurwitz.matrix.infinite.for.poly} is totally nonnegative with nonzero minor $\eta_{n+1}(p)$.
\end{itemize}
\end{theorem}
\noindent Note that the equivalence of $1)$ and $2)$ is the famous
Hurwitz criterion of stability. The implications
$1)\Longrightarrow4)$ and $1)\Longrightarrow5)$ were proved
in~\cite{Asner,{Kemperman}}. The implication $4)\Longrightarrow1)$
was, in fact, proved in~\cite{Asner}. However, the implication
$5)\Longrightarrow1)$ is probably new.

\begin{proof}
Indeed, by
Theorem~\ref{Theorem.Hurwitz.stable.poly.and.sign.regularity}, the
polynomial $p$ is Hurwitz stable if and only if the
inequalities~\eqref{Theorem.Hurwitz.stable.poly.and.sign.regularity.condition.0}
hold. According to the
formul\ae~\eqref{Formulae.Gurwitz.1}--\eqref{Formulae.Gurwitz.2},
these inequalities are equivalent to~\eqref{Hurvitz.det.noneq}.
By~\eqref{Formulae.Gurwitz.3333}, the
inequalities~\eqref{Hurvitz.det.noneq} are equivalent
to~\eqref{Hurvitz.det.noneq.infinite}, since $\eta_1(p)=a_0>0$
(see~\eqref{Hurwitz.matrix.infinite.case.1}--\eqref{Hurwitz.matrix.infinite.case.2}
and Remark~\ref{remark.4.1}).

Furthermore, by Theorem~\ref{Theorem.main.Hurwitz.stability}, the
polynomial $p$ is Hurwitz stable if and only if its associated
function~$\Phi$ is an \textit{R}-function with exactly
$l=\left[\dfrac n2\right]$, all of which are negative, and
$\lim\limits_{u\to\pm\infty}\Phi(u)=\dfrac{a_0}{a_1}=s_{-1}>0$ if
$n=2l+1$. According to
Corollary~\ref{Corol.Hurwitz.Matrix.Total.Nonnegativity.1and2} 
and Remark~\ref{remark.4.1}, this is equivalent to the total
nonnegativity of the~mat\-rix~$H_{\infty}(p)$. Moreover,
$\eta_{n+1}(p)\neq0$ by
Corollary~\ref{Corol.Hurwitz.Matrix.Total.Nonnegativity.1and2}.
Thus, the condition $1)$ is equivalent to the condition~$5)$.

The conditions $1)$ and $4)$ are equivalent also by
Corollary~\ref{Corol.Hurwitz.Matrix.Total.Nonnegativity.1and2} and
Remark~\ref{remark.4.2}.
\end{proof}

\begin{remark}
Note that total nonnegativity of the Hurwitz matrix of a Hurwitz
polynomial was established in~\cite{Asner,{Kemperman},{Holtz1}}.
\end{remark}

Now we are in a position to prove a few another famous criteria of
Hurwitz stability, which are known as the Li\'enard and Chipart
criterion and its modifications
(see~\cite{LienardChipart,{Gantmakher}}). These criteria are
simple consequences of
Theorems~\ref{Th.general.Lienard.Chipart.negative.case.1},~\ref{Th.general.Lienard.Chipart.negative.case.2}
and~\ref{Theorem.main.Hurwitz.stability}.
\begin{theorem}\label{Th.Lienard.Chipart}
The polynomial $p$ given by~\eqref{main.polynomial} is Hurwitz
stable if and only if one of the following conditions holds
\begin{itemize}
\item[$1)$]
$a_n>0,a_{n-2}>0,a_{n-4}>0,\ldots,\qquad\Delta_{n-1}(p)>0,\Delta_{n-3}(p)>0,\Delta_{n-5}(p)>0,\ldots$;
\item[$2)$]
$a_n>0,a_{n-1}>0,a_{n-3}>0,\ldots,\qquad\Delta_{n-1}(p)>0,\Delta_{n-3}(p)>0,\Delta_{n-5}(p)>0,\ldots$;
\item[$3)$]
$a_n>0,a_{n-2}>0,a_{n-4}>0,\ldots,\qquad\Delta_n(p)>0,\Delta_{n-2}(p)>0,\Delta_{n-4}(p)>0,\ldots$;
\item[$4)$]
$a_n>0,a_{n-1}>0,a_{n-3}>0,\ldots,\qquad\Delta_n(p)>0,\Delta_{n-2}(p)>0,\Delta_{n-4}(p)>0,\ldots$
\end{itemize}
\end{theorem}
\begin{proof}
Let $n=2l$. By Definitions~\ref{def.Hurwitz.matrix.finite}
and~\ref{def.Hurwitz.matrix.finite.for.poly} and by
Remark~\ref{remark.4.2}, we have $H_n(p)=H_{2l}(p_0,p_1)$, where
the polynomials $p_0$ and $p_1$ are defined
in~\eqref{app.poly.odd.even}--\eqref{poly1.13}. Obviously,
$\Delta_j(p)=\Delta_j(p_0,p_1)$, $j=1,\ldots,n$. The assertion of
the theorem follows now from
Theorems~\ref{Theorem.main.Hurwitz.stability}
and~\ref{Th.general.Lienard.Chipart.negative.case.1}.

In the case $n=2l+1$, the assertion of the theorem can be proved
in the same way using
Theorem~\ref{Th.general.Lienard.Chipart.negative.case.2} instead
of Theorem~\ref{Th.general.Lienard.Chipart.negative.case.1}.
\end{proof}

We now recall a connection between Hurwitz polynomials and
Stieltjes continued fractions (see~\cite{Gantmakher}).
\begin{theorem}\label{Theorem.Hurwitz.stable.Stieltjes.cont.frac.criteria}
The polynomial $p$ of degree $n\geqslant1$ defined
in~\eqref{main.polynomial} is Hurwitz stable if and only if its
associated function~$\Phi$ has the following Stieltjes continued
fraction expansion:
\begin{equation}\label{Stieltjes.fraction.for.stability}
\Phi(u)=c_0+\dfrac1{c_1u+\cfrac1{c_2+\cfrac1{c_{3}u+\cfrac1{\ddots+\cfrac1{c_{2l-1}u+\cfrac1{c_{2l}}}}}}},\quad\text{with}\quad
c_i>0,\quad i=1,\ldots,2l,
\end{equation}
where $c_0=0$ if $n$ is even, $c_0>0$ if $n$ is odd, and $l$ as
in~\eqref{floor.poly.degree}.
\end{theorem}
\begin{proof}
In fact, by Theorem~\ref{Theorem.main.Hurwitz.stability}, the
polynomial $p$ is Hurwitz stable if and only if the function
$\Phi$ can be represented in the form~\eqref{Mittag.Leffler},
where $\alpha=c_0=0$ if $n=2l$, and
$\alpha=c_0=\dfrac{a_0}{a_1}>0$ if $n=2l+1$. Now the~assertion of
the theorem follows from
Theorem~\ref{Th.R-function.general.properties} and
Corollary~\ref{corol.R-function.Stieltjes.fractions.nonpositive.poles}.
\end{proof}

From the
formul\ae~\eqref{even.coeff.Stieltjes.fraction.main.formula}--\eqref{odd.coeff.Stieltjes.fraction.main.formula}
it follows that the coefficients $c_i$
of~\eqref{Stieltjes.fraction.for.stability} can be found by the
following formul\ae
\begin{equation}\label{even.coeff.Stieltjes.fraction.main.formula.quasi-stability}
c_{2j}=-\dfrac{D_{j}^2(\Phi)}{\widehat{D}_{j-1}(\Phi)\cdot\widehat{D}_{j}(\Phi)},\quad
j=1,2,\ldots,l,
\end{equation}
\begin{equation}\label{odd.coeff.Stieltjes.fraction.main.formula.quasi-stability}
c_{2j-1}=\dfrac{\widehat{D}_{j-1}^2(\Phi)}{D_{j-1}(\Phi)\cdot
D_{j}(\Phi)},\quad j=1,2,\ldots,l,
\end{equation}
where $D_0(\Phi)=\widehat{D}_0(\Phi)\equiv1$. If $n$ is odd, then
$c_0=\dfrac{a_0}{a_1}>0$.

Using the
formul\ae~\eqref{even.coeff.Stieltjes.fraction.main.formula.quasi-stability}--\eqref{odd.coeff.Stieltjes.fraction.main.formula.quasi-stability}
and~\eqref{Formulae.Gurwitz.1}--\eqref{Formulae.Gurwitz.2}, we can
represent the coefficients  $c_i$ in terms of Hurwitz
determinants~$\Delta_i(p)$:
\begin{itemize}
\item[1)]
If $n=2l$, then
\begin{equation}\label{coeff.Stieltjes.fraction.main.formula.quasi-stability.even.degree}
c_{i}=\dfrac{\Delta^2_{i-1}(p)}{\Delta_{i-2}(p)\cdot\Delta_{i}(p)},\quad
i=1,2,\dots,n,
\end{equation}
\item[2)]
If $n=2l+1$, then
\begin{equation}\label{coeff.Stieltjes.fraction.main.formula.quasi-stability.odd.degree}
c_{i}=\dfrac{\Delta^2_{i}(p)}{\Delta_{i-1}(p)\cdot\Delta_{i+1}(p)},\quad
i=0,1,2,\dots,n,
\end{equation}
\end{itemize}
where we set $\Delta_{-1}(p)\equiv\dfrac1{a_0}$ and
$\Delta_{0}(p)\equiv1$.

At last, from Theorems~\ref{Corol.differentiation.of.R-functions}
and~\ref{Theorem.main.Hurwitz.stability} one can obtain the
following simple result.
\begin{theorem}\label{Theorem.stability.with.difeerentiation}
Let $p$ be a Hurwitz polynomial of degree $n\geqslant2$ as
in~\eqref{main.polynomial}. Then all the polynomials
\begin{equation*}\label{Theorem.stability.with.difeerentiation.condition}
p_j(z)=\sum\limits_{i=0}^{n-2j}\left[\dfrac{n-i}2\right]\left(\left[\dfrac{n-i}2\right]-1\right)\cdots\left(\left[\dfrac{n-i}2\right]+j-1\right)a_iz^{n-2j-i},\quad
j=1,\ldots,\left[\dfrac{n}2\right]-1,
\end{equation*}
also are Hurwitz stable.
\end{theorem}
\begin{proof}
By Theorem~\ref{Theorem.main.Hurwitz.stability}, if $p$ is Hurwitz
stable, then the function $\Phi=p_1/p_0$ is an
\textit{R}-function. According to
Theorem~\ref{Corol.differentiation.of.R-functions}, all functions
$\Phi_j=p_1^{(j)}/p_0^{(j)}$,
$j=1,\ldots,\left[\dfrac{n}2\right]-1$, are \textit{R}-functions
with negative zeroes and poles. Now
Theorem~\ref{Theorem.main.Hurwitz.stability} implies that all
polynomials
\begin{equation*}\label{Theorem.stability.with.difeerentiation.proof.1}
p_j(z)=p_0^{(j)}(z^2)+zp_1^{(j)}(z^2),\quad
j=1,\ldots,\left[\dfrac{n}2\right]-1,
\end{equation*}
are Hurwitz stable, as required.
\end{proof}

\subsection{Quasi-stable polynomials}\label{subsection:quasi-stable.poly}

\hspace{4mm} In this section we deal with polynomials whose zeroes
lie in the \textit{closed} left half-plane of the complex~plane.

\begin{definition}\label{def.quasi-stable.poly}
A polynomial $p$ of degree $n$ is called \textit{quasi-stable}
with \textit{degeneracy index} $m$, $0\leqslant m\leqslant n$, if
all its zeroes lie in the \textit{closed} left half-plane of the
complex plane and the number of zeroes of $p$, counting
multiplicities, on the imaginary axis equals~$m$. We call the
number $n-m$ \textit{the stability index} of the polynomial~$p$.
\end{definition}
Throughout this section we use the following notation
\begin{equation}\label{floor.poly.stability.index}
r=\left[\dfrac n2\right]-\left[\dfrac m2\right],
\end{equation}
where $n$ and $m$ are degree and degeneracy index of the
polynomial $p$, respectively.

Obviously, any Hurwitz polynomial is quasi-stable with zero
degeneracy index, that is, it has the smallest degeneracy index
and the largest stability index (which equals the degree of the
polynomial). Note that if the degeneracy index $m$ is even, then
$p(0)\neq0$, and if $m$ is odd, then $p$ must have a zero at $0$.

Moreover, if $p$ is a quasi-stable polynomial, then

\begin{equation*}\label{Th.stable.poly.and.R-function.1}
p(z)=p_0(z^2)+zp_1(z^2)=f(z^2)q(z)=f(z^2)\left[q_0(z^2)+zq_1(z^2)\right],
\end{equation*}
where $f(u)$ is a real polynomial of degree
$\left[\dfrac{m}2\right]$ with nonpositive zeroes, and $q$ is a
Hurwitz stable polynomial if $m$ is even, and it is a quasi-stable
polynomial with degeneracy index $1$ if $m$ is odd. Using this
representation of quasi-stable polynomials, one can extend almost
all results of Section~\ref{subsection:Hurwitz.poly} to
quasi-stable polynomials in the same way, so we state them here
without proofs. We only should take into account that the function
\begin{equation}\label{assoc.function.quasi.stable}
\Phi(u)=\dfrac{p_1(u)}{p_0(u)}=\dfrac{f(u)q_1(u)}{f(u)q_0(u)}=\dfrac{q_1(u)}{q_0(u)}
\end{equation}
has a pole at $0$ whenever $p(0)=0$.

\vspace{3mm}

The analogue of Stodola necessary condition for quasi-stable
polynomials is the following.
\begin{theorem}\label{Th.Stodola.necessary.condition.quasi.stable}
If the polynomial $p$ defined in~\eqref{main.polynomial} is
quasi-stable, then all its coefficients are nonnegative.
\end{theorem}

The next theorem is an extended version of
Theorem~\ref{Theorem.main.Hurwitz.stability}.

\begin{theorem}\label{Th.stable.poly.and.R-function}
The polynomial $p$ of degree $n$ defined
in~\eqref{main.polynomial} is quasi-stable with degeneracy index
$m$ if and only if the~function~$\Phi$ defined
in~\eqref{assoc.function.quasi.stable} is an \textit{R}-function
of negative type with exactly $r$ poles all of which are
nonpositive, and $\displaystyle\lim_{u\to\pm\infty}\Phi(u)$ is
positive if $n=2l+1$. The number $r$ is defined
in~\eqref{floor.poly.stability.index}.
\end{theorem}

If we expand the function $\Phi$ into its Laurent
series~\eqref{app.assoc.function.series}, then the Hankel matrix
$S(\Phi)=\|s_{i+j}\|^{\infty}_{0}$ has rank~$r$, where $r$ as
in~\eqref{floor.poly.stability.index}.

\begin{theorem}\label{Th.quasi-stable.even.poly.and.sign.regularity}
Let the polynomial $p$ be defined in~\eqref{main.polynomial}. The
following conditions are equivalent:
\begin{itemize}
\item[1)] the polynomial $p$ is quasi-stable with \emph{even} degeneracy
index~$m$;
\item[2)] the following hold
\begin{equation*}\label{Th.quasi-stable.even.poly.and.sign.regularity.condition.0}
s_{-1}>0\quad\text{for}\quad n=2l+1,
\end{equation*}
\begin{equation*}\label{Th.quasi-stable.even.poly.and.sign.regularity.condition.1}
D_j(\Phi)>0, \qquad j=1,\ldots,r,
\end{equation*}
\begin{equation*}\label{Th.quasi-stable.even.poly.and.sign.regularity.condition.2}
(-1)^{j}\widehat{D}_{j}(\Phi)>0, \qquad j=1,\ldots,r-1,
\end{equation*}
and $(-1)^r\widehat{D}_r(\Phi)>0$ for even $m$, but
$\widehat{D}_r(\Phi)=0$ for odd $m$. Here $r$ is defined
in~\eqref{floor.poly.stability.index};
\item[3)] the matrix $S(\Phi)$ is sign regular of rank $r$ for even $m$, $S(\Phi)$ is $r-1$-sign regular
of rank $r$ for odd~$m$, and $s_{-1}$ is positive if $n=2l+1$.
\end{itemize}
\end{theorem}

It is also easy to extend
Theorem~\ref{Theorem.Hurwitz.stable.Hurwitz.matrix.criteria} to
quasi-stable polynomials.

\begin{theorem}\label{Th.quasi.stable.Hurwitz.matrix.criteria}
Let the polynomial $p$ of degree $n$ be defined
in~\eqref{main.polynomial}. The following conditions are
equivalent:
\begin{itemize}
\item[1)] the polynomial $p$ is quasi-stable with degeneracy index $m$;
\item[2)] determinants $\Delta_j(p)$ are positive up to order
$n-m$:
\begin{equation*}\label{Th.quasi.stable.Hurwitz.matrix.criteria.condition.1}
\Delta_1(p)>0,\ \Delta_2(p)>0,\dots,\ \Delta_{n-m}(p)>0,\
\Delta_{n-m+1}(p)=\ldots=\Delta_n(p)=0;
\end{equation*}
\item[3)] determinants $\eta_j(p)$ are positive up to order
$n-m+1$:
\begin{equation*}\label{Th.quasi.stable.Hurwitz.matrix.criteria.condition.2}
\eta_1(p)>0,\ \eta_2(p)>0,\dots,\ \eta_{n-m+1}(p)>0,\
\eta_{n-m+i}(p)=0,\quad i=2,3,\ldots;
\end{equation*}
\item[4)] the matrix $H_{\infty}(p)$ is totally nonnegative and
\begin{equation*}\label{Th.quasi.stable.Hurwitz.matrix.criteria.condition.4}
\eta_{n-m+1}(p)\neq0,\,\eta_{n-m+i}(p)=0,\quad i=2,3,\ldots
\end{equation*}
\end{itemize}
\end{theorem}

\noindent The implication $1)\Longrightarrow4)$ was proved,
indeed, in~\cite{Asner,{Kemperman}}. The implication
$4)\Longrightarrow1)$ seems to be new. Thus, we established that
\textit{a polynomial is quasi-stable if and only if its infinite
Hurwitz matrix is totally nonnegative}.

\vspace{2mm}

The following interesting property of quasi-stable polynomials is
a simple consequence of Theorems~\ref{Th.interlacity.preserving}
and~\ref{Theorem.main.Hurwitz.stability}.

\begin{theorem}\label{Th.quasi-stability. preserving.even}
Let the polynomials $p$ of degree $n=2l$ defined
in~\eqref{main.polynomial} be quasi-stable. Given any positive
integer $r(\leqslant n)$, the polynomial
\begin{equation*}\label{Th.quasi-stability. preserving.even.condition}
p_{r}(z)=a_0z^{2k}+a_{2r-1}z^{2k-1}+a_{2r}z^{2k-2}+a_{4r-1}z^{2k-3}+a_{4r}z^{2k-4}+\ldots+a_{2rk-1}z+a_{2rk},
\end{equation*}
where $k=\left[\dfrac{n}r\right]$, also is quasi-stable.
\end{theorem}
\begin{theorem}\label{Th.quasi-stability. preserving.odd}
Let the polynomials $p$ of degree $n=2l+1$ defined
in~\eqref{main.polynomial} be quasi-stable. Given any positive
integer $r(\leqslant n)$, the polynomial
\begin{equation*}\label{Th.quasi-stability. preserving.odd.condition}
p_{r}(z)=a_0z^{2k+1}+a_{1}z^{2k}+a_{2r}z^{2k-1}+a_{2r+1}z^{2k-2}+a_{4r}z^{2k-3}+a_{4r+1}z^{2k-4}+\ldots+a_{2rk}z+a_{2rk+1},
\end{equation*}
where $k=\left[\dfrac{n}r\right]$, also is quasi-stable.
\end{theorem}

Although total nonnegativity of the \textit{infinite} Hurwitz
matrix is equivalent to quasi-stability of polynomials, total
nonnegativity of the \textit{finite} singular Hurwitz matrix is
not equivalent to quasi-stability as was noticed by
Asner~\cite{Asner}. It is clear from
Theorem~\ref{Th.quasi.stable.Hurwitz.matrix.criteria} that the
finite Hurwitz matrix of a quasi-stable polynomial is totally
nonnegative. However, given a real polynomial $p$ of degree $n$,
if $\mathcal{H}_n(p)$ is totally nonnegative, $p$ is not
undertaken to be quasi-stable. In fact~\cite{Holtz_Tyaglov}, $p$
can be represented in the form $p(z)=q(z)g(z^2)$, where the
polynomial $q$ is Hurwitz stable, and the polynomial $g$ is chosen
such that the matrix $\mathcal{H}_n(p)$ is totally nonnegative but
$g(z^2)$ has zeroes in the open \textit{right} half-plane. Such
choice is possible as was mentioned in~\cite{Holtz_Tyaglov}.

\vspace{2mm}

The connection between quasi-stable polynomials and Stieltjes
continued fractions is similar to Hurwitz stable polynomials' one.
Namely, one can prove the following extended version of
Theorem~\ref{Theorem.Hurwitz.stable.Stieltjes.cont.frac.criteria}

\begin{theorem}\label{Th.quasi.stable.Stieltjes.cont.frac.criteria}
The polynomial $p$ of degree $n$ defined
in~\eqref{main.polynomial} is quasi-stable with degeneracy
index~$m$ if and only if the function $\Phi$ has a Stieltjes
continued fraction expansion:
\begin{equation}\label{Stieltjes.fraction.for.quasi-stability}
\Phi(u)=c_0+\dfrac1{c_1u+\cfrac1{c_2+\cfrac1{c_{3}u+\cfrac1{\ddots+\cfrac1{T}}}}},\quad
\text{with}\quad c_i>0,\quad T=\begin{cases}
         &c_{2r},\ \text{if}\ m\ \text{is even},\\
         &c_{2r-1}u,\ \text{if}\ m\ \text{is odd}.
       \end{cases}
\end{equation}
Here $c_0=0$ if $n$ is even, $c_0>0$ if $n$ is odd, and $r$ is
defined in~\eqref{floor.poly.stability.index}.
\end{theorem}
\noindent The coefficients $c_i$
in~\eqref{Stieltjes.fraction.for.quasi-stability} also can be
found by the
formul\ae~\eqref{coeff.Stieltjes.fraction.main.formula.quasi-stability.even.degree}--\eqref{coeff.Stieltjes.fraction.main.formula.quasi-stability.odd.degree}.

Finally, we note that a quasi-stable polynomial $p$ with
degeneracy index $1$ can be represented as follows: $p(z)=zq(z)$,
where $q$ is a Hurwitz stable polynomial. Therefore, for
quasi-stable polynomials with degeneracy index $1$, one can
establish an analogue of Li\'enard and Chipart criterion.
\begin{theorem}\label{Th.Lienard.Chipart.for.quasi-stability}
The polynomial $p$ given by~\eqref{main.polynomial} is
quasi-stable with degeneracy index $1$ if and only if one of the
following conditions holds
\begin{itemize}
\item[$1)$]
$a_n=0,a_{n-1}>0,a_{n-2}>0,a_{n-4}>0,\ldots,\quad\Delta_1(p)>0,\Delta_3(p)>0,\Delta_5(p)>0,\ldots$;
\item[$2)$]
$a_n=0,a_{n-1}>0,a_{n-3}>0,a_{n-5}>0,\ldots,\quad\Delta_1(p)>0,\Delta_3(p)>0,\Delta_5(p)>0,\ldots$;
\item[$3)$]
$a_n=0,a_{n-1}>0,a_{n-2}>0,a_{n-4}>0,\ldots,\quad\Delta_2(p)>0,\Delta_4(p)>0,\Delta_6(p)>0,\ldots$;
\item[$4)$]
$a_n=0,a_{n-1}>0,a_{n-3}>0,a_{n-5}>0,\ldots,\quad\Delta_2(p)>0,\Delta_4(p)>0,\Delta_6(p)>0,\ldots$
\end{itemize}
\end{theorem}

\setcounter{equation}{0}

\section{Self-interlacing polynomials}\label{section:self-interlacings}

\hspace{4mm} In Section~\ref{section:Hurwitz.polys} we established
that the function $\Phi$ (see~\eqref{assoc.function}) associated
with a Hurwitz stable polynomial (quasi-stable polynomial) maps
the upper half-plane of the complex plane to the lower half-plane
and possesses only negative (nonpositive) poles. Now we are in a
position to describe polynomials whose associated function $\Phi$
also maps the upper half-plane to the lower half-plane but has
only positive (nonnegative) poles.

\subsection{General theory}\label{subsection:self-interlacings.general.theory}

\begin{definition}\label{def.self-interlacing.poly}
A real polynomial $p(z)$ is called \textit{self-interlacing} if
all its zeroes are real and simple and interlace zeroes of the
polynomial $p(-z)$.
\end{definition}

In other words, if $\lambda_i$ are the zeroes of a
self-interlacing polynomial $p$, then one of the following holds:
\begin{equation}\label{self-interlacing.zero.distribution.I.type}
0<\lambda_1<-\lambda_2<\lambda_3<\ldots<(-1)^{n-1}\lambda_n,
\end{equation}
\begin{equation}\label{self-interlacing.zero.distribution.II.type}
0<-\lambda_1<\lambda_2<-\lambda_3<\ldots<(-1)^{n}\lambda_n,
\end{equation}
where $n=\deg p$.

\begin{definition}\label{def.self-interlacing.poly.I.type}
A real polynomial $p$ of degree $n$ is called
\textit{self-interlacing of type~I} (\textit{of type~II}) if all
its zeroes are real and simple and satisfy the
inequalities~\eqref{self-interlacing.zero.distribution.I.type}
(respectively,~\eqref{self-interlacing.zero.distribution.II.type}).
\end{definition}

If a polynomial $p$ of degree $n$ is self-interlacing of type~I,
then its minimal absolute value zero is positive. Moreover, if
$n=2l$, then $p$ has exactly $l$ negative zeroes and exactly $l$
positive zeroes. Its zero $\lambda_n$, which has the maximal
absolute value is negative. If $n=2l+1$, then the polynomial $p$
has exactly $l$ negative zeroes and $l+1$ positive zeroes. In this
case, $\lambda_n$ is positive.

Note that a polynomial $p(z)$ is self-interlacing of type~I if and
only if the polynomial $p(-z)$ is self-interlacing of type~II. So
in the sequel, we deal only with self-interlacing polynomials of
type~I.

\begin{theorem}\label{Theorem.main.self-interlacing}
Let $p$ be a real polynomial of degree $n\geqslant1$ as
in~\eqref{main.polynomial}. The polynomial $p$ is self-interlacing
of type~I if and only if its associated function~$\Phi$ defined
in~\eqref{assoc.function} is an \textit{R}-function with
exactly~$l$ poles, all of which are positive, and
$\displaystyle\lim_{u\to\pm\infty}\Phi(u)$ is negative whenever
$n=2l+1$.
\end{theorem}
\begin{proof}
Let $p$ be self-interlacing of type~I. First, we show that the
function
\begin{equation*}\label{Theorem.main.self-interlacing.proof.0}
G(z)=-\dfrac{p(-z)}{p(z)}
\end{equation*}
is an \textit{R}-function. In fact, by
Definition~\ref{def.self-interlacing.poly}, the zeroes of the
polynomials $p(z)$ and $p(-z)$ are real, simple and interlacing,
that is, between any two consecutive zeroes of one polynomial
there lies exactly one zero, counting multiplicity, of the other
polynomial, so $G(z)$ or $-G(z)$ is an \textit{R}-function
according to Theorem~\ref{Th.R-function.general.properties}. By
Corollary~\ref{Corol.monotonicity.R.functions}, $G$ is monotone
between its poles. So it remains to prove that the function $G$ is
decreasing between its poles.

Let $n=2l$, then
by~\eqref{self-interlacing.zero.distribution.I.type}, the maximal
pole of $G$ is $\lambda_{n-1}>0$, but its maximal zero is
$-\lambda_n>0$, which is greater than $\lambda_{n-1}$ according to
Definition~\ref{def.self-interlacing.poly.I.type}
(see~\eqref{self-interlacing.zero.distribution.I.type}).
Therefore, in the interval $(-\lambda_n,+\infty)$, the
function~$G$ has no poles and zeroes. At the same time,
$\displaystyle\lim_{z\to\pm\infty}G(z)=-1$. Consequently, in the
interval~$(-\lambda_n,+\infty)$, $G(z)$ decreases from $0$ to
$-1$. So, $G$ is an \textit{R}-function. In the same way, one can
prove that if $n=2l+1$, then $G$ is also an~\textit{R}-function.

Thus, if $p$ is self-interlacing of type~I, then $G$ maps the
upper half-plane to the lower half-plane. From~\eqref{poly1.1} we
obtain that
\begin{equation}\label{Theorem.main.self-interlacing.proof.1}
z\Phi(z^2)=\dfrac{1+G(z)}{1-G(z)}.
\end{equation}
Since the fractional linear transformation $\displaystyle
z\mapsto\frac{1+z}{1-z}$ conformally maps the lower half-plane to
the lower half-plane:
\begin{equation}\label{Theorem.main.self-interlacing.proof.2}
\Im z<0 \implies\Im\left(\frac{1-z}{1+z}\right)=\dfrac{2\Im z
}{|1+z|^2}<0,
\end{equation}
from~\eqref{Theorem.main.self-interlacing.proof.1}--\eqref{Theorem.main.self-interlacing.proof.2}
it follows that the function $z\Phi(z^2)$ maps the upper
half-plane to the lower half-plane, that is, $z\Phi(z^2)$ is
an~\textit{R}-function.

Since $p$ is self-interlacing of type~I by assumption, the
polynomials $p(z)$ and $p(-z)$ have no common zeroes, therefore
$p_0$ and $p_1$ also have no common zeroes and $p_0(0)\neq0$
by~\eqref{app.poly.odd.even}. Therefore, the number of poles of
the function $z\Phi(z^2)$ equals the number of zeroes of the
polynomial $p_0(z^2)$. If $n=2l$, then by~\eqref{poly1.13}, $\deg
p_0=l$, so $z\Phi(z^2)$ has exactly $2l$ poles. If $n=2l+1$, then
by~\eqref{poly1.12}, $\deg p_1=l$ and $\deg p_0\leqslant l$, so
$z\Phi(z^2)$ has at most $2l$ poles, and it has exactly $2l+1$
zeroes, since $p_0(0)\neq0$ as we mentioned above.
But~$z\Phi(z^2)$ is an~\textit{R}-function, therefore, it has
exactly $2l$ poles by
Theorem~\ref{Th.R-function.general.properties}.

Thus, Theorem~\ref{Th.R-function.general.properties} implies that
the function $z\Phi(z^2)$ can be represented in the
form~\eqref{Mittag.Leffler.1}, where all poles are located
symmetrically with respect to~$0$ and $\beta=0$, since
$z\Phi(z^2)$ is an odd function. Denote the poles of $z\Phi(z^2)$
by $\pm\nu_1,\ldots,\pm\nu_l$ such that
\begin{equation*}\label{Theorem.main.self-interlacing.proof.3}
0<\nu_1<\nu_2<\ldots<\nu_l.
\end{equation*}
Note that $\nu_1\neq0$ since $p_0(0)\neq0$.

So the function $z\Phi(z^2)$ can be represented in the following
form
\begin{equation*}\label{Theorem.main.self-interlacing.proof.4}
z\Phi(z^2)=-\alpha z+\sum_{j=1}^l\frac{\gamma_j}{z-\nu_j}+
\sum_{j=1}^l\frac{\gamma_j}{z+\nu_j}= -\alpha z
+\sum_{j=1}^l\frac{2\gamma_jz}{z^2-\nu_j^2},
\quad\alpha\geqslant0,\,\gamma_j,\nu_j>0.
\end{equation*}
Divide this equality by $z$ and changing variables as follows
$z^2\to u$, $2\gamma_j\to\beta_j$, $\nu_j^2\to\omega_j$, we obtain
the following representation of the function $\Phi$:
\begin{equation}\label{Mittag.Leffler.self-interlacing}
\Phi(u)=-\alpha+\sum_{j=1}^l\frac{\beta_j}{u-\omega_j},
\end{equation}
where $\alpha\geqslant0,\,\beta_j>0$ and
\begin{equation*}\label{Theorem.main.self-interlacing.proof.5}
0<\omega_1<\omega_2<\ldots<\omega_l.
\end{equation*}
Here $\alpha=0$ whenever $n=2l$, and
$-\alpha=\displaystyle\frac{a_0}{a_1}<0$ whenever $n=2l+1$. Since
$\Phi$ can be represented in the
form~\eqref{Mittag.Leffler.self-interlacing}, by
Theorem~\ref{Th.R-function.general.properties}, $\Phi$ is an
\textit{R}-function with exactly $l$ poles. Moreover,
from~\eqref{Mittag.Leffler.self-interlacing} it also follows that
all poles $\Phi$ are positive and
$\displaystyle\lim_{u\to\pm\infty}\Phi(u)<0$ if $n=2l+1$.

\vspace{3mm}

Conversely, let the polynomial $p$ be defined
in~\eqref{main.polynomial}, and let its associated function $\Phi$
be an \textit{R}-function with exactly $l$ poles, all of which are
positive, and let $\displaystyle\lim_{u\to\pm\infty}\Phi(u)<0$
when $n=2l+1$. We will show that $p$ is self-interlacing of
type~I.

By Theorem~\ref{Th.R-function.general.properties}, $\Phi$ can be
represented in the form~\eqref{Mittag.Leffler.self-interlacing},
where $-\alpha=\displaystyle\lim_{u\to\pm\infty}\Phi(u)\leqslant0$
such that $\alpha>0$ if $n=2l+1$, and $\alpha=0$ if $n=2l$. Thus,
the polynomial $p_0$ has only positive zeroes, and the polynomials
$p_0$ and $p_1$ have no common zeroes. Together
with~\eqref{app.poly.odd.even} and~\eqref{assoc.function}, this
means that the set of zeroes of the~polynomial $p$ coincides with
the set of roots of the equation
\begin{equation}\label{Theorem.main.self-interlacing.proof.6}
z\Phi(z^2)=-1.
\end{equation}
At first, we study real zeroes of this equation, so we consider
only real $z$. Since $p_0$ has only positive zeroes, we have
$p(0)=p_0(0)\neq0$. Thus, we put $z\in\mathbb{R}\backslash\{0\}$.
Changing variables as follows $z^2\to u>0$, we rewrite the
equation~\eqref{Theorem.main.self-interlacing.proof.6} in the
following form, which is equivalent
to~\eqref{Theorem.main.self-interlacing.proof.6} for real
nonzero~$z$:
\begin{equation}\label{Theorem.main.self-interlacing.proof.7}
\begin{cases}
\ \Phi(u)=-\dfrac1{\sqrt{u}},\\
\ \Phi(u)=\dfrac1{\sqrt{u}},
\end{cases}\quad u>0.
\end{equation}
Note that all positive roots (if any) of the first
equation~\eqref{Theorem.main.self-interlacing.proof.7} are squares
of the positive roots of
the~equation~\eqref{Theorem.main.self-interlacing.proof.6}, and
all positive roots (if any) of the second
equation~\eqref{Theorem.main.self-interlacing.proof.7} are squares
of the negative roots of
the~equation~\eqref{Theorem.main.self-interlacing.proof.6}.

Let $n=2l$. Then $\lim\limits_{u\to\pm\infty}\Phi(u)=-\alpha=0$.
Consider the function $F_1(u)=\Phi(u)+\dfrac1{\sqrt{u}}$ whose set
of zeroes coincides with the set of roots of the first
equation~\eqref{Theorem.main.self-interlacing.proof.7}. %
The function $F_1(u)$ has the same positive poles as $\Phi(u)$
does, that is, $\omega_1,\omega_2,\ldots,\omega_l$. Moreover,
$F_1(u)$ is decreasing on the intervals
$(0,\omega_1),(\omega_1,\omega_2),\ldots,(\omega_{l-1},\omega_l)$,
$(\omega_l,+\infty)$ as a sum of functions that are decreasing on
those intervals. Since~$F_1$ is decreasing on $(\omega_l,+\infty)$
and $F_1(u)\to+0$ as $u\to+\infty$, we have $F_1(u)>0$
for~$u\in(\omega_l,+\infty)$. Further, $F_1(u)\to-\infty$ as
$u\nearrow\omega_i$ and $F_1(u)\to+\infty$ as $u\searrow\omega_i$,
$i=1,\ldots,l$. Also $F_1(u)\to+\infty$ whenever $u\to+0$, since
$|\Phi(0)|=\left|\dfrac{p_1(0)}{p_0(0)}\right|<\infty$ and
$\dfrac1{\sqrt{u}}\to+\infty$ as $u\to+0$. Thus, on each of the
intervals
$(0,\omega_1),(\omega_1,\omega_2),\ldots,(\omega_{l-1},\omega_l)$,
the function $F_1$ decreases from $+\infty$ to $-\infty$, and
therefore, it has exactly one zero, counting multiplicity, on each
of those intervals. Denote the zero of $F_1$ on the interval
$(\omega_{i-1},\omega_{i})$ by $\mu_i^2$, $i=2,\ldots,l$
($\mu_i>0$). Also we denote by $\mu_1^2$ ($\mu_1>0$) the zero of
$F_1$ on the interval $(0,\omega_1)$. Since $\Phi$ is an
\textit{R}-function by assumption, it has exactly one zero,
counting multiplicity, say $\xi_i$, on each of the intervals
$(\omega_{i},\omega_{i+1})$, $i=1,\ldots,l-1$. But
$F_1(\xi_i)=\dfrac1{\sqrt{\xi_i}}>0$, so we have
\begin{equation}\label{Theorem.main.self-interlacing.proof.8}
\mu^2_1<\omega_1<\xi_1<\mu^2_2<\omega_{2}<\xi_2<\ldots<\xi_{l-1}<\mu_l^2<\omega_l.
\end{equation}
Thus, the first
equation~\eqref{Theorem.main.self-interlacing.proof.7} has exactly
$l$ positive roots $\mu_i$, all of which are simple and satisfy
the inequalities~\eqref{Theorem.main.self-interlacing.proof.8}.

Consider now the function $F_2(u)=\Phi(u)-\dfrac1{\sqrt{u}}$ whose
set of zeroes coincides with the set of roots of the second
equation~\eqref{Theorem.main.self-interlacing.proof.7}. Since
$\Phi(0)=-\displaystyle\sum_{i=1}^{l}\dfrac{\beta_i}{\omega_i}<0$
by~\eqref{Mittag.Leffler.self-interlacing}, and
$\Phi(u)\to-\infty$ as $u\nearrow\omega_1$, the function $\Phi$
decreases from $\Phi(0)<0$ to $-\infty$ on the interval
$(0,\omega_1)$. Consequently, $F_2(u)<0$ for all
$u\in(0,\omega_1)$. It is clear that $F_2(u)\to+\infty$ as
$u\searrow\omega_i$ and $F_2(u)\to-\infty$ as $u\nearrow\omega_i$,
$i=1,\ldots,l$. Since $F_2(\mu^2_i)=-\dfrac1{\mu_i}<0$, the
function $F_2$ has an odd number of zeroes, counting
multiplicities, on each of the intervals $(\omega_i,\mu^2_i)$,
$i=1,\ldots,l-1$. Besides, $F_2(u)=F_1(u)-\dfrac2{\sqrt{u}}<0$
whenever $u\in[\mu^2_i,\omega_{i+1})$, $i=1,\ldots,l-1$, because
$F_1(u)<0$ on those intervals. More exactly, since
$F_2(\xi_i)=\dfrac1{\sqrt{\xi_i}}<0$, $i=1,\ldots,r-1$, the
function $F_2(u)$ has an odd number of zeroes on each interval
$(\omega_i,\xi_i)$ and an even number of zeroes on each interval
$[\xi_i,\mu^2_{i+1})$, $i=1,\ldots,r-1$. So $F_2$ has \textit{at
least} $l-1$ zeroes in the interval $(\omega_1,\omega_l)$.
Consider now the interval $(\omega_l,+\infty)$ and show that
$F_2(u)\to-0$ as $u\to+\infty$. In fact, since the function $\Phi$
is an \textit{R}-function by assumption, we have
$s_0=\lim\limits_{u\to\pm\infty}u\Phi(u)=D_1(\Phi)>0$ according
Theorem~\ref{Th.R-function.general.properties}, where $s_0$ is the
first coefficient\footnote{Recall that $s_{-1}=0$ whenever
$n=2l$.} in the series~\eqref{app.assoc.function.series}.
Therefore, $\Phi(u)\sim u^{-1}$ as $u\to+\infty$. This implies the
following: $F_2(u)\sim-u^{-\tfrac12}$ as $u\to+\infty$. Thus,
$F_2(\omega_l+\varepsilon)>0$ for all sufficiently small
$\varepsilon>0$ and $F_2(u)<0$ for all sufficiently large
positive~$u$. Consequently, $F_2$ has an odd number of zeroes,
counting multiplicities, in the interval $(\omega_l,+\infty)$. So
$F_2$ has \textit{at least} $l$ zeroes in the interval
$(\omega_1,+\infty)$, and it has no zeroes in $(0,\omega_1)$.

Since the first
equation~\eqref{Theorem.main.self-interlacing.proof.7} has exactly
$l$ positive roots and the second
equation~\eqref{Theorem.main.self-interlacing.proof.7} has at
least $l$ positive roots, the
equation~\eqref{Theorem.main.self-interlacing.proof.6} has
\textit{at least} $2l$ real roots, all of which are the zeroes of
the polynomial~$p$. But $\deg p=2l$ by assumption, therefore, the
second equation~\eqref{Theorem.main.self-interlacing.proof.7} also
has exactly $l$ positive roots. Moreover, it has exactly one
simple root in each of the intervals $(\omega_i,\xi_i)$,
$i=1,\dots,l-1$, and exactly one simple root in the interval
$(\omega_l,+\infty)$. Thus, denoting the positive roots of the
second equation~\eqref{Theorem.main.self-interlacing.proof.7} by
$\zeta_i^2$, $i=1,\ldots,l$ ($\zeta_i<0$), we have
\begin{equation}\label{Theorem.main.self-interlacing.proof.9}
\omega_1<\zeta_1^2<\xi_1<\omega_{2}<\zeta_2^2<\xi_2<\ldots<\xi_{l-1}<\omega_l<\zeta_l^2.
\end{equation}
Now
from~\eqref{Theorem.main.self-interlacing.proof.8}--\eqref{Theorem.main.self-interlacing.proof.9}
we obtain
\begin{equation}\label{Theorem.main.self-interlacing.proof.10}
0<\mu_1<-\zeta_1<\mu_2<-\zeta_2\ldots<\mu_{l}<-\zeta_l.
\end{equation}
Recall that $\mu_i>0$ are the positive roots of
the~equation~\eqref{Theorem.main.self-interlacing.proof.6}, and
$\zeta_i<0$ are the negative roots of
the~equation~\eqref{Theorem.main.self-interlacing.proof.6}. Thus,
all roots of the
equation~\eqref{Theorem.main.self-interlacing.proof.6} (and
therefore, all zeroes of the polynomial $p$) are real and simple.
Denote them by $\lambda_i$ and enumerate such that
\begin{equation*}\label{Theorem.main.self-interlacing.proof.11}
0<|\lambda_1|<|\lambda_2|\ldots<|\lambda_n|.
\end{equation*}
Then we have $\lambda_{2i-1}=\mu_i>0$ and
$\lambda_{2i}=\zeta_i<0$, $i=1,\ldots,l$, and the
inequalities~\eqref{Theorem.main.self-interlacing.proof.10}
imply~\eqref{self-interlacing.zero.distribution.I.type}. Thus, $p$
is a self-interlacing polynomial of type~I.

In the same way, one can show that if $n=2l+1$, the function
$\Phi$ is an \textit{R}-function with positive poles, and
$\lim\limits_{u\to\pm\infty}\Phi(u)<0$, then the polynomial $p$ is
self-interlacing of type~I.
\end{proof}

\begin{remark}
Let us note that if the even and odd parts, $p_0$ and $p_1$, of a
given polynomial $p$ have positive interlacing zeroes, then the
polynomials $p_0(z^2)$ and $zp_1(z^2)$ have real interlacing
zeroes, and therefore the polynomial $p(z)=p_0(z^2)+zp_1(z^2)$ has
real zeroes by Theorem~\ref{Th.R-function.general.properties} and
Corollary~\ref{remark.R-func.pos.type} as a linear combination of
polynomials with real interlacing zeroes. However, this notice
does not help to investigate the self-interlacing property of
polynomials.
\end{remark}

\begin{remark}
In the proof of Theorem~\ref{Theorem.main.self-interlacing}, we
also established that the squares of both positive and negative
zeroes of the self-interlacing polynomial $p$ interlace both
zeroes of $p_0$ and zeroes of $p_1$.
\end{remark}

Theorem~\ref{Theorem.main.self-interlacing} allow us to use
properties of \textit{R}-functions in order to obtain additional
criteria of self-interlacing. At first, consider the Hankel matrix
$S(\Phi)=\|s_{i+j}\|^{\infty}_{i,j=0}$ constructed with the
coefficients of the series~\eqref{app.assoc.function.series}. From
Theorems~\ref{Th.R-function.general.properties},~\ref{Theorem.main.self-interlacing}
and~\ref{Theorem.total.nonnetativity.Hankel.matrix} and from
Corollary~\ref{Corol.R-functions.pos.poles} we obtain the
following self-interlacing criteria:
\begin{theorem}\label{Theorem.self.interlacing.poly.and.total.positivity}
Let a real polynomial $p$ be defined by~\eqref{main.polynomial}.
The following conditions are equivalent:
\begin{itemize}
\item[1)] the polynomial $p$ is self-interlacing of type~I;
\item[2)] the following hold
\begin{equation}\label{Theorem.self.interlacing.poly.and.total.positivity.condition.0}
\begin{split}
&s_{-1}<0\quad\text{for}\quad n=2l+1,\\
&D_j(\Phi)>0, \qquad j=1,\ldots,l,\\
& \widehat{D}_{j}(\Phi)>0, \qquad j=1,\ldots,l,
\end{split}
\end{equation}
where $l=\left[\dfrac n2\right]$;
\item[3)] the matrix $S(\Phi)$ is strictly totally positive of rank $l$, and $s_{-1}<0$ whenever $n=2l+1$.
\end{itemize}
The determinants $D_j(\Phi)$ and $\widehat{D}_j(\Phi)$ are defined
in~\eqref{Hankel.determinants.1}
and~\eqref{Hankel.determinants.2}, respectively.
\end{theorem}
\begin{proof}
According to Theorem~\ref{Theorem.main.self-interlacing}, the
polynomial $p$ is self-interlacing of type~I if and only if its
associated function $\Phi$ is an \textit{R}-function with exactly
$l=\left[\dfrac n2\right]$ poles, all of which are positive, and
$\lim\limits_{u\to\pm\infty}\Phi(u)=s_{-1}\leqslant0$. By
Theorem~\ref{Th.R-function.general.properties} and
Corollary~\ref{Corol.R-functions.pos.poles}, this is equivalent to
the
inequalities~\eqref{Theorem.self.interlacing.poly.and.total.positivity.condition.0}.
But these inequalities are equivalent to the strictly total
positivity of the matrix $S(\Phi)$ by
Theorem~\ref{Theorem.total.nonnetativity.Hankel.matrix}.
\end{proof}

From this theorem we obtain the following criterion of
self-interlacing, which is an analogue of the~Hurwitz stability
criterion.
\begin{theorem}\label{Theorem.self-interlacing.Hurwitz.criterion}
Given a real polynomial $p$ of degree $n$ as
in~\eqref{main.polynomial}, the following conditions are
equivalent:
\begin{itemize}
\item[1)] $p$ is self-interlacing of type~I;
\item[2)] the Hurwitz determinants $\Delta_j(p)$ satisfy the
inequalities:
\begin{equation}\label{Hurvitz.det.noneq.self-interlacing.1}
\Delta_{n-1}(p)>0,\ \Delta_{n-3}(p)>0,\ldots,
\end{equation}
\begin{equation}\label{Hurvitz.det.noneq.self-interlacing.2}
\displaystyle(-1)^{\left[\tfrac{n+1}2\right]}\Delta_{n}(p)>0,
\displaystyle(-1)^{\left[\tfrac{n+1}2\right]-1}\Delta_{n-2}(p)>0,\ldots
\end{equation}
\item[3)] the determinants $\eta_j(p)$ up to order $n+1$ satisfy the
inequalities:
\begin{equation}\label{Hurvitz.det.noneq.self-interlacing.3}
\eta_{n}(p)>0,\ \eta_{n-2}(p)>0,\ldots,\ \eta_1(p)=a_0>0,
\end{equation}
\begin{equation}\label{Hurvitz.det.noneq.self-interlacing.4}
\displaystyle(-1)^{\left[\tfrac{n+1}2\right]}\eta_{n+1}(p)>0,
\displaystyle(-1)^{\left[\tfrac{n+1}2\right]-1}\eta_{n-1}(p)>0,\ldots
\end{equation}
\end{itemize}
\end{theorem}
\begin{proof} We establish the equivalence $1)\Longleftrightarrow2)$ while the equivalence
$2)\Longleftrightarrow3)$ follows
from~\eqref{Formulae.Gurwitz.3333}.

By
Theorem~\ref{Theorem.self.interlacing.poly.and.total.positivity},
$p$ is self-interlacing of type~I if and only if the
inequalities~\eqref{Theorem.self.interlacing.poly.and.total.positivity.condition.0}
hold. Now
from~\eqref{Theorem.self.interlacing.poly.and.total.positivity.condition.0}
and~\eqref{Formulae.Gurwitz.1}--\eqref{Formulae.Gurwitz.2} we
obtain that $p$ is self-interlacing of type~I if and only if the
inequalities~\eqref{Hurvitz.det.noneq.self-interlacing.1} hold and
\begin{equation*}\label{Theorem.self-interlacing.Hurwitz.criterion.proof.1}
\begin{array}{l}
\text{for}\qquad n=2l\qquad\qquad\qquad
(-1)^j\Delta_{2j}(p)>0,\qquad j=1,\ldots,l;\\
 \\
\text{for}\qquad
n=2l+1\qquad\qquad\;(-1)^{j+1}\Delta_{2j+1}(p)>0,\qquad
j=0,1,\ldots,l,
\end{array}
\end{equation*}
that is equivalent
to~\eqref{Hurvitz.det.noneq.self-interlacing.2}.
\end{proof}
Note that the
inequalities~\eqref{Hurvitz.det.noneq.self-interlacing.2}
and~\eqref{Hurvitz.det.noneq.self-interlacing.4} are equivalent to
the following ones:
\begin{equation*}\label{Hurvitz.det.noneq.self-interlacing.5}
\Delta_{n-2i}(p)\Delta_{n-2i-2}(p)<0,\quad\text{and}\quad\eta_{n-2i+1}(p)\eta_{n-2i-1}(p)<0,\quad
i=0,1,\ldots,\left[\dfrac{n-1}2\right],
\end{equation*}
where $\Delta_0(p)=\eta_0(p)=a_0\Delta_{-1}(p)\equiv1$.

Analogues of
Theorems~\ref{Th.Lienard.Chipart}--\ref{Theorem.stability.with.difeerentiation}
will be established in Section~\ref{subsection:SI-Lienard-Chipart}
using a connection between Hurwitz stable polynomials and
self-interlacing polynomials of type~I.

\subsection{Interrelation between Hurwitz stable and self-interlacing
polynomials}\label{subsection:connection.Hurwitz.SI}

\hspace{4mm} Comparing
Theorems~\ref{Theorem.main.Hurwitz.stability}
and~\ref{Theorem.main.self-interlacing}, one obtains the following
fact, which has deep consequences and allows us to describe a lot
of properties of self-interlacing polynomials.

\begin{theorem}\label{Theorem.connection.Hurwitz.self-interlacing}
A polynomial $p(z)=p_0(z^2)+zp_1(z^2)$ is self-interlacing of type
I if and only if the polynomial $q(z)=p(-z^2)-zp_1(z^2)$ is
Hurwitz stable.
\end{theorem}
\begin{proof}
By Theorem~\ref{Theorem.main.self-interlacing}, $p$ is
self-interlacing if and only if the function
$\Phi(u)=\dfrac{p_1(u)}{p_0(u)}$ is an \textit{R}-function with
only positive poles, and
$\lim\limits_{u\to\pm\infty}\Phi(u)\leqslant0$ that is equivalent
to the fact that the function $\Psi(u)=-\dfrac{p_1(-u)}{p_0(-u)}$
is an \textit{R}-function with only negative poles, and
$\lim\limits_{u\to\pm\infty}\Psi(u)\geqslant0$. According to
Theorem~\ref{Theorem.main.Hurwitz.stability}, this is equivalent
to the polynomial $q$ being Hurwitz stable, as required.
\end{proof}

\begin{remark}
Thus, we have that there is one-to-one correspondence between
self-interlacing and Hurwitz stable polynomials. Given a real
Hurwitz stable polynomials, we should change sings of some its
coefficients to obtain a self-interlacing polynomial, and vise
versa.
\end{remark}

As an immediate consequence of
Theorem~\ref{Theorem.connection.Hurwitz.self-interlacing}, we
obtain the following analogue of Stodola's theorem,
Theorem~\ref{Th.Stodola.necessary.condition.Hurwitz}.

\begin{theorem}\label{Th.Stodola.necessary.condition.self-interlacing}
If the polynomial
\begin{equation}\label{Th.Stodola.necessary.condition.self-interlacing.poly}
p(z)=a_0z^n+a_1z^{n-1}+\dots+a_n,\qquad
a_1,\dots,a_n\in\mathbb{R},\ a_0>0,
\end{equation}
is self-interlacing of type~I, then
\begin{itemize}
\item[] for $n=2l$,
\begin{equation*}\label{Th.Stodola.necessary.condition.self-interlacing.even}
(-1)^{\tfrac{j(j-1)}2}a_j>0,\qquad j=0,1,\ldots,n;
\end{equation*}
\item[] for $n=2l+1$,
\begin{equation*}\label{Th.Stodola.necessary.condition.self-interlacing.odd}
(-1)^{\tfrac{j(j+1)}2}a_j>0,\qquad j=0,1,\ldots,n.
\end{equation*}
\end{itemize}
\end{theorem}
\begin{proof}
Since $p$ is self-interlacing of type~I, by
Theorem~\ref{Theorem.connection.Hurwitz.self-interlacing}, the
polynomial\footnote{We put the factor $(-1)^{\tfrac{n(n+1)}2}$ in
order to make the leading coefficient of the polynomial $q$ equal
to $a_0>0$.}
\begin{equation}\label{Th.Stodola.necessary.condition.self-interlacing.poly.proof}
q(z)=(-1)^{\tfrac{n(n+1)}2}[p_0(-z^2)-zp_1(-z^2)]=b_0z^n+b_1z^{n-1}+b_2z^{n-2}+\dots+b_n,\quad
b_0=a_0>0,
\end{equation}
is Hurwitz stable. It is easy to see that
\begin{equation}\label{Th.Stodola.necessary.condition.self-interlacing.even.proof}
b_j=(-1)^{\tfrac{j(j-1)}2}a_j,\qquad
j=0,1,\ldots,n\qquad\qquad\qquad\quad\quad\text{for}\qquad n=2l;
\end{equation}
and
\begin{equation}\label{Th.Stodola.necessary.condition.self-interlacing.odd.proof}
\quad\,\,\,\,\qquad\qquad\qquad
b_j=(-1)^{\tfrac{j(j+1)}2}a_j,\qquad
j=0,1,\ldots,n\qquad\qquad\qquad\qquad\text{for}\qquad n=2l+1.
\end{equation}
By Stodola's theorem,
Theorem~\ref{Th.Stodola.necessary.condition.Hurwitz}, all $b_j$
are positive.
\end{proof}

Thus, a necessary form of a self-interlacing polynomial with
positive leading coefficient is as follows:\\
for $n=2l$
\begin{equation*}
p(z)=b_0z^n+b_1z^{n-1}-b_2z^{n-2}-b_3z^{n-3}+b_4z^{n-4}+b_5z^{n-5}-b_6z^{n-6}-\dots,\qquad
b_0,\dots,b_n>0,
\end{equation*}
for $n=2l+1$
\begin{equation*}
p(z)=b_0z^n-b_1z^{n-1}-b_2z^{n-2}+b_3z^{n-3}+b_4z^{n-4}-b_5z^{n-5}-b_6z^{n-6}+\dots,\qquad
b_0,\dots,b_n>0.
\end{equation*}
We also notice that if the polynomial $p$ defined
in~\eqref{Th.Stodola.necessary.condition.self-interlacing.poly} is
self-interlacing, then
\begin{equation}\label{self-interlacing.last.coeffs}
a_{n-1}a_n<0,
\end{equation}
since $\Psi(0)=-\Phi(0)=-{p_1(0)}/{p_0(0)}=-{a_{n-1}}/{a_n}>0$
(see the proof of
Theorem~\ref{Th.Stodola.necessary.condition.self-interlacing}).

\begin{remark}
In fact,
Theorem~\ref{Th.Stodola.necessary.condition.self-interlacing} was
proved in~\cite{Fisk} by another methods.
\end{remark}

Let us point out at one more interesting connection between
Hurwitz stable and self-interlacing polynomials. Let the
polynomial $p$ be self-interlacing. If $p(z)=p_0(z^2)+zp_1(z^2)$,
then $p(iz)=p_0(-z^2)+izp_1(-z^2)$, where $i=\sqrt{-1}$.
Consequently, the Hurwitz stable polynomials
$q(z)=p_0(-z^2)-zp_1(-z^2)$ can be represented as follows
\begin{equation*}\label{self-in.stable.connection.one.more}
q(z)=\dfrac{p(iz)+p(-iz)}2-\dfrac{p(iz)-p(-iz)}{2i}=p(iz)\dfrac{1+i}2+p(-iz)\dfrac{1-i}2.
\end{equation*}
So one can establish the following theorem.
\begin{theorem}\label{Theorem.SI_connection.second}
Let the polynomial $p(z)=p_0(z^2)+zp_1(z^2)$ be self-interlacing
of type~I and let its dual polynomial $q(z)=p_0(-z^2)-zp_1(-z^2)$
be the Hurwitz stable polynomial associated with $p$. Then
\begin{equation*}
p(\lambda)=0\Longleftrightarrow\arg
q(i\lambda)=\dfrac\pi4\,\,\,\text{or}\,\,\,\dfrac{5\pi}4;
\end{equation*}
and, respectively,
\begin{equation*}
q(\mu)=0\Longleftrightarrow\arg
p(i\mu)=\dfrac\pi4\,\,\,\text{or}\,\,\,\dfrac{5\pi}4.
\end{equation*}
\end{theorem}
\begin{proof}
In fact, $p(\lambda)=0$ if and only if $-\dfrac{\lambda
p_1(\lambda^2)}{p_0(\lambda^2)}=1$. At the same time,
$q(i\lambda)=p_0(\lambda^2)-i\lambda p_1(\lambda^2)$.
Consequently, $\arg q(i\lambda)=\arctan\left(-\dfrac{\lambda
p_1(\lambda^2)}{p_0(\lambda^2)}\right)=\arctan1=\dfrac\pi4\,\text{or}\,\dfrac{5\pi}4$.

The second assertion of the theorem can be proved analogously.
\end{proof}

\subsection{Li\'enard and Chipart criterion, Stieltjes continued fractions and
the signs of Hurwitz minors}\label{subsection:SI-Lienard-Chipart}

\hspace{4mm} Now we consider the polynomials
$p(z)=p_0(z^2)+zp_1(z^2)$ and $q(z)=p_0(-z^2)-zp_1(-z^2)$ and
their associated functions $\Phi(u)=\dfrac{p_1(u)}{p_0(u)}$ and
$\Psi(u)=-\dfrac{p_1(-u)}{p_0(-u)}$. Let
\begin{equation*}\label{app.assoc.function.series.self-int.1}
\Phi(u)=s_{-1}+\frac{s_0}u+\frac{s_1}{u^2}+\frac{s_2}{u^3}+\frac{s_3}{u^4}+\dots
\end{equation*}
Then we have
\begin{equation*}\label{app.assoc.function.series.self-int.2}
\Psi(u)=-\dfrac{p_1(-u)}{p_0(-u)}=t_{-1}+\frac{t_0}w+\frac{t_1}{w^2}+\frac{t_2}{w^3}+\frac{t_3}{w^4}+\frac{t_3}{w^5}+\dots=
-s_{-1}+\frac{s_0}w-\frac{s_1}{w^2}+\frac{s_2}{w^3}-\frac{s_3}{w^4}+\frac{s_3}{w^5}-\dots
\end{equation*}
Thus, the connection between $s_j$ and $t_j$ is as follows
\begin{equation*}\label{app.assoc.function.series.self-int.3}
t_j=(-1)^js_j,\qquad j=-1,0,1,2,\ldots
\end{equation*}
Consider the two infinite Hankel matrices
$S=\|s_{i+j}\|_{0}^{\infty}$ and $T=\|t_{i+j}\|_{0}^{\infty}$.

In~\cite{Holtz_Tyaglov} there was proved the following formula.
\begin{theorem}\label{Theorem.self-int.Hankel.matrices.connection}
Minors of the matrices $S$ and $T$ are connected as follows
\begin{equation*}\label{Theorem.self-int.Hankel.matrices.connection.equal}
T\begin{pmatrix}
    i_1 &i_2 &\dots &i_m\\
    j_1 &j_2 &\dots &j_m\\
\end{pmatrix}=
(-1)^{\sum\limits_{k=1}^{m}i_k+\sum\limits_{k=1}^{m}j_k}S
\begin{pmatrix}
    i_1 &i_2 &\dots &i_m\\
    j_1 &j_2 &\dots &j_m\\
\end{pmatrix},\quad m=1,2,\ldots
\end{equation*}
\end{theorem}

From this theorem one can easily obtain the following consequence.
\begin{corol}\label{corol.self-int.Hankel.matrices.connection}
Let rational functions $\Phi(u)$ and $\Psi(u)$ be such that
$\Psi(u)=-\Phi(-u)$. Then the following equalities hold
\begin{equation*}\label{SI_Stab.connection.1}
D_j(\Phi)=D_j(\Psi),\quad j=1,2,\ldots
\end{equation*}
\begin{equation*}\label{SI_Stab.connection.2}
\widehat{D}_j(\Phi)=(-1)^j\widehat{D}_j(\Psi),\quad j=1,2,\ldots
\end{equation*}
\end{corol}

From this corollary and
from~\eqref{Formulae.Gurwitz.1}--\eqref{Formulae.Gurwitz.2} one
obtains
\begin{equation}\label{SI_Stab.connection.Hurwitz.1}
\Delta_{n+1-2j}(q)=\Delta_{n+1-2j}(p),\qquad
j=1,\ldots,\left[\dfrac{n}2\right],
\end{equation}
\begin{equation}\label{SI_Stab.connection.Hurwitz.2}
\Delta_{n-2j}(q)=\displaystyle(-1)^{\left[\tfrac{n+1}2\right]-j}\Delta_{n-2j}(p),\qquad
j=0,1,\ldots,\left[\dfrac{n}2\right],
\end{equation}
where the polynomials $p$ and $q$ are defined
in~\eqref{Th.Stodola.necessary.condition.self-interlacing.poly}
and~\eqref{Th.Stodola.necessary.condition.self-interlacing.poly.proof},
respectively.

Now we are in position to prove analogues of
Theorems~\ref{Th.Lienard.Chipart}--\ref{Theorem.stability.with.difeerentiation}
for self-interlacing polynomials.

\begin{theorem}\label{Th.Lienard.Chipart.for.self-interlacing}
The polynomial $p$ defined
in~\eqref{Th.Stodola.necessary.condition.self-interlacing.poly} is
self-interlacing of type~I if and only if
\begin{equation}\label{antypoly.4}
(-1)^{\left[\tfrac{n+1}2\right]}a_n>0,\,(-1)^{\left[\tfrac{n+1}2\right]-1}a_{n-2}>0,\,(-1)^{\left[\tfrac{n+1}2\right]-2}a_{n-4}>0,\,\dots
\end{equation}
or
\begin{equation}\label{antypoly.5}
(-1)^{\left[\tfrac{n+1}2\right]}a_n>0,\,(-1)^{\left[\tfrac{n+1}2\right]-1}a_{n-1}>0,\,(-1)^{\left[\tfrac{n+1}2\right]-2}a_{n-3}>0,\,\dots
\end{equation}
and one of the following two conditions holds
\begin{itemize}
\item[1)]
\begin{equation}\label{antypoly.6}
\Delta_{n-1}(p)>0,\,\Delta_{n-3}(p)>0,\,\Delta_{n-5}(p)>0,\,\dots;
\end{equation}
\item[2)]
\begin{equation}\label{antypoly.7}
\displaystyle(-1)^{\left[\tfrac{n+1}2\right]}\Delta_{n}(p)>0,\,
\displaystyle(-1)^{\left[\tfrac{n+1}2\right]-1}\Delta_{n-2}(p)>0,\,\displaystyle(-1)^{\left[\tfrac{n+1}2\right]-2}\Delta_{n-4}(p)>0,\,\ldots
\end{equation}
\end{itemize}
\end{theorem}
\begin{proof}
If the polynomial $p$ is self-interlacing of type~I, then by
Theorems~\ref{Theorem.self-interlacing.Hurwitz.criterion}
and~\ref{Th.Stodola.necessary.condition.self-interlacing}, all the
conditions~\eqref{antypoly.4}--\eqref{antypoly.7} hold.

Let the conditions~\eqref{antypoly.4} and~\eqref{antypoly.6} hold.
Consider the polynomial $q$ defined
in~\eqref{Th.Stodola.necessary.condition.self-interlacing.poly.proof}.
From~\eqref{antypoly.4}
and~\eqref{Th.Stodola.necessary.condition.self-interlacing.even.proof}--\eqref{Th.Stodola.necessary.condition.self-interlacing.odd.proof}
it follows that $b_n>0,b_{n-2}>0,b_{n-4}>0,\ldots$, and
from~\eqref{antypoly.6},~\eqref{Formulae.Gurwitz.1}--\eqref{Formulae.Gurwitz.2}
and~\eqref{SI_Stab.connection.Hurwitz.1} we obtain that
$\Delta_{n-1}(q)>0,\Delta_{n-3}(q)>0,\ldots$ Thus, the polynomial
$q$ satisfies the condition $1)$ of
Theorem~\ref{Th.Lienard.Chipart}. Therefore, $q$ is Hurwitz
stable, so $p$ is self-interlacing of type~I according to
Theorem~\ref{Theorem.connection.Hurwitz.self-interlacing}.

Analogously,
using~\eqref{Th.Stodola.necessary.condition.self-interlacing.even.proof}--\eqref{Th.Stodola.necessary.condition.self-interlacing.odd.proof},~\eqref{Formulae.Gurwitz.1}--\eqref{Formulae.Gurwitz.2},~\eqref{SI_Stab.connection.Hurwitz.1}--\eqref{SI_Stab.connection.Hurwitz.2}
and Theorem~\ref{Th.Lienard.Chipart} one can show that the
conditions~\eqref{antypoly.5} and~\eqref{antypoly.6},
or~\eqref{antypoly.4} and~\eqref{antypoly.7},
or~\eqref{antypoly.5} and~\eqref{antypoly.7} imply Hurwitz
stability of the polynomial $q$, so by
Theorem~\ref{Theorem.connection.Hurwitz.self-interlacing}, $p$ is
self-interlacing of type~I.
\end{proof}

The following theorem presents a relation between self-interlacing
polynomials and continued fractions of Stieltjes type.
\begin{theorem}\label{Theorem.self-interlacing.Stieltjes.cont.frac.criteria}
The polynomial $p$ of degree $n$ defined
in~\eqref{Th.Stodola.necessary.condition.self-interlacing.poly} is
self-interlacing of type I if and only if its associated
function~$\Phi$ has the following Stieltjes continued fraction
expansion:
\begin{equation}\label{Stieltjes.fraction.for.self-interlacing}
\Phi(u)=c_0+\dfrac1{c_1u+\cfrac1{c_2+\cfrac1{c_{3}u+\cfrac1{\ddots+\cfrac1{c_{2l-1}u+\cfrac1{c_{2l}}}}}}},\quad\text{with}\quad
(-1)^{i-1}c_i>0,\quad i=1,\ldots,2l,
\end{equation}
where $c_0=0$ if $n$ is even, and $c_0>0$ if $n$ is odd, and $l$
as in~\eqref{floor.poly.degree}.
\end{theorem}
\begin{proof}
In fact, by Theorem~\ref{Theorem.main.self-interlacing}, the
polynomial $p$ is self-interlacing of type~I if and only if the
function $\Phi$ can be represented in the
form~\eqref{Mittag.Leffler.self-interlacing}, where
$-\alpha=c_0=0$ if $n=2l$, and $-\alpha=c_0=\dfrac{a_0}{a_1}<0$ if
$n=2l+1$. Now the~assertion of the theorem follows from
Theorem~\ref{Th.R-function.general.properties} and
Corollary~\ref{corol.R-function.Stieltjes.fractions.nonnegative.poles}.
\end{proof}

Theorems~\ref{Theorem.stability.with.difeerentiation}
and~\ref{Theorem.connection.Hurwitz.self-interlacing} imply the
following theorem.
\begin{theorem}\label{Theorem.SI.with.differentiation}
Let $p$ be a self-interlacing polynomial of type~I of degree
$n\geqslant2$ as
in~\eqref{Th.Stodola.necessary.condition.self-interlacing.poly}.
Then all the polynomials
\begin{equation*}\label{Theorem.SI.with.differentiation.condition}
p_j(z)=\sum\limits_{i=0}^{n-2j}\left[\dfrac{n-i}2\right]\left(\left[\dfrac{n-i}2\right]-1\right)\cdots\left(\left[\dfrac{n-i}2\right]+j-1\right)a_iz^{n-2j-i},\quad
j=1,\ldots,\left[\dfrac{n}2\right]-1,
\end{equation*}
also are self-interlacing of type~I \emph{(}if $j$ is even\emph{)}
or of type~II \emph{(}if $j$ is odd\emph{)}.
\end{theorem}

\vspace{3mm}

Let again $p$ be a self-interlacing polynomial of type~Iof degree
$n$ and let $q$ be its associated Hurwitz stable polynomial, that
is, $p(z)=p_0(z^2)+zp_1(z^2)$ and
$q(z)=(-1)^{\tfrac{n(n+1)}2}[p_0(-z^2)-zp_1(-z^2)]$. We are in a
position to find an interrelation between Hurwitz minors of these
polynomials.

Let $n=2l$. By
Theorem~\ref{Theorem.connection.Hurwitz.self-interlacing}, the
Hurwitz matrix of the polynomial $q$ has the form
\begin{equation*}
\mathcal{H}_n(q)=
\begin{pmatrix}
a_1&-a_3& a_5&-a_7&\dots&0\\
a_0&-a_2& a_4&-a_6&\dots&0\\
0  & a_1&-a_3& a_5&\dots&0\\
0  & a_0&-a_2& a_4&\dots&0\\
\dots&\dots&\dots&\dots&\dots&\dots\\
0&0&0&0&\dots&(-1)^{\left[\frac{n}2\right]}a_n\\
\end{pmatrix}.
\end{equation*}
It is easy to see that the matrix $\mathcal{H}_n(q)$ can be
factorized as follows
\begin{equation}\label{SI_Stab.connection.7}
\mathcal{H}_n(q)=\widetilde{C}_n\mathcal{H}_n(p)\widetilde{E}_n,
\end{equation}
where the $n\times n$ matrices $\widetilde{E}_n$ and
$\widetilde{C}_n$ have the forms
\begin{equation*}
\widetilde{E}_n=
\begin{pmatrix}
    1 &0  &0  & 0 & 0 &\dots\\
    0 & -1 &0  & 0 & 0 &\dots\\
    0 &0  &1 & 0 & 0 &\dots\\
    0 &0  &0  &-1 & 0 &\dots\\
    0 &0  &0  & 0 & 1 &\dots\\
    \vdots&\vdots&\vdots&\vdots&\vdots&\ddots
\end{pmatrix},
\qquad
\widetilde{C}_n=
\begin{pmatrix}
    1 &0  &0  & 0 & 0 &\dots\\
    0 & 1 &0  & 0 & 0 &\dots\\
    0 &0  &-1 & 0 & 0 &\dots\\
    0 &0  &0  &-1 & 0 &\dots\\
    0 &0  &0  & 0 & 1 &\dots\\
    \vdots&\vdots&\vdots&\vdots&\vdots&\ddots
\end{pmatrix}.
\end{equation*}
All non-principal minors of these matrices are equal to zero. So
since $\widetilde{e}_{jj}=(-1)^{j-1}$ and
$\widetilde{c}_{jj}=(-1)^{\tfrac{(j-1)(j-2)}2}$, we have
\begin{equation}\label{Matrix.Unique.2.minors}
\widetilde{C}_n
\begin{pmatrix}
    i_1 &i_2 &\dots &i_m\\
    i_1 &i_2 &\dots &i_m
\end{pmatrix}=
(-1)^{\sum\limits_{k=1}^m\tfrac{(i_k-1)(i_k-2)}2},
\qquad
 \widetilde{E}_n
\begin{pmatrix}
    i_1 &i_2 &\dots &i_m\\
    i_1 &i_2 &\dots &i_m
\end{pmatrix}=(-1)^{\sum\limits_{k=1}^{m}i_k-m},
\end{equation}
where $1\leqslant i_1<i_2<\ldots<i_m\leqslant n$. Thus, the
Cauchy--Binet formula together with~\eqref{Matrix.Unique.2.minors}
and~\eqref{SI_Stab.connection.7} implies, for $n=2l$,

\vspace{3mm}

\begin{equation}\label{SI_Stab.connection.5}
\mathcal{H}_n(q)
\begin{pmatrix}
    i_1 &i_2 &\dots &i_m\\
    j_1 &j_2 &\dots &j_m
\end{pmatrix}=(-1)^{\sum\limits_{k=1}^m\tfrac{(i_k-1)(i_k-2)}2+\sum\limits_{k=1}^mj_k-m}
\mathcal{H}_n(p)
\begin{pmatrix}
    i_1 &i_2 &\dots &i_m\\
    j_1 &j_2 &\dots &j_m
\end{pmatrix}.
\end{equation}
where $1\leqslant
\begin{array}{c}
i_1<i_2<\ldots<i_m\\
j_1<j_2<\ldots<j_m
\end{array}
\leqslant n$,

\vspace{3mm}

Let now $n=2l+1$. Then by
Theorem~\ref{Theorem.connection.Hurwitz.self-interlacing}, the
Hurwitz matrix of the polynomial $q$ has the form
\begin{equation*}
\mathcal{H}_n(q)=
\begin{pmatrix}
-a_1& a_3&-a_5& a_7&\dots&0\\
 a_0&-a_2& a_4&-a_6&\dots&0\\
0  &-a_1& a_3&-a_5&\dots&0\\
0  & a_0&-a_2& a_4&\dots&0\\
\dots&\dots&\dots&\dots&\dots&\dots\\
0&0&0&0&\dots&(-1)^{\left[\frac{n+1}2\right]}a_n\\
\end{pmatrix}.
\end{equation*}
In this case, $\mathcal{H}_n(q)$ also can be factorized:
\begin{equation}\label{SI_Stab.connection.8}
\mathcal{H}_n(q)=\widehat{C}_n\mathcal{H}_n(p)(-\widetilde{E}_n),
\end{equation}
where the $n\times n$ matrix $\widehat{C}_n$ is as follows
\begin{equation*}\label{Matrix.Unique.1}
\widehat{C}_n=
\begin{pmatrix}
    1 &0  &0  & 0 & 0 & 0 &\dots\\
    0 &-1 &0  & 0 & 0 & 0 &\dots\\
    0 &0  &-1 & 0 & 0 & 0 &\dots\\
    0 &0  &0  & 1 & 0 & 0 &\dots\\
    0 &0  &0  & 0 & 1 & 0 &\dots\\
    0 &0  &0  & 0 & 0 &-1 &\dots\\
    \vdots&\vdots&\vdots&\vdots&\vdots&\vdots&\ddots
\end{pmatrix}.
\end{equation*}
All non-principal minors of the matrices $\widehat{C}_n$ and
$-\widetilde{E}_n$ equal zero. The principal minors of these
matrices can be easily calculated:
\begin{equation}\label{Matrix.Unique.1.minors}
\widehat{C}_n
\begin{pmatrix}
    i_1 &i_2 &\dots &i_m\\
    i_1 &i_2 &\dots &i_m
\end{pmatrix}=
(-1)^{\sum\limits_{k=1}^m\tfrac{i_k(i_k-1)}2},
\qquad
\widehat{E}_n\begin{pmatrix}
    i_1 &i_2 &\dots &i_m\\
    i_1 &i_2 &\dots &i_m\\
\end{pmatrix}
=(-1)^{\sum\limits_{k=1}^{m}i_k},
\end{equation}
where $1\leqslant i_1<i_2<\ldots<i_m\leqslant n$. Thus, the
Cauchy--Binet formula together with~\eqref{Matrix.Unique.1.minors}
and~\eqref{SI_Stab.connection.8} implies, for $n=2l+1$,

\vspace{3mm}

\begin{equation}\label{SI_Stab.connection.6}
\mathcal{H}_n(q)
\begin{pmatrix}
    i_1 &i_2 &\dots &i_m\\
    j_1 &j_2 &\dots &j_m
\end{pmatrix}=(-1)^{\sum\limits_{k=1}^m\tfrac{i_k(i_k-1)}2+\sum\limits_{k=1}^mj_k}
\mathcal{H}_n(p)
\begin{pmatrix}
    i_1 &i_2 &\dots &i_m\\
    j_1 &j_2 &\dots &j_m
\end{pmatrix},
\end{equation}
where $1\leqslant
\begin{array}{c}
i_1<i_2<\ldots<i_m\\
j_1<j_2<\ldots<j_m
\end{array}
\leqslant n$.

\vspace{3mm}

Since $q$ is Hurwitz stable, from
Theorem~\ref{Theorem.Hurwitz.stable.Hurwitz.matrix.criteria} and
from the formul\ae~\eqref{SI_Stab.connection.5}
and~\eqref{SI_Stab.connection.6} we obtain the following theorem.
\begin{theorem}
Let $p$ be a self-interlacing polynomial of type~I of degree $n$
and let $\mathcal{H}_n(p)$ be its Hurwitz matrix defined
in~\eqref{HurwitzMatrix}. Then
\begin{itemize}
\item[] for $n=2l$,
\begin{equation*}
(-1)^{\sum\limits_{k=1}^m\tfrac{(i_k-1)(i_k-2)}2+\sum\limits_{k=1}^mj_k-m}
\mathcal{H}_n(p)
\begin{pmatrix}
    i_1 &i_2 &\dots &i_m\\
    j_1 &j_2 &\dots &j_m
\end{pmatrix}\geqslant0;
\end{equation*}
\item[] for $n=2l+1$,
\begin{equation*}
(-1)^{\sum\limits_{k=1}^m\tfrac{i_k(i_k-1)}2+\sum\limits_{k=1}^mj_k}
\mathcal{H}_n(p)
\begin{pmatrix}
    i_1 &i_2 &\dots &i_m\\
    j_1 &j_2 &\dots &j_m
\end{pmatrix}\geqslant0,
\end{equation*}
\end{itemize}
where  $1\leqslant
\begin{array}{c}
i_1<i_2<\ldots<i_m\\
j_1<j_2<\ldots<j_m
\end{array}
\leqslant n$.
\end{theorem}

\vspace{4mm}

\subsection{The second proof of the Hurwitz self-interlacing
criterion}\label{subsection:Second.proof.of.SI.criterion}

\hspace{4mm} Given a polynomial
\begin{equation}\label{SI.second.proof.poly}
p(z)=a_0z^n+a_1z^{n-1}+\dots+a_n,\qquad a_1,\dots,a_n\in\mathbb
R,\ a_0>0,
\end{equation}
we consider the following rational functions
\begin{equation}\label{Rat.func.2}
R(z)=\dfrac{(-1)^np(-z)}{p(z)}\quad\text{and}\quad
F(z)=\dfrac1{R(z)}=\dfrac{p(z)}{(-1)^np(-z)}
\end{equation}
and expand them into their Laurent series at $\infty$:
\begin{equation*}
R(z)=1+\frac{s_0}z+\frac{s_1}{z^2}+\frac{s_2}{z^3}+\dots,\qquad
F(z)=1+\frac{t_0}z+\frac{t_1}{z^2}+\frac{t_2}{z^3}+\dots
\end{equation*}
%
As we mentioned in Section~\ref{subsection:Hankel}, ranks of the
matrices $S=\|s_{j+k}\|_0^{\infty}$ and $T=\|t_{j+k}\|_0^{\infty}$
are equal to the number of poles of the functions $R$ and $F$,
respectively. It is easy to see that rank of each matrix
equals~$n$ if the polynomials $p(z)$ and $p(-z)$ have no common
zeroes. In the rest of this section we consider only such
polynomials, so in this section, ranks of the matrices $S$ and $T$
always equal $n$. Denoting $$q(z)=(-1)^np(-z),$$ we obtain
$R=\dfrac{q}{p}$ and $F=\dfrac{p}{q}$.
\begin{lemma}\label{lem.Determinants.relations.1}
For the functions $R$ and $F$ defined in~\eqref{Rat.func.2}, the
following formul\ae~are valid:
\begin{equation}\label{Determinants.relations.1}
\nabla_{2j}(p,q)=a_0^{2j}D_j(R)=(-1)^{\frac{j(j+1)}2}2^ja_0\Delta_{j-1}(p)\Delta_j(p),\quad
j=1,2,\ldots
\end{equation}
and
\begin{equation}\label{Determinants.relations.2}
\nabla_{2j}(q,p)=a_0^{2j}D_j(F)=(-1)^{\frac{j(j-1)}2}2^ja_0\Delta_{j-1}(p)\Delta_j(p),\quad
j=1,2,\ldots,
\end{equation}
where the determinants $\nabla_{2j}(p,q)$ are the leading
principal minors of the matrix $\mathcal{H}_{2n+1}(p,q)$ of order
$2j$, and $\Delta_j(p)$ are defined in~\eqref{delta},
$\Delta_0(p)\equiv1$.
\end{lemma}
\begin{proof}
We prove the formul\ae~\eqref{Determinants.relations.1}. The
formul\ae~\eqref{Determinants.relations.2} follow
from~\eqref{Determinants.relations.1} since
$\nabla_{2j}(q,p)=(-1)^j\nabla_{2j}(p,q)$.

Assuming $a_j=0$ for $j>n$, we have
\begin{equation*}
\nabla_{2j}(p,q)=
\begin{vmatrix}
    a_0 &a_1 &a_2 &a_3 &\dots &a_{j-1} & a_{j}  &\dots &a_{2j-2} & a_{2j-1}\\
    a_0 &-a_1&a_2 &-a_3&\dots &-a_{j-1}& a_{j}  &\dots &a_{2j-2} &-a_{2j-1}\\
     0  &a_0 &a_1 &a_2 &\dots &a_{j-2} & a_{j-1}&\dots &a_{2j-3} &a_{2j-2}\\
     0  &a_0 &-a_1&a_2 &\dots &a_{j-2} &-a_{j-1}&\dots &-a_{2j-3}&a_{2j-2}\\
     0  & 0  &a_0 &a_1 &\dots &a_{j-3} & a_{j-2}&\dots &a_{2j-4} & a_{2j-3}\\
     0  & 0  &a_0 &-a_1&\dots &-a_{j-3}& a_{j-2}&\dots &a_{2j-4} &-a_{2j-3}\\
    \dots&\dots&\dots&\dots&\dots&\dots&\dots&\dots&\dots&\dots\\
     0  &  0 &  0 &  0 &\dots &a_{0} & a_{1}&\dots &a_{j-1} &a_{j}\\
     0  &  0 &  0 &  0 &\dots &a_{0} &-a_{1}&\dots &-a_{j-1}&a_{j}\\
\end{vmatrix}=
\end{equation*}
\begin{equation*}
=2^j\begin{vmatrix}
    a_0 &a_1 &a_2 &a_3 &\dots &a_{j-1} & a_{j}  &\dots &a_{2j-2} & a_{2j-1}\\
    a_0 &  0 &a_2 &  0 &\dots &   0    & a_{j}  &\dots &a_{2j-2} &    0    \\
     0  &a_0 &a_1 &a_2 &\dots &a_{j-2} & a_{j-1}&\dots &a_{2j-3} &a_{2j-2}\\
     0  &a_0 &  0 &  0 &\dots &a_{j-2} &    0   &\dots &    0    &a_{2j-2}\\
     0  & 0  &a_0 &a_1 &\dots &a_{j-3} & a_{j-2}&\dots &a_{2j-4} & a_{2j-3}\\
     0  & 0  &a_0 &  0 &\dots &    0   & a_{j-2}&\dots &a_{2j-4} &    0    \\
    \dots&\dots&\dots&\dots&\dots&\dots&\dots&\dots&\dots&\dots\\
     0  &  0 &  0 &  0 &\dots &a_{0} & a_{1}&\dots &a_{j-1} &a_{j}\\
     0  &  0 &  0 &  0 &\dots &a_{0} &   0  &\dots &    0   &a_{j}\\
\end{vmatrix}=
\end{equation*}
\begin{equation*}
=(-2)^j\begin{vmatrix}
    a_0 &  0 &a_2 &  0 &\dots &   0    & a_{j}  &\dots &a_{2j-2} &    0    \\
     0  &a_1 &  0 &a_3 &\dots &a_{j-1} &   0    &\dots &    0    & a_{2j-1}\\
     0  &a_0 &  0 &a_2 &\dots &a_{j-2} &   0    &\dots &    0    &a_{2j-2}\\
     0  & 0  &a_1 &  0 &\dots &   0    & a_{j-1}&\dots &a_{2j-3} &    0    \\
     0  & 0  &a_0 &  0 &\dots &   0    & a_{j-2}&\dots &a_{2j-4} &    0    \\
     0  & 0  &  0 &a_1 &\dots &a_{j-3} &   0    &\dots &    0    &a_{2j-3}\\
    \dots&\dots&\dots&\dots&\dots&\dots&\dots&\dots&\dots&\dots\\
     0  &  0 &  0 &  0 &\dots &a_{0} &   0  &\dots &    0   &a_{j}\\
     0  &  0 &  0 &  0 &\dots &   0  & a_{1}&\dots &a_{j-1} &  0  \\
\end{vmatrix}=
\end{equation*}
\begin{equation*}
=(-2)^ja_0\begin{vmatrix}
    a_1 &  0 &a_3 &  0  &\dots &a_{j-1} &   0    &\dots &    0    & a_{2j-1}\\
    a_0 &  0 &a_2 &  0  &\dots &a_{j-2} &   0    &\dots &    0    &a_{2j-2}\\
      0 &  0 &a_1 &  0  &\dots &a_{j-3} &   0    &\dots &    0    &a_{2j-3}\\
      0 &  0 &a_0 &  0  &\dots &a_{j-4} &   0    &\dots &    0    &a_{2j-4}\\
    \dots&\dots&\dots&\dots&\dots&\dots&\dots&\dots&\dots&\dots\\
      0 &  0 &  0 &  0  &\dots &a_{1}   &   0    &\dots &    0    &a_{j+1}\\
      0 &  0 &  0 &  0  &\dots &a_{0}   &   0    &\dots &    0    &a_{j}\\
      0 &a_1 &  0 &a_3  &\dots &   0    & a_{j-1}&\dots &a_{2j-3} &    0    \\
      0 &a_0 &  0 &a_2  &\dots &   0    & a_{j-2}&\dots &a_{2j-4} &    0    \\
    \dots&\dots&\dots&\dots&\dots&\dots&\dots&\dots&\dots&\dots\\
      0 &  0 &  0 &  0  &\dots &   0    & a_{2}  &\dots &a_{j}   &  0  \\
      0 &  0 &  0 &  0  &\dots &   0    & a_{1}  &\dots &a_{j-1} &  0  \\
\end{vmatrix}=
\end{equation*}
\begin{equation*}
=(-2)^ja_0(-1)^{\frac{j(j-1)}2}\begin{vmatrix}
    a_1 &a_3 &a_5 &\dots &a_{2j-1} &   0    &   0    &\dots &    0    &   0\\
    a_0 &a_2 &a_4 &\dots &a_{2j-2} &   0    &   0    &\dots &    0    &   0\\
      0 &a_1 &a_3 &\dots &a_{2j-3} &   0    &   0    &\dots &    0    &   0\\
      0 &a_0 &a_2 &\dots &a_{2j-4} &   0    &   0    &\dots &    0    &   0\\
    \dots&\dots&\dots&\dots&\dots&\dots&\dots&\dots&\dots&\dots\\
      0 &  0 &  0 &\dots &a_{j+1} &   0    &   0    &\dots &    0    &   0\\
      0 &  0 &  0 &\dots &a_{j}   &   0    &   0    &\dots &    0    &   0\\
      0 &  0 &  0 &\dots &   0    & a_1    & a_{3}  &\dots &a_{2j-5} &a_{2j-3}\\
      0 &  0 &  0 &\dots &   0    & a_0    & a_{2}  &\dots &a_{2j-6} &a_{2j-4}\\
    \dots&\dots&\dots&\dots&\dots&\dots&\dots&\dots&\dots&\dots\\
      0 &  0 &  0 &\dots &  0     &   0    &    0   &\dots &a_{j-2} &a_{j}   \\
      0 &  0 &  0 &\dots &  0     &   0    &    0   &\dots &a_{j-3} &a_{j-1}\\
\end{vmatrix}=
\end{equation*}
\begin{equation*}
=(-1)^{\frac{j(j+1)}2}2^ja_0\Delta_{j-1}(p)\Delta_j(p).
\end{equation*}
\end{proof}

Using this lemma it is easy to prove the equivalence of the
conditions $1)$ and $2)$ of
Theorem~\ref{Theorem.self-interlacing.Hurwitz.criterion}.
\begin{theorem}\label{Theorem.SI.Hurwitz.criterion.2}
The polynomial $p$ defined in~\eqref{SI.second.proof.poly} is
self-interlacing of type~I if and only if
\begin{equation}\label{Hurvitz.det.noneq.self-interlacing.1111}
\Delta_{n-1}(p)>0,\ \Delta_{n-3}(p)>0,\dots,
\end{equation}
\begin{equation}\label{Hurvitz.det.noneq.self-interlacing.22222}
\displaystyle(-1)^{\left[\tfrac{n+1}2\right]}\Delta_{n}(p)>0,
\displaystyle(-1)^{\left[\tfrac{n+1}2\right]-1}\Delta_{n-2}(p)>0,\ldots
\end{equation}
where the Hurwitz minors $\Delta_i(p)$ are defined
in~\eqref{delta}.
\end{theorem}
\begin{proof}
In the proof of Theorem~\ref{Theorem.main.self-interlacing}, it
was established that if $p$ is self-interlacing of type~I, then
the function $-\dfrac{p(-z)}{p(z)}$ is an \textit{R}-function with
exactly $n$ poles.

\vspace{2mm}

Let $n=2l+1$. Then the function $R$ defined in~\eqref{Rat.func.2}
is an \textit{R}-function. By
Theorem~\ref{Th.R-function.general.properties}, we have $D_j(R)>0$
for $j=1,\ldots,n$. These inequalities together with the
formul\ae~\eqref{Determinants.relations.1} imply
\begin{equation}\label{Second.proof.relations}
(-1)^{\frac{j(j+1)}2}\Delta_{j-1}(p)\Delta_j(p)>0,\qquad
j=1,\ldots,n.
\end{equation}
Multiplying the inequalities~\eqref{Second.proof.relations} for
$j=2m$ and $j=2m-1$, we obtain
\begin{equation*}
\Delta_{2m-1}^2\Delta_{2m}\Delta_{2m-2}>0.
\end{equation*}
Consequently, the minors $\Delta_{2i}(p)$, $i=1,\ldots,l$, are
positive, so the
inequalities~\eqref{Hurvitz.det.noneq.self-interlacing.1111} are
proved for odd~$n$.

If we multiply the inequalities~\eqref{Second.proof.relations} for
$j=2m$ and $j=2m+1$, we get
\begin{equation*}
-\Delta_{2m}^2(p)\Delta_{2m-1}(p)\Delta_{2m+1}(p)>0.
\end{equation*}
These inequalities
imply~\eqref{Hurvitz.det.noneq.self-interlacing.22222} for
odd~$n$.

The converse assertion can be proved in the same way. That is, the
inequalities~\eqref{Hurvitz.det.noneq.self-interlacing.1111}--\eqref{Hurvitz.det.noneq.self-interlacing.22222}
imply the inequalities~\eqref{Second.proof.relations} which in
turn imply the inequalities $D_j(R)>0$, $j=1,\ldots,n$, according
to~\eqref{Determinants.relations.1}. By
Theorem~\ref{Th.R-function.general.properties}, the function
$-\dfrac{p(-z)}{p(z)}$ is an \textit{R}-function, so $p$ is a
self-interlacing polynomial of type~I (see the proof of
Theorem~\ref{Theorem.main.self-interlacing}).

\vspace{2mm}

If $n$ is even, the
inequalities~\eqref{Hurvitz.det.noneq.self-interlacing.1111}--\eqref{Hurvitz.det.noneq.self-interlacing.22222}
can be proved analogously using the function~$F$ defined
in~\eqref{Rat.func.2} instead of the function $R$.
\end{proof}

From Lemma~\ref{lem.Determinants.relations.1} and from Hurwitz's
stability criterion (see
Theorem~\ref{Theorem.Hurwitz.stable.Hurwitz.matrix.criteria}) it
follows one more stability criterion.

\begin{theorem}\label{Theorem.new.stability.criterion}
$\quad$
\begin{itemize}
\item[1)] The polynomial $p$ is Hurwitz stable if and only if for
the function $R$ defined in~\eqref{Rat.func.2}, the following
inequalities hold
\begin{equation*}
(-1)^{\frac{j(j+1)}2}D_j(R)>0,\quad j=1,2,\ldots,n.
\end{equation*}
\item[2)] The polynomial $p$ is Hurwitz stable if and only if for
the function $F$ defined in~\eqref{Rat.func.2}, the following
inequalities hold
\begin{equation*}
(-1)^{\frac{j(j-1)}2}D_j(F)>0,\quad j=1,2,\ldots,n.
\end{equation*}
\end{itemize}
\end{theorem}

Finally, we establish a relationships between the Hurwitz minors
$\Delta_j(p)$ of a given polynomial $p$ and the Hankel minors
$\widehat{D}_j(R)$ and $\widehat{D}_j(F)$ (for the definition of
the minors $\widehat{D}_j$ see~\eqref{Hankel.determinants.2}). To
do this, consider the functions
\begin{equation}\label{s.R(z)}
zR(z)=z+\dfrac{h_1(z)}{p(z)}=z+\dfrac{-2a_1z^n-2a_3z^{n-2}-2a_5z^{n-4}-\dots}{a_0z^n+a_1z^{n-1}+\dots+a_n}=z+s_0+\frac{s_1}{z}+\frac{s_2}{z^2}+\frac{s_3}{z^3}+\dots
\end{equation}
\begin{equation}\label{s.F(z)}
zF(z)=z+\dfrac{h_2(z)}{(-1)^np(-z)}=z+\dfrac{2a_1z^n+2a_3z^{n-2}+2a_5z^{n-4}+\dots}{a_0z^n-a_1z^{n-1}+\dots+(-1)^na_n}=z+t_0+\frac{t_1}{z}+\frac{t_2}{z^2}+\frac{t_3}{z^3}+\dots
\end{equation}
where $s_0=t_0=1$, and note that $\widehat{D}_j(R)=D_j(zR)$ and
$\widehat{D}_j(F)=D_j(zF)$. This allows us to establish
the~following lemma.
\begin{lemma}\label{lem.Determinants.relations.2}
Given the polynomial $p$, for the functions $R$ and $F$ defined
in~\eqref{Rat.func.2}, the following relationships hold
\begin{equation}\label{Determinants.relations.3}
\nabla_{2j}(p,h_1)=a_0^{2j}\widehat{D}_j(R)=(-1)^{\frac{j(j-1)}2}2^j\Delta_{j}^2(p),\quad
j=1,2,\ldots,
\end{equation}
and
\begin{equation}\label{Determinants.relations.4}
\nabla_{2j}(p,h_2)=a_0^{2j}\widehat{D}_j(F)=(-1)^{\frac{j(j-1)}2}2^j\Delta_{j}^2(p),\quad
j=1,2,\ldots,
\end{equation}
where the determinants\footnote{The polynomials $h_1$ and $h_2$
are defined in~\eqref{s.R(z)} and~\eqref{s.F(z)}, respectively.}
$\nabla_{2j}(p,h_1)$ and $\nabla_{2j}(p,h_2)$ are the leading
principal minors  of order $2j$ of the~matrices
$\mathcal{H}_{2n+1}(p,h_1)$ and $\mathcal{H}_{2n+1}(p,h_2)$
respectively, and $\Delta_j(p)$ are the Hurwitz minors of the
polynomial~$p$, $\Delta_0(p)\equiv1$.
\end{lemma}
\begin{proof}
Assuming $a_j=0$ for $j>n$, from~\eqref{s.R(z)} we obtain
\begin{equation*}
\nabla_{2j}(p,h_1)=(-2)^j
\begin{vmatrix}
    a_0 &a_1 &a_2 &a_3 &\dots &a_{j-1} & a_{j}  &\dots &a_{2j-2} & a_{2j-1}\\
    a_1 &  0 &a_3 &  0 &\dots &    0   & a_{j+1}&\dots &a_{2j-1} &    0    \\
     0  &a_0 &a_1 &a_2 &\dots &a_{j-2} & a_{j-1}&\dots &a_{2j-3} &a_{2j-2}\\
     0  &a_1 &  0 &a_3 &\dots &a_{j-1} &    0   &\dots &    0    &a_{2j-1}\\
     0  & 0  &a_0 &a_1 &\dots &a_{j-3} & a_{j-2}&\dots &a_{2j-4} & a_{2j-3}\\
     0  & 0  &a_1 &  0 &\dots &    0   & a_{j-1}&\dots &a_{2j-3} &    0    \\
    \dots&\dots&\dots&\dots&\dots&\dots&\dots&\dots&\dots&\dots\\
     0  &  0 &  0 &  0 &\dots &a_{0}   & a_{1}  &\dots &a_{j-1}  &a_{j}\\
     0  &  0 &  0 &  0 &\dots &a_{1}   &     0  &\dots &    0    &a_{j+1}\\
\end{vmatrix}=
\end{equation*}
\begin{equation*}
=2^j\begin{vmatrix}
    a_1 &  0 &a_3 &  0  &\dots &   0    &a_{j+1} &\dots & a_{2j-1}&   0    \\
    a_0 &  0 &a_2 &  0  &\dots &   0    &a_{j}   &\dots & a_{2j-2}&   0    \\
      0 &a_1 &  0 & a_3 &\dots &a_{j-1} &   0    &\dots &    0    &a_{2j-1}\\
      0 &a_0 &  0 & a_2 &\dots &a_{j-2} &   0    &\dots &    0    &a_{2j-2}\\
      0 &  0 &a_1 &  0  &\dots &   0    &a_{j-1} &\dots &a_{2j-3} &   0    \\
      0 &  0 &a_0 &  0  &\dots &   0    &a_{j-2} &\dots &a_{2j-4} &   0    \\
    \dots&\dots&\dots&\dots&\dots&\dots&\dots&\dots&\dots&\dots\\
      0 &  0 &  0 &  0  &\dots &   0    & a_{2}  &\dots &a_{j}    &   0    \\
      0 &  0 &  0 &  0  &\dots & a_{1}  &   0    &\dots &    0    &a_{j+1} \\
      0 &  0 &  0 &  0  &\dots & a_{0}  &   a_1    &\dots &a_{j-1}  &a_{j}   \\
\end{vmatrix}=
\end{equation*}
\begin{equation*}
=2^j(-1)^{\frac{j(j-1)}2}\begin{vmatrix}
    a_1 &a_3 &a_5 &\dots &a_{2j-1} &   0    &   0    &\dots &    0    &   0\\
    a_0 &a_2 &a_4 &\dots &a_{2j-2} &   0    &   0    &\dots &    0    &   0\\
      0 &a_1 &a_3 &\dots &a_{2j-3} &   0    &   0    &\dots &    0    &   0\\
      0 &a_0 &a_2 &\dots &a_{2j-4} &   0    &   0    &\dots &    0    &   0\\
    \dots&\dots&\dots&\dots&\dots&\dots&\dots&\dots&\dots&\dots\\
      0 &  0 &  0 &\dots &a_{j+1} &   0    &   0    &\dots &    0    &   0\\
      0 &  0 &  0 &\dots &a_{j}   &   0    &   0    &\dots &    0    &   0\\
      0 &  0 &  0 &\dots &   0    & a_1    & a_{3}  &\dots &a_{2j-3} &a_{2j-1}\\
      0 &  0 &  0 &\dots &   0    & a_0    & a_{2}  &\dots &a_{2j-4} &a_{2j-2}\\
    \dots&\dots&\dots&\dots&\dots&\dots&\dots&\dots&\dots&\dots\\
      0 &  0 &  0 &\dots &  0     &   0    &    0   &\dots &a_{j-1}  &a_{j+1}   \\
      0 &  0 &  0 &\dots &a_{j-1} &   0    &    0   &\dots &a_{j-2}  &a_{j}\\
\end{vmatrix}=
\end{equation*}
\begin{equation*}
=(-1)^{\frac{j(j-1)}2}2^j\Delta_j^2(p).
\end{equation*}
The relationships~\eqref{Determinants.relations.4} can be proved
analogously.
\end{proof}

Thus, we obtain that for any real polynomial $p$ of degree $n$
with nonzero Hurwitz minors and, in particular, for Hurwitz stable
and for self-interlacing polynomials, the following inequalities
hold
\begin{equation*}
(-1)^{\frac{j(j-1)}2}\widehat{D}_j(R)=(-1)^{\frac{j(j-1)}2}\widehat{D}_j(F)>0,\quad
j=1,2,\ldots,n.
\end{equation*}

Last, let us note that the
formul\ae~\eqref{Determinants.relations.1}--\eqref{Determinants.relations.2}
and~\eqref{Determinants.relations.3}--\eqref{Determinants.relations.4}
can be obtained (overcoming certain difficulties) from some
theorems of the book~\cite{Wall}. But it was more simple to deduce
them directly as we did in
Lemmata~\ref{lem.Determinants.relations.1}
and~\ref{lem.Determinants.relations.2}.

\subsection{Almost and quasi- self-interlacing
polynomials}\label{subsection:Alomst.SI}

\hspace{4mm} In this section, we describe polynomials, which are
dual (in the sense of
Theorem~\ref{Theorem.connection.Hurwitz.self-interlacing}) to the
quasi-stable polynomials. Because of the mentioned duality we just
give the definition of these polynomials.

\begin{definition}\label{def.Alomst.SI.poly}
A polynomial $p(z)=p_0(z^2)+zp_1(z^2)$ of degree $n$ is called
\textit{quasi-self-interlacing} of type~I (of type~II) with
\textit{degeneracy index} $m$, $0\leqslant m\leqslant n$, if the
polynomial $p_0(-z^2)-zp_1(-z^2)$ (resp. the polynomial
$p_0(-z^2)+zp_1(-z^2)$) is quasi-stable with degeneracy index $m$.
\end{definition}
In other words, the polynomial $p$ is quasi-self-interlacing if
non-common zeroes f the polynomials $p(z)$ and $p(-z)$ are real,
simple and interlacing.

Obviously, all results regarding quasi-stable polynomials can be
easily reformulated for the quasi-stable polynomials. So we leave
such reformulation to the reader.

\vspace{3mm}

In Section~\ref{section:Generalized.Hurwitz} we use the
quasi-stable polynomials with degeneracy index $1$, that is, the
polynomials of the form $p(z)=zq(z)$, where $q(z)$ is a
self-interlacing polynomial. Such polynomials we call
\textit{almost self-interlacing} polynomials. Note that if
$p(z)=p_0(z^2)+zp_1(z^2)$ is almost self-interlacing, then the
function $\Phi=\dfrac{p_1}{p_0}$ is an \textit{R}-function with
all poles positive except one, which is zero.

\subsection{''Strange'' polynomials}\label{subsection:strange.poly}

\hspace{4mm} This section is devoted to short description of some
numeric results on the distribution of zeroes of some polynomials
closely connected to the Hurwitz stable and self-interlacing
polynomials. We call these polynomials ''strange'' because of
their curious zero location.

\vspace{2mm}

Let $p(z)=p_0(z^2)+zp_1(z^2)$ be a Hurwitz stable polynomial of
degree $n$. Consider the polynomial $q(z)=p_0(-z^2)+zp_1(z^2)$. It
is evident that $q(0)=p(0)\neq0$. Calculations with Maple Software
showed that the polynomial $q$ has exactly
$\left[\dfrac{n+1}2\right]$ simple zeroes in the \textit{open}
right half-plane and exactly $\left[\dfrac{n}2\right]$ simple
zeroes in the \textit{open} left half-plane, but it has no zeroes
on the imaginary axis. At least one zero of $q$ is nonreal.

Let us denote by $\mu_1,\mu_2,\ldots$ the \textit{distinct}
absolute values of the zeroes of the polynomial $q$ lying in
the~open left half-plane such that
\begin{equation*}
0<\mu_1<\mu_2<\mu_3<\ldots
\end{equation*}
And denote by $\lambda_1,\lambda_2,\ldots$ the \textit{distinct}
absolute values of the zeroes of the polynomial $q$ lying in
the~open right half-plane such that
\begin{equation*}
0<\lambda_1<\lambda_2<\lambda_3<\ldots
\end{equation*}
So numerical experiments showed that
\begin{equation*}
0<\lambda_1<\mu_1<\lambda_2<\mu_2<\lambda_3<\ldots
\end{equation*}
The polynomial $p_0(z^2)+zp_1(-z^2)$ possesses similar properties
according to calculations.

However, it is not known for sure if all the polynomials of the
type $p_0(-z^2)+zp_1(z^2)$, where the polynomial
$p_0(z^2)+zp_1(z^2)$ is Hurwitz stable, have the zero location
described above.

\subsection{Matrices with self-interlacing
spectrum}\label{subsection:matrices.SI}

\hspace{4mm} In this section, we consider some classes of
\textit{real} matrices with self-interlacing spectrum and develop
a method of constructing such kind of matrices from a given
totally positive matrix. But at first, we recall some definitions
and statements from the book~\cite{KreinGantmaher}.

\begin{definition}[\cite{KreinGantmaher}]\label{def.sign definite.matrix}
A square matrix $A=\|a_{ij}\|_1^n$ is called \textit{sign
definite} of class $n$ if for any $k\leqslant n$, all the non-zero
minors of order $k$ have the same sign $\varepsilon_k$. The
sequence $\{\varepsilon_1,\varepsilon_2,\ldots,\varepsilon_n\}$ is
called the \textit{signature sequence} of the matrix $A$.

A sign definite matrix of class $n$ is called \textit{strictly
sign definite} of class $n$ if all its minors are different from
zero.
\end{definition}
\begin{definition}[\cite{KreinGantmaher}]
A square sign definite matrix $A=\|a_{ij}\|_1^n$ of class $n$ is
called the \textit{matrix of class}~$n^{+}$ if some its power is a
strictly sign definite matrix of class $n$.
\end{definition}

Note that a sign definite (strictly sign definite) matrix of class
$n$ with the signature sequence
$\varepsilon_1=\varepsilon_2=\ldots=\varepsilon_n=1$ is totally
nonnegative (\textit{strictly totally positive}). Also a sign
definite matrix of class $n^{+}$ with the signature sequence
$\varepsilon_1=\varepsilon_2=\ldots=\varepsilon_n=1$ is an
\textit{oscillating} matrix (see~\cite{KreinGantmaher}), that is,
a totally nonnegative matrix whose certain power is strictly
totally positive. It is clear from the Binet-Cauchy formula that
the square of a sign definite matrix is totally nonnegative.

In~\cite{KreinGantmaher} it was established the following lemma.
\begin{lemma}\label{lem.SI_spectra.matrix}
Let the matrix $A=\|a_{ij}\|_1^n$ be totally nonnegative. Then the
matrices $B=\|a_{n-i+1,j}\|_1^n$ and $C=\|a_{i,n-j+1}\|_1^n$ are
sign definite of class $n$. Moreover, the signature sequence of
the matrices $B$ and $C$ is as follows:
\begin{equation}\label{SI_matrix.eps}
\varepsilon_k=(-1)^{\tfrac{k(k-1)}2},\quad k=1,2,\ldots,n.
\end{equation}
\end{lemma}

Note that the matrices $B$ and $C$ can be represented as follows
\begin{equation*}
B=JA,\qquad\text{and}\qquad C=AJ,
\end{equation*}
where
\begin{equation}\label{Matrice.J}
J=
\begin{pmatrix}
    0 & 0 &\dots& 0 & 1 \\
    0 & 0 &\dots& 1 & 0 \\
    \vdots&\vdots&\cdot&\vdots&\vdots\\
    0 & 1 &\dots& 0 & 0 \\
    1 & 0 &\dots& 0 & 0 \\
\end{pmatrix}.
\end{equation}

It is easy to see that the matrix $J$ is sign definite of class
$n$ (but not of class $n^{+}$) with the signature sequence of the
form~\eqref{SI_matrix.eps}. So by the Binet-Cauchy
formula~\cite{Gantmakher.1} we obtain the following statement.
\begin{theorem}\label{teorem.SI_oscill.spectra.matrix}
The matrix $A=\|a_{ij}\|_1^n$ is totally nonnegative if and only
if the matrix $JA$ (or the matrix~$AJ$) is sign definite of class
$n$ with the signature sequence~\eqref{SI_matrix.eps}. The matrix
$J$ is defined in~\eqref{Matrice.J}.
\end{theorem}

Obviously, the converse statement is also true.
\begin{theorem}
The matrix $A=\|a_{ij}\|_1^n$ is a sign definite matrix of class
$n$ with the signature sequence~\eqref{SI_matrix.eps} if and only
if the matrix $JA$ (or the matrix~$AJ$) is totally nonnegative.
\end{theorem}


In the sequel, we need the following two theorems established in
the book~\cite{KreinGantmaher}.
\begin{theorem}\label{teorem.signreg.matrix}
Let the matrix $A=\|a_{ij}\|_1^n$ be a sign definite of class
$n^{+}$ with the signature sequence $\varepsilon_k$,
$k=1,2,\ldots,n$. Then all the eigenvalues $\lambda_k$,
$k=1,2,\ldots,n$, of the matrix $A$ are nonzero real and simple,
and if
\begin{equation}\label{SI_matrix.modules}
|\lambda_1|>|\lambda_2|>\ldots>|\lambda_n|>0,
\end{equation}
then
\begin{equation}\label{SI_matrix.signs.eigvals}
\emph{sign}\,\lambda_k=\dfrac{\varepsilon_k}{\varepsilon_{k-1}},\qquad
k=1,2,\ldots,n,\,\varepsilon_0=1.
\end{equation}
\end{theorem}
\begin{theorem}\label{teorem.oscill.matrix}
A totally nonnegative matrix $A=\|a_{ij}\|_1^n$ is oscillating if
and only if $A$ is nonsingular and the following inequalities hold
\begin{equation*}
a_{j,j+1}>0,\qquad\text{and}\qquad a_{j+1,j}>0,\qquad
j=1,2,\ldots,n-1.
\end{equation*}
\end{theorem}

\vspace{2mm}

Let us introduce the following definition.
\begin{definition}
A matrix $A$ is said to have the \textit{self-interlacing
spectrum} if its eigenvalues are real and simple and satisfy the
following inequalities
\begin{equation}\label{SI_matrix.real.spectra}
\lambda_1>-\lambda_2>\lambda_3>\ldots>(-1)^{n-1}\lambda_n>0,
\end{equation}
or
\begin{equation}\label{SI_matrix.real.spectra.2}
-\lambda_1>\lambda_2>-\lambda_3>\ldots>(-1)^n\lambda_n>0.
\end{equation}
\end{definition}

Now we are in a position to complement
Lemma~\ref{lem.SI_spectra.matrix}.
\begin{theorem}\label{teorem.oscill_signreg.matrix.1}
Let all the entries of a nonsingular matrix $A=\|a_{ij}\|_1^n$ be
nonnegative\footnote{The matrices with nonnegative entries are
usually called \textit{nonnegative} matrices.} and for each $i$,
$i=1,2,\ldots,n-1$, there exist number $r_1$ and $r_2$,
$1\leqslant r_1,r_2\leqslant n$, such that
\begin{equation}\label{matrix.entries.1}
a_{n-i,r_1}\cdot a_{n+1-r_1,i}>0, a_{n+1-i,r_2}\cdot
a_{n+1-r_2,i+1}>0
\end{equation}
%
%
\begin{equation}\label{matrix.entries.2}
(\text{or}\qquad a_{i,n+1-r_1}\cdot a_{r_1,n-i}>0,
a_{i+1,n+1-r_2}\cdot a_{r_2,n+1-i}>0).
\end{equation}
The matrix $A$ is totally nonnegative if and only if the matrix
$B=JA=\|b_{ij}\|_1^n$ \emph{(}or, respectively, the matrix
$C=AJ$\emph{)} is sign definite of class~$n^{+}$ with the
signature sequence defined in~\eqref{SI_matrix.eps}. Moreover, the
matrix~$B$ \emph{(}or, respectively, the matrix $C$\emph{)}
possesses a self-interlacing spectrum of the
form~\eqref{SI_matrix.real.spectra}.
\end{theorem}
\begin{proof}
We prove the theorem in the case when the
condition~\eqref{matrix.entries.1} holds. The case of the
condition~\eqref{matrix.entries.2} can be established analogously.

Let $A$ be a nonsingular totally nonnegative matrix and the
condition~\eqref{matrix.entries.1} holds. From
Theorem~\ref{teorem.SI_oscill.spectra.matrix} it follows that the
matrix $B=JA$ is sign definite of class $n$ with the signature
sequence~\eqref{SI_matrix.eps}. In order to the matrix be be sign
definite of class $n^{+}$ it is necessary and sufficient that a
certain power of this matrix $B$ be strictly sign definite of
class~$n$. Since the entries of the matrix $J$ have the form
\begin{equation*}
(J)_{ij}=\begin{cases}
&1,\qquad i=n+1-j;\\
&0,\qquad i\neq n+1-j;
\end{cases}
\end{equation*}
the entries of the matrix $B$ can be represented as follows:
$$
b_{ij}=\sum\limits_{k=1}^{n}(J)_{ik}a_{kj}=a_{n+1-i,j}.
$$

Consider the totally nonnegative matrix $B^2$. Its entries have
the form
\begin{equation*}
(B^2)_{ij}=\sum\limits_{k=1}^{n}b_{ik}b_{kj}=\sum\limits_{k=1}^{n}a_{n+1-i,k}a_{n+1-k,j}.
\end{equation*}
From these formul\ae~and from~\eqref{matrix.entries.1} it follows
that all the entries of the matrix $B^2$ above and under the main
diagonal are positive, that is, $(B^2)_{i,i+1}>0$ ¨
$(B^2)_{i+1,i}>0$, $i=1,2,\ldots,n-1$. By
Theorem~\ref{teorem.oscill.matrix}, $B^2$ is an oscillating
matrix. According to the definition of oscillating matrices, a
certain power of $B^2$ is strictly totally positive. Thus, we
proved that a certain power of the matrix $B$ is strictly sign
definite, so $B$ is a sign definite matrix of class $n^{+}$ with
the signature sequence~\eqref{SI_matrix.eps} according to
Lemma~\ref{lem.SI_spectra.matrix}. By
Theorem~\ref{teorem.signreg.matrix}, all eigenvalues of the matrix
$B$ are nonzero real and simple. Moreover, if we enumerate the
eigenvalues in order of decreasing absolute values as
in~\eqref{SI_matrix.modules}, then
from~\eqref{SI_matrix.signs.eigvals} and~\eqref{SI_matrix.eps} we
obtain that the spectrum of $B$ is of the
form~\eqref{SI_matrix.real.spectra}.

The converse assertion of the theorem follows from
Theorem~\ref{teorem.SI_oscill.spectra.matrix}.
\end{proof}

\begin{remark}
For $n=2l+1$ (for $n=2l$), the characteristic polynomial of the
matrix $B$ described in
Theorem~\ref{teorem.oscill_signreg.matrix.1} is self-interlacing
of type~I (of type II).
\end{remark}
\begin{remark}
If the matrix $-A$ is totally nonnegative and the
conditions~\eqref{matrix.entries.1} (or the
conditions~\eqref{matrix.entries.2}) hold, then the matrix $B=JA$
(or, respectively, $C=AJ$) has a spectrum of the
form~\eqref{SI_matrix.real.spectra.2}.
\end{remark}
\begin{remark}
One can obtain another types of totally nonnegative matrices which
result matrices with self-interlacing spectra after multiplication
by the matrix $J$. To do this we need to change the
conditions~\eqref{matrix.entries.1}--\eqref{matrix.entries.2} by
another ones such that, for instance, the matrix $B^4$, or $B^6$,
(or $B^8$ etc.) becomes oscillating.
\end{remark}

Theorem~\ref{teorem.oscill_signreg.matrix.1} implies the following
corollary.
\begin{corol}
A nonsingular matrix $A$ with positive entries is totally
nonnegative if and only if the matrix $B=JA$ (or the matrix
$C=AJ$) is sign definite of class~$n^{+}$ with the signature
sequence~\eqref{SI_matrix.eps}. Moreover, the matrix $B$ has a
self-interlacing spectrum of the
form~\eqref{SI_matrix.real.spectra}.
\end{corol}

Consider a particular case of conditions~\eqref{matrix.entries.1}.
Suppose that all diagonal entries if the matrix $A$ be positive
$a_{jj}>0$, $j=1,2,\ldots,n$. It is easy to see that in this case
the conditions~\eqref{matrix.entries.1} hold if $a_{j,j+1}>0$,
$j=1,2,\ldots,n-1$. If all remaining entries of the matrix $A$ are
nonnegative, then by Theorem~\ref{teorem.oscill_signreg.matrix.1}
the matrices $B=JA$ and $C=AJ$ have self-interlacing spectra if
$A$ is totally nonnegative.

If all other entries of the matrix $A$ (that is, all entries
except $a_{jj}$ and $a_{j,j+1}$, which are positive) equal zero,
then $A$ is a bidiagonal matrix with positive entries on and under
the main diagonal. Clearly, $A$ is totally nonnegative. Then by
Theorem~\ref{teorem.oscill_signreg.matrix.1} the matrix $B=JA$ has
a self-interlacing spectrum. Thus we obtain the following
statement established in~\cite{H4}:
\begin{theorem}
Any anti-bidiagonal $n\times n$ matrix with positive entries
\begin{equation}\label{Matrix.antybidiag}
B=
\begin{pmatrix}
    0 & 0 &0&\dots&  0   & 0 & b_n \\
    0 &0&0&\dots&  0 & b_{n-2} & b_{n-1} \\
     0  & 0&0&\dots& b_{n-4}& b_{n-3}& 0 \\
    \vdots&\vdots&\vdots&\ddots&\vdots&\vdots\\
     0& 0 &c_{n-4}&\dots&0&0&0\\
      0 &c_{n-2}&c_{n-3}&\dots&0&0&0\\
     c_n&c_{n-1}&  0  &\dots&0&0&0\\
\end{pmatrix}
\end{equation}
has a self-interlacing spectrum of the
form~\eqref{SI_matrix.real.spectra}. Here the entries $b_j>0$,
$j=2,3,\ldots,n$, lie above the main diagonal, the entries
$c_j>0$, $j=2,3,\ldots,n$, lie under the main diagonal and the
only entry on the main diagonal is $a_1>0$.
\end{theorem}

\vspace{2mm}

In~\cite{H4} it was also proved that the spectrum of the
matrix~\eqref{Matrix.antybidiag} coincides with the spectrum of
the following tridiagonal matrix
\begin{equation*}
K=
\begin{pmatrix}
    a_1 & b_2 &  0 &\dots&   0   & 0 \\
    c_2 & 0 &b_3 &\dots&   0   & 0 \\
     0  &c_3 & 0 &\dots&   0   & 0 \\
    \vdots&\vdots&\vdots&\ddots&\vdots&\vdots\\
     0  &  0  &  0  &\dots&0&b_n\\
     0  &  0  &  0  &\dots&c_n&0\\
\end{pmatrix}.
\end{equation*}
It is well-known~\cite{KreinGantmaher} that the spectrum of this
matrix does not depend on the entries $b_j$ and $c_j$ separately.
It depends on products $b_jc_j$, $j=2,3,\ldots,n$. So in order to
the matrices~\eqref{Matrix.antybidiag} and $K$ to have
self-interlacing spectra, it is sufficient that the inequalities
$a_1>0$ and $b_jc_j>0$, $j=2,3,\ldots,n$, hold.

Finally, consider a tridiagonal matrix
\begin{equation*}
M_{J}=
\begin{pmatrix}
    a_1 & b_1 &  0 &\dots&   0   & 0 \\
    c_1 & a_2 &b_2 &\dots&   0   & 0 \\
     0  &c_2 & a_3 &\dots&   0   & 0 \\
    \vdots&\vdots&\vdots&\ddots&\vdots&\vdots\\
     0  &  0  &  0  &\dots&a_{n-1}&b_{n-1}\\
     0  &  0  &  0  &\dots&c_{n-1}&a_n\\
\end{pmatrix},
\end{equation*}
where $a_k,b_k,c_k\in\mathbb{R}$ and $c_kb_k\neq0$.
In~\cite{KreinGantmaher}, there was proved the following fact.
\begin{theorem}\label{teorem.SI_Matrix.Jacobi}
The matrix $M_{J}$ is oscillating if and only if all the entries
$b_k$ and $c_k$ are positive and all the leading principal minors
of $M_{J}$ are also positive:
\begin{equation}\label{Jacobi.ineq}
a_1>0,\;
\begin{vmatrix}
a_1&b_1\\
c_1&a_2
\end{vmatrix}>0,\;
\begin{vmatrix}
a_1&b_1&0\\
c_1&a_2&b_2\\
0&c_2&a_3
\end{vmatrix}>0,\;
\ldots\;
\begin{vmatrix}
    a_1 & b_1 &  0 &\dots&   0   & 0 \\
    c_1 & a_2 &b_2 &\dots&   0   & 0 \\
     0  &c_2 & a_3 &\dots&   0   & 0 \\
    \vdots&\vdots&\vdots&\ddots&\vdots&\vdots\\
     0  &  0  &  0  &\dots&a_{n-1}&b_{n-1}\\
     0  &  0  &  0  &\dots&c_{n-1}&a_n\\
\end{vmatrix}>0.
\end{equation}
\end{theorem}

This theorem together with
Theorem~\ref{teorem.oscill_signreg.matrix.1} implies the following
statement.
\begin{theorem}
The anti-tridiagonal matrix
\begin{equation*}
A_{J}=
\begin{pmatrix}
    0 & 0 &0&\dots&  0   & b_1 & a_1 \\
    0 &0&0&\dots&  b_2 & a_2 & c_1 \\
     0  & 0&0&\dots& a_3& c_2& 0 \\
    \vdots&\vdots&\vdots&\ddots&\vdots&\vdots\\
     0&b_{n-2}&a_{n-2}&\dots&0&0&0\\
     b_{n-1}&a_{n-1}&c_{n-2}&\dots&0&0&0\\
     a_n&c_{n-1}&  0  &\dots&0&0&0\\
\end{pmatrix},
\end{equation*}
where $a_j,b_j,c_j>0$ for $j=1,2,\ldots,n-1$, is sign definite of
class $n^{+}$ and has a self-interlacing
spectrum~\eqref{SI_matrix.real.spectra} if and only if the
following inequalities hold:
\begin{equation}\label{antyJacobi.ineq}
(-1)^{\frac{k(k-1)}2}A_{J}
\begin{pmatrix}
    1&2&\dots&k\\
    n+1-k&n+2-k&\dots&n\\
\end{pmatrix}>0,\quad k=1,2,\ldots,n.
\end{equation}
\end{theorem}
\begin{proof}
If the matrix $A_{J}$ is sign definite of class $n^{+}$ and has a
self-interlacing spectrum of the
form~\eqref{SI_matrix.real.spectra}, then according to
Theorem~\ref{teorem.signreg.matrix}, the signs of its nonzero
minors can be calculated by the formula~\eqref{SI_matrix.eps}.
This implies the inequalities~\eqref{antyJacobi.ineq}.

Conversely, let the inequalities~\eqref{antyJacobi.ineq} hold.
Then we have that the inequalities~\eqref{Jacobi.ineq} hold for
the~matrix $M_J=JA_J$. By Theorem~\ref{teorem.SI_Matrix.Jacobi},
the matrix $M_J$ is oscillating and, in particular, totally
nonnegative. Now notice that $(M_J)_{ii}>0$, $i=1,\ldots,n$, and
$(M_J)_{k,k+1}>0$, $k=1,\ldots,k-1$, so $M_J$ satisfies the
condition~\eqref{matrix.entries.1} of
Theorem~\ref{teorem.oscill_signreg.matrix.1}. Therefore, the
matrix $A_J$ is sign definite of class $n^{+}$ and has a
self-interlacing spectrum of the
form~\eqref{SI_matrix.real.spectra}.
\end{proof}


\setcounter{equation}{0}

\section{Generalized Hurwitz polynomials}\label{section:Generalized.Hurwitz}

\hspace{4mm} In this section, we describe the class of real
polynomials whose associated function $\Phi$ is an
\textit{R}-function. In particular, this class includes all
Hurwitz stable and self-interlacing polynomials as extremal cases.

\subsection{General theory}\label{subsection:General.Hurwitz.general.theory}

\hspace{4mm} Let us consider a polynomial
\begin{equation}\label{general.Hurwitz.poly}
p(z)=a_0z^n+a_1z^{n-1}+\dots+a_n,\qquad a_1,\dots,a_n\in\mathbb
R,\ a_0>0.
\end{equation}

\begin{definition}\label{def.general.Hurw.poly}
The polynomial $p$ is called \textit{generalized Hurwitz}
polynomial of type~I of order $k$, where $1\leqslant
k\leqslant\left[\dfrac{n+1}2\right]$, if it has exactly $k$ zeroes
in the closed right half-plane, all of which are nonnegative and
simple:
\begin{equation}\label{Hurwitz.poly.def}
0\leqslant\mu_1<\mu_2<\cdots<\mu_k,
\end{equation}
such that $p(-\mu_i)\neq0$, $i=1,\ldots,k$, and $p$ has an odd
number of zeroes, counting multiplicities, on each interval
$(-\mu_{k},-\mu_{k-1}),\ldots,(-\mu_3,-\mu_2),(-\mu_2,-\mu_1)$.
Moreover, the number of zeroes of $p$ on the interval~$(-\mu_1,0)$
(if any) is even, counting multiplicities. The other real zeroes
lie on the interval $(-\infty,-\mu_{k})$: an odd number of zeroes,
counting multiplicities, when $n=2l$, and an even number of
zeroes, counting multiplicities, when $n=2l+1$. All nonreal zeroes
of $p$ (if any) are located in the open left half-plane of the
complex plane.
\end{definition}

\begin{definition}\label{def.general.Hurw.poly.II}
If $p(z)$ is a generalized Hurwitz polynomials of type~I, then the
polynomial $p(-z)$ is called generalized Hurwitz polynomial of
type~II.
\end{definition}

It clear that all results obtained for the generalized Hurwitz
polynomials of type~I can be easily reformulated for  the
generalized Hurwitz polynomials of type~I. Thus, in the rest of
the section we consider only generalized Hurwitz polynomials of
type~I unless explicitly stipulated otherwise.

Since the generalized Hurwitz polynomials of order $k$ have
exactly $k$ zeroes in the closed right half-plane, the generalized
Hurwitz polynomials of order $0$ have no zeroes in the closed
right half-plane, so they are \textit{Hurwitz stable}.

Analogously, if $p$ is a generalized Hurwitz polynomial of order
$k=\left[\dfrac{n+1}2\right]$ without a root at zero, then $p$ is
\textit{self-interlacing} of type~I. In fact, let
$0<\mu_1<\mu_2<\ldots<\mu_k$ be its positive zeroes. Then by
definition, $p$ has an odd number of zeroes, counting
multiplicities, (\textit{at least} one) on each interval
$(-\mu_{j+1},-\mu_j)$, $j=1,\ldots,k-1$. Moreover, if $n=2l$, then
$p$ has at least one zero on the interval $(-\infty,-\mu_{k})$. It
is easy to see that $p$ can not have nonreal zeroes and has
exactly one simple zero on each interval $(-\mu_{j+1},-\mu_j)$,
$j=1,\ldots,k-1$, and on the interval $(-\infty,-\mu_k)$ for
$n=2l$. Besides, the zeroes of~$p$ are distributed as
in~\eqref{self-interlacing.zero.distribution.I.type}.

Thus, Hurwitz stable and self-interlacing polynomials are
generalized Hurwitz polynomials of minimal and maximal orders,
respectively. We note that generalized Hurwitz polynomials of
maximal order with a root at zero are almost self-interlacing.

Now we are in a position to establish the main theorem of the
theory of generalized Hurwitz polynomials. This theorem is a
direct generalization of
Theorems~\ref{Theorem.main.Hurwitz.stability}
and~\ref{Theorem.main.self-interlacing}.

\begin{theorem}\label{Theorem.main.general.Hurwitz}
Let $p$ be a given real polynomial of degree $n\geqslant1$ as
in~\eqref{general.Hurwitz.poly}. The polynomial $p$ is generalized
Hurwitz if and only if its associated function~$\Phi$ defined
in~\eqref{assoc.function} is an \textit{R}-function with
exactly~$l=\left[\dfrac{n}2\right]$ \emph{(}for $a_1\neq0$\emph{)}
or $l-1$ \emph{(}for $a_1=0$\emph{)} poles\footnote{Let us note
that if $n=2l$ and $\Phi$ is an \textit{R}-function, then
$a_1\neq0$, so $\Phi$ has exactly $l$ poles.}.
Moreover,
\begin{itemize}
\item if $n=2l$ or if $n=2l+1$ with
\begin{equation}\label{general.Hurwitz.poly.0.1}
a_0a_1>0,
\end{equation}
then the number of nonnegative poles of the function $\Phi$ equals
order of $p$;
\item if $n=2l+1$ with
\begin{equation}\label{general.Hurwitz.poly.0.2}
a_0a_1\leqslant0,
\end{equation}
then the number of nonnegative poles of the function $\Phi$ equals
order of $p$ minus one.
\end{itemize}
\end{theorem}
\begin{proof}
As we said above, for $k=0$ and $k=\left[\dfrac{n+1}2\right]$ this
theorem is true.

\noindent\textit{Sufficiency}. Let the function $\Phi$ be an
\textit{R}-function with exactly $r$ nonpositive poles,
$0\leqslant r\leqslant\left[\dfrac{n}2\right]$.

At first, assume that $a_1\neq0$ and the function $\Phi$ has no
pole at zero. Then $\Phi$ has exactly $l=\left[\dfrac{n}2\right]$
poles and can be represented in the following form by
Theorem~\ref{Th.R-function.general.properties}:

\begin{equation}\label{Mittag.Leffler.general.Hurwitz}
\Phi(u)=\alpha+\sum_{j=1}^r\frac{\alpha_j}{u-\omega_j}+\sum_{j=r+1}^l\frac{\alpha_j}{u+\omega_j},
\end{equation}
where $\alpha\in\mathbb{R}$ ($\alpha=0$ if and only if $n=2l$),
$\alpha_j>0,\omega_j>0$ for $j=1,\ldots,l$, all poles  are
distinct. Let the positive poles are enumerated as follows
$0<\omega_1<\omega_2<\ldots<\omega_r$.

As in the proofs of Theorems~\ref{Theorem.main.Hurwitz.stability}
and~\ref{Theorem.main.self-interlacing}, we note that the set of
zeroes of the polynomial $p$ coincides with the set of roots of
the following equation
\begin{equation}\label{Theorem.main.general.Hurwitz.proof.1}
\Phi(z^2)=-\dfrac1{z}.
\end{equation}
So if we suppose that $p(\lambda)=0$ and $\Im\lambda\neq0$ (that
is, $\lambda$ is a nonreal zero of $p$), then
from~\eqref{Theorem.main.general.Hurwitz.proof.1} we obtain
\begin{equation}\label{Theorem.main.general.Hurwitz.proof.2}
\displaystyle\frac{\Im\Phi(\lambda^2)}{\Im(\lambda^2)}=\frac{1}{|\lambda|^2}\frac{\Im\lambda}{\Im
\lambda^2}= \frac{1}{2|\lambda|^2}\frac{1}{\Re\lambda}<0,
\end{equation}
since $\Phi$ is an \textit{R}-function by assumption.
From~\eqref{Theorem.main.general.Hurwitz.proof.2} it follows that
all nonreal zeroes of $p$ lie in the open left half-plane. Thus,
all zeroes of $p$ are located in the open left half-plane or on
the nonnegative real half-axis.

As in the proof of Theorem~\ref{Theorem.main.self-interlacing}, we
investigate the distribution of the real zeroes of the
polynomial~$p$. By assumption, $|\Phi(0)|<\infty$, so the
polynomial $p$ does not vanish at zero. In fact, if $p(0)=0$, then
$p_0(0)=0$, but it contradicts with $|\Phi(0)|<\infty$. Thus, we
suppose that $z\in\mathbb{R}\setminus\{0\}$ and change the
variables as follows $z^2=u>0$. Then for real nonzero $z$, the
equation~\eqref{Theorem.main.general.Hurwitz.proof.1} is
equivalent to the following system
\begin{equation}\label{Theorem.main.general.Hurwitz.proof.3}
\begin{cases}
\ \Phi(u)=-\dfrac1{\sqrt{u}},\\
\ \Phi(u)=\dfrac1{\sqrt{u}},
\end{cases}\quad u>0.
\end{equation}
We also note (as in the proof of
Theorem~\ref{Theorem.main.self-interlacing}) that all positive
roots (if any) of the first
equation~\eqref{Theorem.main.general.Hurwitz.proof.3} are squares
of positive roots of
the~equation~\eqref{Theorem.main.general.Hurwitz.proof.1}, and all
positive roots (if any) of the second
equation~\eqref{Theorem.main.general.Hurwitz.proof.3} are squares
of negative roots of
the~equation~\eqref{Theorem.main.general.Hurwitz.proof.1}. The
equation~\eqref{Theorem.main.general.Hurwitz.proof.1} may have
nonreal roots, but as we showed above all these roots are located
in the open left half-plane.

Let $n=2l$. Then $\lim\limits_{u\to\pm\infty}\Phi(u)=\alpha=0$.
The set of zeroes the function $F_1(u)=\Phi(u)+\dfrac1{\sqrt{u}}$
coincides with the set of roots of the first
equation~\eqref{Theorem.main.general.Hurwitz.proof.3}. In the same
way as in the proof of
Theorem~\ref{Theorem.main.self-interlacing}, one can show that the
function $F_1$ has exactly one simple zero, say~$\mu_i^2$
($\mu_i>0$), on each interval $(\omega_{i-1},\omega_i)$,
$i=2,\ldots,r$ and exactly one simple zero~$\mu_1^2$ ($\mu_1>0$)
on the interval $(0,\omega_1)$. Moreover, denoting by $\xi_i$ a
unique zero\footnote{Recall that $\Phi$ is an \textit{R}-function
by assumption. Consequently, it has exactly one simple zero on
each interval $(\omega_i,\omega_{i+1})$.} of the function $\Phi$
on the interval~$(\omega_i,\omega_{i+1})$, $i=1,\ldots,r-1$, we
have
\begin{equation}\label{Theorem.main.general.Hurwitz.proof.4}
\mu_i^2<\omega_i<\xi_i<\mu_{i+1}^2<\omega_{i+1},\qquad
i=1,\ldots,r-1,
\end{equation}
since $F_1(\xi_i)=\dfrac1{\sqrt{\xi_i}}>0$, and $F_1$ is
decreasing between its poles on the positive half-axis as a sum of
decreasing functions.

Thus, the first
equation~\eqref{Theorem.main.general.Hurwitz.proof.3} has exactly
$r$ positive roots $\mu_i^2$ all of which are simple and satisfy
the inequalities~\eqref{Theorem.main.general.Hurwitz.proof.4}.

Now consider the function $F_2(u)=\Phi(u)-\dfrac1{\sqrt{u}}$ whose
set of zeroes coincides with the set of roots of the second
equation~\eqref{Theorem.main.general.Hurwitz.proof.3}. In the same
way as in the proof of
Theorem~\ref{Theorem.main.self-interlacing}, one can show that the
function~$F_2$ has an odd number of zeroes, counting
multiplicities, on each interval $(\omega_{i},\mu^2_{i+1})$,
$i=1,\ldots,r-1$, and on the interval $(\omega_r,+\infty)$.
Moreover, $F_2(u)<0$ on the intervals
$[\mu^2_{i+1},\omega_{i+1})$, $i=1,\ldots,r-1$. Let us now turn
our attention to the interval $(0,\omega_1)$. Since $|\Phi(0)|<0$
by assumption and $-\dfrac1{\sqrt{u}}\to-\infty$ as $u\to+0$, we
have $F_2(u)\to-\infty$ as $u\to+0$. Besides,
$F_2(u)=F_1(u)-\dfrac2{\sqrt{u}}$, so
$F_2(\mu_1^2)=-\dfrac2{\mu_1}<0$. Therefore, $F_2$~has an even
number of zeroes, counting multiplicities, on the interval
$(0,\mu_1^2)$. However, since $F_1(u)$ is monotone decreasing on
$(0,\omega_1)$ as a sum of two monotone decreasing functions, we
have $F_1(u)<0$ on the interval $(\mu_1^2,\omega_1)$.
Consequently, $F_2(u)$ is also negative for
$u\in(\mu_1^2,\omega_1)$, so its zeroes on the interval
$(0,\omega_1)$ (if any) lie, indeed, on the interval
$(0,\mu_1^2)$.

Thus, we obtain that in the closed right half-plane, the
polynomial $p$ has only~$r$ zeroes, all of which are positive and
simple:
\begin{equation*}\label{Theorem.main.general.Hurwitz.proof.5}
0<\mu_1<\mu_2<\ldots<\mu_r.
\end{equation*}
Moreover, $p$ has an even number of zeroes, counting
multiplicities, on the interval~$(-\mu_1,0)$ and an odd number of
zeroes, counting multiplicities, on each interval
$(-\mu_{i+1},-\mu_{i})$, $i=1,\ldots,r-1$. Also $p(-\mu_{i})\neq0$
for $i=1,\ldots,r$. As we showed above, all nonreal zeroes of $p$
are in the open left-half-plane, so $p$ is a generalized Hurwitz
polynomial of order~$r$.

\vspace{2mm}

Let $n=2l+1$, and the condition~\eqref{general.Hurwitz.poly.0.1}
holds. The only difference between this case and the case $n=2l$
is behaviour of the function $\Phi$ at infinity. Thus, as in the
case $n=2l$, $F_1$ has a unique simple zero, say~$\mu_i^2$
($\mu_i>0$), $i=1,\ldots,r$, on each interval $(0,\omega_1)$,
$(\omega_1,\omega_2)$, \ldots, $(\omega_{r-1},\omega_r)$. These
zeroes are distributed as
in~\eqref{Theorem.main.general.Hurwitz.proof.4}. Furthermore, we
have $F_1(u)\to+\infty$ as $u\searrow\omega_r$ and
$F_1(u)\to\alpha=\dfrac{a_0}{a_1}>0$ as $u\to+\infty$, so
$F_1(u)>0$ for $u\in(\omega_r,+\infty)$ (recall that $\Phi$ is a
monotone decreasing function on the interval $(\omega_r,+\infty)$
and $\Phi(u)\to+\infty$ as $u\searrow\omega_r$). As in the case
$n=2l$, the function $F_2(u)$ has an even number of zeroes,
counting multiplicities, on the interval $(0,\mu_1^2)$ and an odd
number of zeroes, counting multiplicities, on each interval
$(\omega_i,\mu^2_{i+1})$ for $i=1,\ldots,r-1$. Moreover, $F_2$ has
no zeroes on the intervals $(\mu_1^2,\omega_1)$ and
$(\mu^2_{i+1},\omega_{i+1})$, $i=1,\ldots,r-1$. At the same time,
we have $F_2(u)\to+\infty$ as $u\searrow\omega_r$ and
$F_2(u)\to\alpha>0$ as $u\to+\infty$. Since the function $F_2(u)$
is not monotone, it has an even number of zeroes, counting
multiplicities, on the interval $(\omega_r,+\infty)$. Recall that
these zeroes are squares of positive zeroes of the polynomial~$p$
on the interval $(\omega_r,+\infty)$. Since the other zeroes of
the functions $F_1$ and $F_2$ are distributed as in the case
$n=2l$, we obtain that $p$ is a generalized Hurwitz polynomial of
order~$r$.

\vspace{2mm}

Let now $n=2l+1$, and let $a_0a_1<0$. This case also differs from
the case $n=2l$ by behaviour of the function $\Phi$ at infinity:
$\lim\limits_{u\to+\infty}\Phi(u)=\alpha=\dfrac{a_0}{a_1}<0$. As
in the case $n=2l$, the function
$F_1(u)=\Phi(u)+\dfrac1{\sqrt{u}}$ has a unique simple zero,
say~$\mu_i^2$ ($\mu_i>0$), $i=1,\ldots,r$, on each interval
$(0,\omega_1)$, $(\omega_1,\omega_2)$, \ldots,
$(\omega_{r-1},\omega_r)$. These zeroes are distributed as
in~\eqref{Theorem.main.general.Hurwitz.proof.4}. Furthermore, we
have $F_1(u)\to+\infty$ as $u\searrow\omega_r$ and
$F_1(u)\to\alpha<0$ as $u\to+\infty$. Since $F_1(u)$ is decreasing
on the interval $(\omega_r,+\infty)$ as a sum of two decreasing
functions on this interval, it has a unique simple zero,
say~$\mu_{r+1}^2$ ($\mu_{r+1}>0$) on $(\omega_r,+\infty)$ such
that
\begin{equation*}\label{Theorem.main.general.Hurwitz.proof.6}
\omega_r<\xi_r<\mu^2_{r+1},
\end{equation*}
where $\xi_r$ is a zero\footnote{Recall that this zero is unique,
since $\Phi$ is a monotone function on $(\omega_r,+\infty)$.} of
$\Phi$ on the interval $(\omega_r,+\infty)$. As in the case
$n=2l$, the function $F_2(u)$ has an even number of zeroes,
counting multiplicities, on the interval $(0,\mu_1^2)$ and an odd
number of zeroes, counting multiplicities, on each interval
$(\omega_i,\mu^2_{i+1})$ for $i=1,\ldots,r-1$. Moreover, $F_2$ has
no zeroes on the intervals $(\mu^2_{i},\omega_{i})$,
$i=1,\ldots,r$. Now consider the interval $(\omega_r,+\infty)$.
Since $F_2(u)\to+\infty$ as $u\searrow\omega_r$ and
$F_2(\mu^2_{r+1})=-\dfrac1{\mu_{r+1}}<0$, we obtain that $F_2$ has
an odd number of zeroes, counting multiplicities, on the interval
$(\omega_r,\mu^2_{r+1})$. Moreover,
$F_2(u)=F_1(u)-\dfrac2{\sqrt{u}}<0$ for $(\mu^2_{r+1},+\infty)$,
because $F_1$ is negative on that interval. Since the positive
roots of the first
equation~\eqref{Theorem.main.general.Hurwitz.proof.3} are squares
of the positive zeroes of the polynomial $p$ and the roots of the
second equation~\eqref{Theorem.main.general.Hurwitz.proof.3} are
squares of the negative zeroes of the polynomial~$p$, we obtain
that that $p$ has only~$r+1$ zeroes  in the closed right
half-plane, all of which are positive and simple:
\begin{equation*}\label{Theorem.main.general.Hurwitz.proof.7}
0<\mu_1<\mu_2<\ldots<\mu_{r+1}.
\end{equation*}
Moreover, $p$ has an even number of zeroes, counting
multiplicities, on the interval~$(-\mu_1,0)$ and an odd number of
zeroes, counting multiplicities, on each interval
$(-\mu_{i+1},-\mu_i)$, $i=1,\ldots,r$. All nonreal zeroes of $p$
are in the open left half-plane, and $p(-\mu_i)\neq0$ for
$i=1,\ldots,r+1$, so $p$ is a generalized Hurwitz polynomial of
order~$r+1$.

\vspace{3mm}

If $n=2l+1$, and $a_1=0$, then by
Theorem~\ref{Th.R-function.general.properties}, the function
$\Phi$ has exactly $l-1$ poles and can be represented in the form
\begin{equation}\label{Mittag.Leffler.general.Hurwitz.second}
\Phi(u)=\beta u+\gamma
+\sum_{j=1}^r\frac{\alpha_j}{u-\omega_j}+\sum_{j=r+1}^{l-1}\frac{\alpha_j}{u+\omega_j},
\end{equation}
where $\beta<0$, $\gamma\in\mathbb{R}$, $\alpha_j>0,\omega_j>0$
for $j=1,\ldots,l$, all poles are distinct. In the same way as in
the case $n=2l+1$ with $a_0a_1<0$, one can show that $p$ is a
generalized Hurwitz polynomial of order~$r+1$, since $\Phi(u)$ is
also negative for sufficiently large $u$.

\vspace{3mm}

Suppose now that the function $\Phi$ has a pole at zero. Then
$p_0(0)=0$, so $p(0)=a_n=0$. At the same time,
$p_1(0)=a_{n-1}=p'(0)\neq0$, since $\Phi$ is an
\textit{R}-function by assumption, so it has interlacing poles and
zeroes (see Theorem~\ref{Th.R-function.general.properties}). Thus,
$p$ has a simple root at zero. Note that the set of zeroes of the
polynomial~$p$, but the root at zero, coincides with the set of
roots of the
equation~\eqref{Theorem.main.general.Hurwitz.proof.1}.

If $n=2l$ or $n=2l+1$ with~\eqref{general.Hurwitz.poly.0.1}, then
as above, one can show that in the closed right half-plane $p$ has
exactly $r$ zeroes, say~$\mu_i$, $i=1\ldots,r$, all of which are
nonnegative and simple:
\begin{equation*}\label{Theorem.main.general.Hurwitz.proof.8}
0=\mu_1<\mu_2<\mu_2<\ldots<\mu_r,
\end{equation*}
and $p(-\mu_i)\neq0$ for $i=2,\ldots,r$. Also $p$ has an odd
number of zeroes, counting multiplicities, on each interval
$(-\mu_{i},-\mu_{i-1})$ and an odd (even) number of zeroes on the
interval $(-\infty,-\mu_r)$ if $n=2l$ ($n=2l+1$). All nonreal
zeroes of $p$ lie in the open left half-plane. Thus, $p$ is
generalized Hurwitz of order $r$.

If $n=2l+1$ and the condition~\eqref{general.Hurwitz.poly.0.2}
holds, then in the same way as above, one can show that $p$ is
generalized Hurwitz of order $r+1$.

\vspace{4mm}

\textit{Necessity}. Let $p$ be a generalized Hurwitz polynomial of
order $k$. At first, suppose that $p(0)\neq0$.

\vspace{2mm}

\noindent Let $n=2l$. Then by definition, $p(z)$ has only $k$
zeroes in the closed right half-plane, all of which are positive
and simple: $0<\mu_1<\mu_2<\ldots<\mu_k$, and $p(-\mu_i)\neq0$.
Moreover, $p$ has an odd number of zeroes, counting
multiplicities, on each interval $(-\infty,-\mu_k)$,
$(-\mu_k,-\mu_{k-1})$,\ldots, $(-\mu_2,-\mu_1)$, and it has an
even number of zeroes, counting multiplicities, on the interval
$(-\mu_1,0)$. All nonreal zeroes of $p$ lie in the open left
half-plane. Thus, we can represent the polynomial $p$ as a product
$p(z)=h(z)g(z)$, where $h$ is a self-interlacing
polynomial\footnote{More exactly, $h$ is self-interlacing of
type~I.}, of degree $2k$ whose positive zeroes are $\mu_i$,
$i=1,\ldots,k$, and $g$ is a Hurwitz stable polynomial of degree
$2(l-k)$, which has an even number of zeroes, counting
multiplicities, on each interval $(-\infty,-\mu_k)$,
$(-\mu_k,-\mu_{k-1})$,\ldots, $(-\mu_2,-\mu_1)$, $(-\mu_1,0)$. Let
\begin{equation*}\label{Theorem.main.general.Hurwitz.proof.9}
h(z)=h_0(z^2)+zh_1(z^2),
\end{equation*}
\begin{equation*}\label{Theorem.main.general.Hurwitz.proof.10}
g(z)=g_0(z^2)+zg_1(z^2).
\end{equation*}
Then by Theorem~\ref{Theorem.main.self-interlacing}, the function
$H:=\dfrac{h_1}{h_0}$ is an \textit{R}-function with exactly $k$
poles, all of which are positive. Analogously, by
Theorem~\ref{Theorem.main.Hurwitz.stability}, the function
$G:=\dfrac{g_1}{g_0}$ is an $R$ function with exactly $l-k$ poles,
all of which are negative. If $p(z)=p_0(z^2)+zp_1(z^2)$, then
\begin{equation}\label{general.Hurwitz.poly.4}
p_1(u)=h_0(u)g_1(u)+h_1(u)g_0(u),
\end{equation}
\begin{equation}\label{general.Hurwitz.poly.5}
p_0(w)=h_0(u)g_0(u)+uh_1(u)g_1(u).
\end{equation}
Thus,
\begin{equation*}\label{Theorem.main.general.Hurwitz.proof.11}
\displaystyle\Phi(u)=\frac{p_1(u)}{p_0(u)}=\frac{h_0(u)g_1(u)+h_1(u)g_0(u)}{h_0(u)g_0(u)+uh_1(u)g_1(u)}=
\frac{H(u)+G(u)}{1+uH(u)G(u)}.
\end{equation*}
Note that the function
\begin{equation}\label{Theorem.main.general.Hurwitz.proof.12}
\displaystyle\frac{h_1(u)}{h_0(u)}+\frac{g_1(u)}{g_0(u)}=\frac{p_1(u)}{h_0(u)g_0(u)}
\end{equation}
is an \textit{R}-function as a sum of \textit{R}-functions.
Consequently, all zeroes of $p_1$ are real and simple and
interlace zeroes of the polynomial $h_0g_0$, which has $k$
positive simple zeroes and $l-k$ negative simple zeroes.

Furthermore, the function
\begin{equation}\label{general.Hurwitz.poly.6}
\displaystyle
F(u)=-\frac{h_0(u)}{h_1(u)}-\frac{ug_1(u)}{g_0(u)}=-\dfrac1{H(u)}-uG(u)=-\frac{p_0(u)}{h_1(u)g_0(u)}
\end{equation}
is an \textit{R}-function by
Theorem~\ref{Theorem.properties.R-functions}, since the functions
$H$ and $G$ are \textit{R}-functions. Consequently, all zeroes of
the polynomial $p_0$ are real and simple and interlace zeroes of
the polynomial~$h_1g_0$ (see
Theorem~\ref{Th.R-function.general.properties}). Thus,
from~\eqref{Theorem.main.general.Hurwitz.proof.12}--\eqref{general.Hurwitz.poly.6}
it follows that all poles and zeroes of the function
$\Phi=\dfrac{p_1}{p_0}$ are real and simple. Therefore, $\Phi$ can
be represented as follows
\begin{equation}\label{Theorem.main.general.Hurwitz.proof.13}
\Phi(u)=\dfrac{p_1(u)}{p_0(u)}=\sum\limits_{j=1}^{l}\dfrac{\alpha_j}{u-\omega_j},
\end{equation}
where $\omega_i\in\mathbb{R}$, $\omega_i\neq\omega_j$ for $i\neq
j$, and
\begin{equation*}\label{Theorem.main.general.Hurwitz.proof.14}
\displaystyle\alpha_j=\frac{p_1(\omega_j)}{p'_0(\omega_j)}\neq0,\qquad
j=1,2,\ldots,l.
\end{equation*}
According to Theorem~\ref{Th.R-function.general.properties}, now
it is sufficient to establish positivity of the numbers $\alpha_j$
to prove that $\Phi$ is an \textit{R}-function.

Let $\omega_j$ be a zero of the polynomial $p_0$. Since
$p(0)\neq0$ by assumption, we have $p_0(0)\neq0$, so
$\omega_j\neq0$, $j=1,2,\ldots,l$. It is clear that
$h_1(\omega_j)g_0(\omega_j)\neq0$, because the zeroes of the
polynomial $p_0$ interlace the zeroes of $h_1g_0$.
From~\eqref{general.Hurwitz.poly.5} it follows that
\begin{equation*}\label{Theorem.main.general.Hurwitz.proof.15}
\displaystyle
h_0(\omega_j)=-\frac{\omega_jh_1(\omega_j)g_1(\omega_j)}{g_0(\omega_j)},
\end{equation*}
so
\begin{equation*}\label{Theorem.main.general.Hurwitz.proof.16}
\begin{array}{l}
\displaystyle
p_1(\omega_j)=h_1(\omega_j)g_0(\omega_j)-\frac{\omega_jh_1(\omega_j)g_1^2(\omega_j)}{g_0(\omega_j)}=\\=
\dfrac{h_1(\omega_j)}{g_0(\omega_j)}\left[g_0^2(\omega_j)-\omega_jg_1^2(\omega_j)\right]
=\dfrac{h_1(\omega_j)}{g_0(\omega_j)}g(\sqrt{\omega_j})g(-\sqrt{\omega_j}).
\end{array}
\end{equation*}
According to~\eqref{general.Hurwitz.poly.6}, we have
$p_0(u)=-h_1(u)g_0(u)F(u)$, and $F(\omega_j)=0$, therefore
\begin{equation*}
p_0'(\omega_j)=-h_1(\omega_j)g_0(\omega_j)F'(\omega_j).
\end{equation*}
So we obtain
\begin{equation*}\label{Theorem.main.general.Hurwitz.proof.17}
\alpha_j=\dfrac{p_1(\omega_j)}{p'_0(\omega_j)}=
\dfrac{h_1(\omega_j)}{g_0(\omega_j)}\cdot\dfrac{g(\sqrt{\omega_j})g(-\sqrt{\omega_j})}{-h_1(\omega_j)g_0(\omega_j)F'(\omega_j)}=
-\dfrac{g(\sqrt{\omega_j})g(-\sqrt{\omega_j})}{g_0^2(\omega_j)F'(\omega_j)}.
\end{equation*}
Since $F$ is an \textit{R}-function, we have $F'(u)<0$, so
$g_0^2(\omega_j)F'(\omega_j)<0$. Thus, we obtain that
\begin{equation*}\label{Theorem.main.general.Hurwitz.proof.18}
\textrm{sign}(\alpha_j)=\textrm{sign}\left(g(\sqrt{\omega_j})g(-\sqrt{\omega_j})\right).
\end{equation*}

If $\omega_j<0$, then
$g(\sqrt{\omega_j})g(-\sqrt{\omega_j})=g_0^2(\omega_j)-\omega_jg_1^2(\omega_j)>0$,
since $g_0$ and $g_1$ are real polynomials. So $\alpha_j>0$ in
this case.

Let now $\omega_j>0$. The polynomial $g$ is Hurwitz stable with
positive leading coefficient by construction, so $g(z)>0$ for all
$z\geqslant0$. Therefore,
$\textrm{sign}(\alpha_j)=\textrm{sign}\left(g(-\sqrt{\omega_j})\right)$.
Since $\deg p_0=l$, $\deg g_0=l-k$, $\deg h_1=k-1$, the function
$F$ defined in~\eqref{general.Hurwitz.poly.6} has exactly $l$
zeroes and $l-1$ poles. Consequently, the maximal positive zero of
the polynomial $p_0$ is greater than the maximal positive zero of
the polynomial~$h_1$ (recall that all zeroes of $h_1$ are
positive). Besides, $F$ is decreasing between its poles, so
by~\eqref{general.Hurwitz.poly.5}, we obtain
$F(0)=-\dfrac{p_0(0)}{h_1(0)g_0(0)}=-\dfrac{h_0(0)}{h_1(0)}>0$,
since $h$ is self-interlacing by construction,
and~\eqref{self-interlacing.last.coeffs} holds for it. This means
that the function $F$ (and the polynomial $p_0$) has a positive
zero, which is less than the minimal positive zero of $h_1$. Thus,
the polynomial $p_0$ has exactly $k$ positive zeroes and,
respectively, exactly $l-k$ negative zeroes.

Let us enumerate the positive zeroes of $p_0$ as follows
\begin{equation*}\label{Theorem.main.general.Hurwitz.proof.19}
0<\omega_1<\omega_{2}<\ldots<\omega_k.
\end{equation*}
Then denoting the zeroes of the polynomial $h_1$ by $\gamma_j$
$(j=1,2,\ldots,k-1)$, we obtain
\begin{equation*}
0<\omega_1<\gamma_1<\omega_2<\gamma_2<\ldots<\omega_k<\gamma_{k-1}<\omega_k.
\end{equation*}

Let us now show that $\omega_1>\mu_1^2$. Indeed, since $\mu_1$ is
the minimal positive zero of $h$ and $h(0)\neq0$, the polynomial
$h$ does not change its sign on the interval $[0,\mu_1)$.
Moreover, both polynomials $h_0$ and $h_1$ have no zeroes  on the
interval $[0,\mu_1^2]$ (see the proof of
Theorem~\ref{Theorem.main.self-interlacing}). Consequently, we get
$\textrm{sign}(h(0))=\textrm{sign}\left(h_0(0)\right)=-\textrm{sign}\left(h_1(0)\right)=
\textrm{sign}\left(h_0(\mu_1^2)\right)=-\textrm{sign}\left(h_1(\mu_1^2)\right)$
that follows from the self-interlacing of $h$. Note that
$\textrm{sign}\left(p_0(0)\right)=\textrm{sign}\left(h_0(0)g_0(0)\right)=\textrm{sign}\left(h_0(0)\right)$,
since $g$ is Hurwitz stable with the positive leading coefficient,
so $g_0(0)>0$. The number $\mu_1$ is a zero of the polynomial $h$,
therefore $h_0(\mu_1^2)=-\mu_1h_1(\mu_1^2)$, so
from~\eqref{general.Hurwitz.poly.5} we have
\begin{equation*}
p_0(\mu_1^2)=h_0(\mu_1^2)g_0(\mu_1^2)-\mu_1h_0(\mu_1^2)g_1(\mu_1^2)=h_0(\mu_1^2)g(-\mu_1).
\end{equation*}

By construction, the polynomial $g$ has an even number of zeroes,
counting multiplicities, on each interval
$(-\infty,-\mu_{k}),(-\mu_{k},-\mu_{k-1}),\ldots,(-\mu_2,-\mu_1),(-\mu_1,0)$,
therefore
$\textrm{sign}\left(g(0)\right)=\textrm{sign}\left(g(-\mu_j)\right)$,
$j=1,2,\ldots,k$, that is, $g(-\mu_j)>0$. This implies that
$\textrm{sign}\left(p_0(\mu_1^2)\right)=\textrm{sign}\left(h_0(\mu_1^2)\right)=\textrm{sign}\left(h_0(0)\right)=
\textrm{sign}\left(p_0(0)\right)$. Consequently, the polynomial
$p_0$ has an even number of zeroes, counting multiplicities, on
the interval $(0,\mu_1^2)$. But the positive zeroes of the
polynomials $p_0$ and $h_1$ interlace, and $h_1$ has no zeroes in
the interval $(0,\mu_1^2)$, therefore $p_0$ has also no zeroes on
that interval, that is, $\omega_1>\mu_1^2$.

Note that (see the proof of
Theorem~\ref{Theorem.main.self-interlacing})
\begin{equation*}
0<\mu_1^2<\gamma_1<\mu_2^2<\gamma_2<\ldots<\mu_{k-1}^2<\gamma_{k-1}<\mu_k^2.
\end{equation*}

Thus, we already established that
$0<\mu_1^2<\omega_1<\gamma_1<\mu_2^2$. Now we prove that the
polynomial $p_0$ has no zeroes on each interval
$(\gamma_j,\mu_{j+1}^2)$, $j=1,2,\ldots,k-1$. In fact, since
$\textrm{sign}\left(p_0(0)\right)=\textrm{sign}\left(h_0(0)\right)$
and both polynomials $p_0$ and $h_0$ change their signs on each
interval $(\gamma_j,\gamma_{j+1})$, $j=1,2,\ldots,k-1$, we obtain
$\textrm{sign}\left(p_0(\gamma_j)\right)=\textrm{sign}\left(h_0(\gamma_j)\right)$
for all $j=1,2,\ldots,k-1$. By
Theorem~\ref{Theorem.main.self-interlacing}, the following
inequalities hold
\begin{equation}\label{general.Hurwitz.poly.7}
0<\mu_1^2<\beta_1<\gamma_1<\mu_2^2<\beta_2<\gamma_2<\ldots<\mu_{k-1}^2<\beta_{k-1}<\gamma_{k-1}<\mu_k^2<\beta_k,
\end{equation}
where $\beta_j$ are the zeroes of the polynomial $h_0$.
From~\eqref{general.Hurwitz.poly.7} it follows that $h_0$ does not
change its sign on each interval $(\gamma_j,\mu_{j+1}^2)$,
$j=1,2,\ldots,k-1$, that is,
$\textrm{sign}\left(h_0(\gamma_j)\right)=\textrm{sign}\left(h_0(\mu_{j+1}^2)\right)$
for all $j=1,2,\ldots,k-1$, in particular. Since $h(\mu_{j+1})=0$
and therefore $h_0(\mu_{j+1}^2)=-\mu_jh_1(\mu_{j+1}^2)$, we have
\begin{equation*}
p_0(\mu_{j+1}^2)=h_0(\mu_{j+1}^2)g_0(\mu_{j+1}^2)-\mu_{j+1}h_0(\mu_{j+1}^2)g_1(\mu_{j+1}^2)=h_0(\mu_{j+1}^2)g(-\mu_{j+1}).
\end{equation*}
The inequalities $g(-\mu_{j+1})>0$, $j=1,\ldots,k-1$, proved above
imply
$\textrm{sign}\left(p_0(\mu_{j+1}^2)\right)=\textrm{sign}\left(h_0(\mu_{j+1}^2)\right)=
\textrm{sign}\left(h_0(\gamma_j)\right)=\textrm{sign}\left(p_0(\gamma_j)\right)$.
Consequently, the polynomial $p_0$ has an even number of zeroes,
counting multiplicities, on each interval
$(\gamma_j,\mu_{j+1}^2)$, $j=1,\ldots,k-1$. But the zeroes of the
polynomials $p_0$ and $h_1$ interlace, and $h_1$ has no zeroes on
the intervals $(\gamma_j,\mu_{j+1}^2)$, $j=1,\ldots,k-1$,
therefore $p_0$ also has no zeroes on those intervals. So we
proved the following inequalities
\begin{equation}\label{general.Hurwitz.poly.8}
\mu_j^2<\omega_j<\gamma_j<\mu_{j+1}^2\qquad\text{for}\qquad
j=1,\ldots,k-1.
\end{equation}

Let us consider the intervals $(\mu_j^2,\omega_j)$, $j=1,\ldots,k$
and prove that the polynomial $g(-\sqrt{u})$ has no zeroes on
those interval. Suppose that it is not true. Then there exists a
number $\beta\in(\mu_j^2,\omega_j)$ such that
$g(-\sqrt{\beta})=0$. Then $g_0(\beta)=\sqrt{\beta}g_1(\beta)>0$,
and from the equality
\begin{equation*}
p_0(\beta)=h_0(\beta)g_0(\beta)+\beta
h_1(\beta)g_1(\beta)=g_0(\beta)h(\sqrt{\beta}),
\end{equation*}
we obtain that
$\textrm{sign}\left(p_0(\beta)\right)=\textrm{sign}\left(h(\sqrt{\beta})\right)$.
It contradicts the equality
$\textrm{sign}\left(p_0(\beta)\right)=-\textrm{sign}\left(h(\sqrt{\beta})\right)$,
which follows from~\eqref{general.Hurwitz.poly.8} and from the
equalities
$\textrm{sign}\left(p_0(0)\right)=\textrm{sign}\left(h_0(0)\right)=\textrm{sign}\left(h(0)\right)$
established above.

Thus, we obtain that the polynomial $g$ has real zeroes (en even
number of zeroes, counting multiplicities) only on the intervals
$(-\infty,\sqrt{\omega_k})$, $(-\mu_k,-\sqrt{\omega_{k-1}})$,
\ldots, $(-\mu_2,-\sqrt{\omega_1}),(-\mu_1,0)$. Consequently, the
following equalities hold
$\textrm{sign}\left(g(0)\right)=\textrm{sign}\left(g(-\mu_j)\right)=\textrm{sign}\left(g(-\sqrt{\omega_j})\right)=\textrm{sign}(\alpha_j)>0$
for all $j=1,2,\ldots,k$. So we showed that if $n=2l$ and
$p(0)\neq0$, then the function $\Phi(u)=\dfrac{p_1}{p_0}$ can be
represented as in~\eqref{Theorem.main.general.Hurwitz.proof.13}
with positive~$\alpha_j$ and with exactly $k$ positive poles, as
required.

\vspace{4mm}

Let now $n=2l+1$ and $p(0)\neq0$. Then by definition, the
polynomial $p$ has only $k$ zeroes in the closed right half-plane,
all of which are positive: $0<\mu_1<\mu_2<\ldots<\mu_k$. Moreover,
$p$ has an odd number of zeroes of $p$, counting multiplicities,
on each interval $(-\mu_{k},-\mu_{k-1})$, \ldots,
$(-\mu_3,-\mu_2)$, $(-\mu_2,-\mu_1)$, and an even number of
zeroes, counting multiplicities, on the intervals
$(-\infty,-\mu_{k})$ and $(-\mu_1,0)$. All nonreal zeroes of $p$
are located in the open left half-plane.

As in the case $n=2l$, we represent the polynomial $p$ as a
product $p(z)=g(z)h(z)$, where $h$ is self-interlacing polynomial
(of type~I) of degree $2k-1$ whose positive zeroes are
$\mu_1,\mu_2,\ldots,\mu_k$, and $g$ is Hurwitz stable polynomial
of degree $2(l+1-k)$, which has an even number of zeroes, counting
multiplicities, on each interval
$(-\infty,-\mu_{k}),(-\mu_{k},-\mu_{k-1}),\ldots$,
$(-\mu_2,-\mu_1),(-\mu_1,0)$. In the same way as above one can
prove that the function $\Phi=\dfrac{p_1}{p_0}$ is an
\textit{R}-function. The only difference between this case and the
case of $n=2l$ is the number of positive poles of the function
$\Phi$. In fact, let us consider the function $F$ defined
in~\eqref{general.Hurwitz.poly.6}. In this case, we have $\deg
p_0=l$ if $a_1\neq0$ (or $\deg p_0=l-1$ if $a_1=0$), $\deg
h_1=k-1$,   $\deg g_0=l+1-k$, therefore $\deg(h_1g_0)=l$. Thus,
the number of poles of the function $F$ is equal to (or greater
than) the number of its zeroes.

Let
\begin{equation*}
g(z)=b_0z^{2(l-k)+2}+c_1z^{2(l-k)+1}+\ldots,
\end{equation*}
\begin{equation*}
h(z)=c_0z^{2k-1}+c_1z^{2k-2}+\ldots
\end{equation*}
Then the leading coefficient of the polynomial $h_1$ equals $c_0$,
and the leading coefficient of the polynomial~$g_0$ equals $b_0$
according to~\eqref{poly1.13}--\eqref{poly1.12}. Consequently, the
leading coefficient of the polynomial $h_1g_0$ is equal to
$c_0b_0=a_0$. Moreover, the leading coefficient of the polynomial
$p_0$ is equal to $a_1$ or\footnote{If $a_1=0$.} $a_3$
by~\eqref{poly1.12}, so we have
$\displaystyle\lim_{u\to\infty}F(u)=-\frac{a_1}{a_0}$ or
$-\dfrac{a_3}{a_0}$. The function $F$ is an \textit{R}-function
(that can be proved in the same way as above), so it is decreasing
between its poles, and $F(0)=-\dfrac{h_0(0)}{h_1(0)}>0$
(see~\eqref{self-interlacing.last.coeffs}). Therefore, the minimal
positive zero of the polynomial $p_0$ is less than the minimal
positive zero of the polynomial~$h_1$.

If the condition~\eqref{general.Hurwitz.poly.0.1} holds, then
$\displaystyle\lim_{u\to\infty}F(u)=-\dfrac{a_1}{a_0}<0$, so the
maximal positive zero of $p_0$ (the maximal positive zero of the
function $F$) is greater than the maximal positive zero of the
polynomial $h_1$. This means that the polynomial $p_0$ has exactly
$\deg h_1+1=k$ positive zeroes, so the function $\Phi$ has exactly
$k$ positive poles.

If the condition~\eqref{general.Hurwitz.poly.0.2} holds, then
either $F(u)\to+0$ as $u\to+\infty$ if $a_1=0$ or
$F(u)\to-\dfrac{a_1}{a_0}>0$ as $u\to+\infty$ if $a_1\neq0$, so we
have $\displaystyle\lim_{u\to+\infty}F(u)>0$. This implies that
the maximal positive zero of $h_1$ is greater than the maximal
positive zero of $p_0$, so the polynomial $p_0$ has exactly
$k-1=\deg h_1$ positive zeroes. Consequently, $\Phi$~has exactly
$k-1$ positive poles. Thus, for $n=2l+1$ and $p(0)\neq0$, we
obtain that the function $\Phi=\dfrac{p_1}{p_0}$ can be
represented as in~\eqref{Theorem.main.general.Hurwitz.proof.13}
with positive~$\alpha_j$ and with exactly $k$
(if~\eqref{general.Hurwitz.poly.0.1} holds) or $k-1$
(if~\eqref{general.Hurwitz.poly.0.2} holds) positive poles, as
required.

\vspace{4mm}

Let $n=2l$ and $p(0)=0$. Then by definition, the polynomial $p$
has exactly $k$ zeroes in the \textit{closed} right half-plane:
$0=\mu_1<\mu_2<\ldots<\mu_k$. Moreover, $p$ has  an odd number of
zeroes on each interval $(-\infty,-\mu_{k})$,
$(-\mu_{k},-\mu_{k-1})$, \ldots, $(-\mu_3,-\mu_2)$, $(-\mu_2,0)$,
counting multiplicities. All nonreal zeroes of $p$ are located in
the open left half-plane.

As above, we represent the polynomial $p$ as a product:
$p(z)=g(z)h(z)=p_0(z^2)+zp_1(z^2)$, where
$h(z)=h_0(z^2)+zh_1(z^2)$ is an \textit{almost} self-interlacing
polynomial\footnote{See Section~\ref{subsection:Alomst.SI}.} (of
type~I) of degree $2k$ whose nonnegative zeroes are
$\mu_1,\mu_2,\ldots,\mu_k$, and $g(z)=g_0(z^2)+zg_1(z^2)$ is a
Hurwitz stable polynomial of degree $2(l-k)$, which has an even
number of zeroes, counting multiplicities, on each interval
$(-\infty,-\mu_{k})$, $(-\mu_{k},-\mu_{k-1})$, \ldots,
$(-\mu_3,-\mu_2)$, $(-\mu_2,0)$.

In the same way as above, one can establish that $\Phi$ is an
\textit{R}-function. However, $p_0$ has a root at zero in this
case. The polynomial $h_0$ also has a root at zero, and all its
other roots are positive. Since $\deg h=2k$, we have $\deg
h_0=\deg h_1+1=k$. Besides, the zeroes of $h_0$ interlace the
zeroes of $h_1$, so all zeroes of the polynomial $h_1$ are
positive. Consequently, the \textit{R}-function $F$ defined
in~\eqref{general.Hurwitz.poly.6} has $k-1$ positive poles and
$l-k$ negative poles. It is easy to see that $\deg(g_0h_1)=l-1$
and $\deg p_0=l$, so the maximal positive zero of the polynomial
$p_0$ (and of the function $F$) is greater than the maximal
positive zero of the polynomial $h_1$. Also the minimal positive
of $p_0$ is also greater than the minimal positive zero of $h_1$,
since $F$ is decreasing between its poles and $F(0)=0$. Thus, we
obtain that $p_0$ has exactly $k-1$ positive zeroes, since the
zeroes of $p_0$ interlace the zeroes of $g_0h_1$. Consequently,
the function $\Phi$ has $l-k$ negative poles, $k-1$ positive poles
and a pole at zero.

\vspace{4mm}

If $n=2l+1$ and the polynomial $p(z)=p_0(z^2)+zp_1(z^2)$ has a
root at zero, then in the same as above, one can prove that $\Phi$
is an \textit{R}-function with either exactly $k-1$ positive poles
(if the condition~\eqref{general.Hurwitz.poly.0.1} holds), or
exactly $k-2$ positive poles (if the
condition~\eqref{general.Hurwitz.poly.0.2} holds), since we have
$F(0)=0$ for the \textit{R}-function $F$ defined
in~\eqref{general.Hurwitz.poly.6}.
\end{proof}

Theorem~\ref{Theorem.main.general.Hurwitz} allow us to use
properties of \textit{R}-functions to obtain additional criteria
for polynomials from the class of generalized the Hurwitz
polynomials. In particular, we can generalize Hurwitz and
Li\'enard and Chipart criteria.

Let us consider the Hankel matrix
$S(\Phi)=\|s_{i+j}\|^{\infty}_{i,j=0}$ constructed with the
coefficients of the series~\eqref{app.assoc.function.series}
or~\eqref{app.assoc.function.series.Phi.1}. From
Theorems~\ref{Theorem.main.general.Hurwitz},~\ref{Th.R-function.general.properties},
and~\ref{Th.number.of.negative.poles.of.R-frunction} we obtain the
following criteria for generalized Hurwitz polynomials.
\begin{theorem}\label{Theorem.general.Hurwitz.poly.via.Phi}
The polynomial $p$ defined in~\eqref{general.Hurwitz.poly}, with
$a_1\neq0$ if $n=2l+1$, is generalized Hurwitz if and~only~if
\begin{equation}\label{Theorem.general.Hurwitz.poly.via.Phi.condition.1}
D_j(\Phi)>0, \qquad j=1,\ldots,l.
\end{equation}
Its order $k$ equals:
\begin{itemize}
\item[1\emph{)}] if $n=2l$, then
\begin{equation}\label{general.Hurwitz.poly.1}
k=l-\SCF(1,\widehat{D}_1(\Phi),\widehat{D}_2(\Phi),\ldots,\widehat{D}_r(\Phi)),
\end{equation}
\item[2\emph{)}] if $n=2l+1$, then
\begin{equation}\label{general.Hurwitz.poly.2}
k=l-\SCF(-s_{-1},1,\widehat{D}_1(\Phi),\widehat{D}_2(\Phi),\ldots,\widehat{D}_r(\Phi))+1,
\end{equation}
\end{itemize}
where $r=l$ for $\widehat{D}_{l}(\Phi)\neq0$, and $r=l-1$
for\footnote{By Corollary~\ref{corol.zero.pole}, if
$\widehat{D}_l(\Phi)=0$, then $\widehat{D}_{l-1}(\Phi)\neq0$.}
$\widehat{D}_{l}(\Phi)=0$. Here $s_{-1}=\dfrac{a_0}{a_1}$.

The number of sign changes in the sequence
$1,\widehat{D}_1(\Phi),\widehat{D}_2(\Phi),\ldots,\widehat{D}_r(\Phi)$
must be calculated according to the Frobenius rule provided by
Theorem~\ref{Th.Frobenius}.
\end{theorem}
\begin{proof}
According to Theorem~\ref{Theorem.main.general.Hurwitz}, the
polynomial $p$, with $a_1\neq0$ if $n=2l+1$, is generalized
Hurwitz if and only if its associated function $\Phi$ is an
\textit{R}-function with exactly $l=\left[\dfrac n2\right]$ poles.
But by Theorem~\ref{Th.R-function.general.properties}, this is
equivalent to the
inequalities~\eqref{Theorem.general.Hurwitz.poly.via.Phi.condition.1}.

If $n=2l$, or if $n=2l+1$ with~\eqref{general.Hurwitz.poly.0.1},
then by Theorem~\ref{Theorem.main.general.Hurwitz} the number of
nonnegative poles of the function $\Phi$ is equal to order of the
polynomial $p$. From
Theorem~\ref{Th.number.of.negative.poles.of.R-frunction} we
obtain~\eqref{general.Hurwitz.poly.1}, which is equivalent
to~\eqref{general.Hurwitz.poly.2} if $n=2l+1$
and~\eqref{general.Hurwitz.poly.0.1} holds.

If $n=2l+1$ and $a_0a_1<0$, then by
Theorem~\ref{Theorem.main.general.Hurwitz}, the number of
nonnegative poles of $\Phi$ is equal to order of $p$ minus one.
Thus, by Theorem~\ref{Th.number.of.negative.poles.of.R-frunction}
we have
\begin{equation*}
k-1=l-\SCF(1,\widehat{D}_1(\Phi),\widehat{D}_2(\Phi),\ldots,\widehat{D}_r(\Phi)),
\end{equation*}
that is equivalent to~\eqref{general.Hurwitz.poly.2} in this case.
\end{proof}

In the same way as above one can establish the corresponding
theorem for polynomials of odd degree with $a_1=0$.

\begin{theorem}\label{Theorem.general.Hurwitz.poly.via.Phi.2}
Let the polynomial $p$ of degree $n=2l+1$ be defined
in~\eqref{general.Hurwitz.poly} and let $a_1=0$. The polynomial
$p$ is generalized Hurwitz if and~only~if $a_0a_3<0$ and
\begin{equation}\label{Theorem.general.Hurwitz.poly.via.Phi.2.condition.1}
D_j(\Phi)>0, \qquad j=1,\ldots,l-1.
\end{equation}
Its order $k$ equals:
\begin{equation}\label{general.Hurwitz.poly.3}
k=l-\SCF(1,\widehat{D}_1(\Phi),\widehat{D}_2(\Phi),\ldots,\widehat{D}_r(\Phi)),
\end{equation}
where $r=l-1$ for $\widehat{D}_{l-1}(\Phi)\neq0$, and $r=l-2$
for\footnote{By Corollary~\ref{corol.zero.pole}, if
$\widehat{D}_{l-1}(\Phi)=0$, then $\widehat{D}_{l-2}(\Phi)\neq0$.}
$\widehat{D}_{l-1}(\Phi)=0$.

The number of sign changes in the sequence
$1,\widehat{D}_1(\Phi),\widehat{D}_2(\Phi),\ldots,\widehat{D}_r(\Phi)$
must be calculated according to the Frobenius rule provided by
Theorem~\ref{Th.Frobenius}.
\end{theorem}
\begin{proof}
In fact, let the polynomial $p$ is generalized Hurwitz of
order~$k$. Then by Theorem~\ref{Theorem.main.general.Hurwitz}, the
function $\Phi=\dfrac{p_1}{p_0}$ is an \textit{R}-function. In
particular, this means that $|\deg p_0-\deg p_1|\leqslant1$ (see
e.g.~\cite{ChebotarevMeiman,Holtz_Tyaglov}), therefore, $\deg
p_1=l$ and $\deg p_0=l-1$, since $\deg p=2l+1$ and $a_1=0$ by
assumption. So by Theorems~\ref{Th.R-function.general.properties}
and~\ref{Theorem.main.general.Hurwitz} the function $\Phi$ has the
form
\begin{equation}\label{Theorem.general.Hurwitz.poly.via.Phi.2.proof.1}
\Phi(u)=\dfrac{a_0}{a_3}\,u+\dfrac{a_2a_3-a_0a_5}{a_3^2}+\sum\limits_{j=1}^{k-1}\dfrac{\alpha_j}{u-\omega_j}+\sum\limits_{j=k}^{l-1}\dfrac{\alpha_j}{u+\nu_j},
\end{equation}
where $\dfrac{a_0}{a_3}<0$, $\alpha_j>0$ for $j=1,\ldots,l-1$,
$\omega_j\geqslant0$, $j=1\ldots,k-1$, and $\nu_j>0$ for
$j=k,\ldots,l-1$. The
inequalities~\eqref{Theorem.general.Hurwitz.poly.via.Phi.2.condition.1}
and the formula~\eqref{general.Hurwitz.poly.3} now follow from
Theorems~\ref{Th.R-function.general.properties}
and~\ref{Th.number.of.negative.poles.of.R-frunction}.

Conversely, from the
conditions~\eqref{Theorem.general.Hurwitz.poly.via.Phi.2.condition.1}--\eqref{general.Hurwitz.poly.3},
$a_1=0$ and from the inequality $a_0a_3<0$ it follows that the
function $\Phi$ has the
form~\eqref{Theorem.general.Hurwitz.poly.via.Phi.2.proof.1} by
Theorems~\ref{Th.R-function.general.properties}
and~\ref{Th.number.of.negative.poles.of.R-frunction}. Now
Theorem~\ref{Theorem.main.general.Hurwitz} implies that $p$ is a
generalized Hurwitz polynomial of order $k$.
\end{proof}

From Theorems~\ref{Theorem.general.Hurwitz.poly.via.Phi}
and~\ref{Theorem.general.Hurwitz.poly.via.Phi.2} and from the
formul\ae~\eqref{Formulae.Gurwitz.1}--\eqref{Formulae.Gurwitz.2}
and~\eqref{Formulae.Gurwitz.3} we obtain the following theorem,
which is an analogue of Hurwitz stability criterion for
generalized Hurwitz polynomials.
\begin{theorem}[Generalized Hurwitz theorem]\label{Theorem.general.Hurwitz.criterion}
The polynomial $p$ given in~\eqref{general.Hurwitz.poly} is
generalized Hurwitz if and~only~if
\begin{equation}\label{general.Hurwitz.main.det.noneq}
\Delta_{n-1}(p)>0,\ \Delta_{n-3}(p)>0,\ \Delta_{n-5}(p)>0,\ldots
\end{equation}
The order $k$ of the polynomial $p$ equals
\begin{equation}\label{general.Hurwitz.poly.10}
k=\SCF(\Delta_{n}(p),\Delta_{n-2}(p),\ldots,1)\qquad\text{if}\quad
p(0)\neq0,
\end{equation}
or
\begin{equation}\label{general.Hurwitz.poly.11}
k=\SCF(\Delta_{n-2}(p),\Delta_{n-4}(p),\ldots,1)+1\qquad\text{if}\quad
p(0)=0.
\end{equation}
Here the number of sign changes in~\eqref{general.Hurwitz.poly.10}
and~\eqref{general.Hurwitz.poly.11} must be calculated according
to the Frobenius rule provided by Theorem~\ref{Th.Frobenius}.
\end{theorem}
\begin{proof}
Let the real polynomial $p$ be of degree $n$ such that the
coefficient $a_1$ is nonzero for $n=2l+1$. By
Theorem~\ref{Theorem.general.Hurwitz.poly.via.Phi}, the polynomial
$p$ is generalized Hurwitz if and only if the
inequalities~\eqref{Theorem.general.Hurwitz.poly.via.Phi.condition.1}
hold. These inequalities are equivalent to the
inequalities~\eqref{general.Hurwitz.main.det.noneq} according to
the
formul\ae~\eqref{Formulae.Gurwitz.1}--\eqref{Formulae.Gurwitz.2}.

If the polynomial $p$ is of odd degree and $a_1=0$, then by
Theorem~\ref{Theorem.general.Hurwitz.poly.via.Phi.2}, the
polynomial $p$ is generalized Hurwitz if and only if the
inequalities~\eqref{Theorem.general.Hurwitz.poly.via.Phi.2.condition.1}
and $-a_0a_3>0$ hold. Those inequalities are equivalent to the
inequalities~\eqref{general.Hurwitz.main.det.noneq} according to
the
formul\ae~\eqref{Formulae.Gurwitz.3}--\eqref{Formulae.Gurwitz.333}.

We now prove that order of the polynomial can be calculated by the
formul\ae~\eqref{general.Hurwitz.poly.10}--\eqref{general.Hurwitz.poly.11}.

\vspace{2mm}

Let $n=2l$, then from~\eqref{general.Hurwitz.poly.1}
and~\eqref{Formulae.Gurwitz.1} we obtain
\begin{equation*}\label{Theorem.general.Hurwitz.criterion.proof.1}
\begin{array}{c}
k=l-\SCF(1,\widehat{D}_1(\Phi),\widehat{D}_2(\Phi),\ldots,\widehat{D}_r(\Phi))=l-\SCF(1,-\Delta_2(p),\Delta_4(p),\ldots,(-1)^r\Delta_{2r}(p))=\\
 \\
=\SCF(1,\Delta_2(p),\Delta_4(p),\ldots,\Delta_{2r}(p))+l-r.
\end{array}
\end{equation*}
This is exactly~\eqref{general.Hurwitz.poly.10} for $r=l$ and
exactly~\eqref{general.Hurwitz.poly.11} for $r=l-1$.

\vspace{2mm}

Let now $n=2l+1$ and $a_0a_1>0$. The
formul\ae~\eqref{general.Hurwitz.poly.2}
and~\eqref{Formulae.Gurwitz.1} imply
\begin{equation}\label{Theorem.general.Hurwitz.criterion.proof.2}
\begin{array}{c}
k=l-\SCF(-s_{-1},1,\widehat{D}_1(\Phi),\widehat{D}_2(\Phi),\ldots,\widehat{D}_r(\Phi))+1=\\
 \\
=l+1-\SCF(-s_{-1},1)+\SCF(1,\widehat{D}_1(\Phi),\widehat{D}_2(\Phi),\ldots,\widehat{D}_r(\Phi))=\\
 \\
=l-\SCF(1,-a_1^{-3}\Delta_3(p),a_1^{-5}\Delta_5(p),\ldots,(-1)^ra_1^{-2r-1}\Delta_{2r+1}(p))=\\
 \\
=\SCF(1,\Delta_1(p),\Delta_3(p),\Delta_5(p),\ldots,\Delta_{2r+1}(p))+l-r,
\end{array}
\end{equation}
since $a_1=\Delta_1(p)>0$. One can see
that~\eqref{Theorem.general.Hurwitz.criterion.proof.2}
implies~\eqref{general.Hurwitz.poly.10} for $r=l$
and~\eqref{general.Hurwitz.poly.11} for $r=l-1$.

\vspace{2mm}

Let $n=2l+1$ and $a_0a_1<0$. The
formul\ae~\eqref{general.Hurwitz.poly.2}
and~\eqref{Formulae.Gurwitz.1} imply
\begin{equation}\label{Theorem.general.Hurwitz.criterion.proof.3}
\begin{array}{c}
k=l-\SCF(-s_{-1},1,\widehat{D}_1(\Phi),\widehat{D}_2(\Phi),\ldots,\widehat{D}_r(\Phi))+1=\\
 \\
=l-\SCF(-s_{-1},1)+\SCF(1,\widehat{D}_1(\Phi),\widehat{D}_2(\Phi),\ldots,\widehat{D}_r(\Phi))+1=\\
 \\
=1+l-\SCF(1,-a_1^{-3}\Delta_3(p),a_1^{-5}\Delta_5(p),\ldots,(-1)^ra_1^{-2r-1}\Delta_{2r+1}(p))=\\
 \\
=1+\SCF(a_1,a_1^{-2}\Delta_3(p),a_1^{-4}\Delta_5(p),\ldots,a_1^{-2r}\Delta_{2r+1}(p))+l-r=\\
 \\
=\SCF(1,\Delta_1(p))+\SCF(a_1,a_1^{-2}\Delta_3(p),a_1^{-4}\Delta_5(p),\ldots,a_1^{-2r}\Delta_{2r+1}(p))+l-r=\\
 \\
=\SCF(1,\Delta_1(p),\Delta_3(p),\Delta_5(p),\ldots,\Delta_{2r+1}(p))+l-r,
\end{array}
\end{equation}
since $a_1=\Delta_1(p)<0$. Now
from~\eqref{Theorem.general.Hurwitz.criterion.proof.3} we
obtain~\eqref{general.Hurwitz.poly.10} for $r=l$
and~\eqref{general.Hurwitz.poly.11} for $r=l-1$.

\vspace{2mm}

Finally, let $n=2l+1$ and $a_1=\Delta_1(p)=0$. The
formul\ae~\eqref{general.Hurwitz.poly.3}
and~\eqref{Formulae.Gurwitz.3}--\eqref{Formulae.Gurwitz.333} yield
\begin{equation}\label{Theorem.general.Hurwitz.criterion.proof.4}
\begin{array}{c}
k=l-\SCF(1,\widehat{D}_1(\Phi),\widehat{D}_2(\Phi),\ldots,\widehat{D}_{r}(\Phi))=\\
 \\
=l-\SCF\left(1,-\dfrac{\Delta_5(p)}{\Delta_3(p)},\dfrac{\Delta_7(p)}{\Delta_3(p)},\ldots,(-1)^{r}\dfrac{\Delta_{2r+3}(p)}{\Delta_3(p)}\right)=\\
 \\
=l-1-r+\SCF(1,0,-a_0a_3^2)+\SCF\left(1,\dfrac{\Delta_5(p)}{\Delta_3(p)},\dfrac{\Delta_7(p)}{\Delta_3(p)},\ldots,\dfrac{\Delta_{2r+3}(p)}{\Delta_3(p)}\right)=\\
\\
=l-1-r+\SCF(1,\Delta_1(p),\Delta_3(p))+\SCF\left(\Delta_3(p),\Delta_5(p),\Delta_7(p),\ldots,\Delta_{2r+3}(p)\right)=\\
 \\
=l-1-r+\SCF(1,\Delta_1(p)\Delta_3(p),\Delta_5(p),\Delta_7(p),\ldots,\Delta_{2r+3}(p)),
\end{array}
\end{equation}
where $r=l-1$ for $p(0)\neq0$, and $r=l-2$ for $p(0)=0$.
Now~\eqref{Theorem.general.Hurwitz.criterion.proof.4} is
equivalent to~\eqref{general.Hurwitz.poly.10} for $r=l-1$ and
to~\eqref{general.Hurwitz.poly.11} for $r=l-2$.
\end{proof}
\begin{remark}
Using the
formul\ae~\eqref{Formulae.Gurwitz.1}--\eqref{Formulae.Gurwitz.2}
and~\eqref{Formulae.Gurwitz.3}--\eqref{Formulae.Gurwitz.333}, one
can easily reformulate
Theorem~\ref{Theorem.general.Hurwitz.criterion} in terms of the
minors $\eta_j(p)$ defined in
Definition~\ref{def.Hurwitz.matrix.infinite.for.poly}.
\end{remark}

Now from Theorems~\ref{Th.generalized.Lienard.Chipart}
and~\ref{Theorem.main.general.Hurwitz} we obtain the following
theorem, which is a generalization of the famous Li\'enard and
Chipart criterion of stability for the class of generalized
Hurwitz polynomials.

\begin{theorem}[Generalized Li\'enard and Chipart theorem]\label{Theorem.Lienard.Chipart.Generalized}
The polynomial $p$ given in~\eqref{general.Hurwitz.poly} is
generalized Hurwitz if and~only~if the
inequalities~\eqref{general.Hurwitz.main.det.noneq} hold. The
order $k$ of the polynomial $p$ equals
\begin{equation}\label{general.Hurwitz.poly.12}
k=v(a_n,a_{n-2},\ldots,1)=v(a_n,a_{n-1},a_{n-3},\ldots,1)\qquad\text{if}\qquad
a_n\neq0,
\end{equation}
or\footnote{If $a_n=0$, then always $a_{n-1}\neq0$, since the
polynomials $p_1$ and $p_0$, the numerator and the denominator of
the \textit{R}-function~$\Phi$, have no common roots.}
\begin{equation}\label{general.Hurwitz.poly.13}
k=v(a_{n-2},a_{n-4},\ldots,1)+1=v(a_{n-1},a_{n-3},\ldots,1)+1,\qquad\text{if}\qquad
a_n=0,
\end{equation}
where $v$ is the number of strong sign changes (see
Definition~\ref{Def.strong.sign.changes}).
\end{theorem}
\begin{proof}
The fact that the polynomial $p$ is generalized Hurwitz if and
only if the inequalities~\eqref{general.Hurwitz.main.det.noneq}
hold follow from Theorem~\ref{Theorem.general.Hurwitz.criterion}.
Now we prove the
formul\ae~\eqref{general.Hurwitz.poly.12}--\eqref{general.Hurwitz.poly.13}.

\vspace{3mm}

Let $n=2l$. By Theorem~\ref{Theorem.main.general.Hurwitz}, $p$ is
a generalized Hurwitz polynomial if and only if the function
\begin{equation*}\label{Theorem.Lienard.Chipart.Generalized.proof.1.5}
\Phi(u)=\dfrac{p_1(u)}{p_0(u)}=\dfrac{a_1u^{l-1}+a_3u^{l-2}+\ldots+a_{2l-3}u+a_{2l-1}}{a_0u^{l}+a_2u^{l-1}+\ldots+a_{2l-2}u+a_{2l}}
\end{equation*}
is an \textit{R}-function. The number of its nonnegative poles is
equal to order $k$ of the polynomial $p$. Let $\Phi$ has exactly
$m$ positive poles\footnote{Clearly, $m=k$ if $a_n\neq0$, and
$m=k-1$ if $a_n=0$.}. We prove that
\begin{equation}\label{Theorem.Lienard.Chipart.Generalized.proof.1}
m=v(a_n,a_{n-2},\ldots,1).
\end{equation}
It will be equivalent to to the first equalities
in~\eqref{general.Hurwitz.poly.12}--\eqref{general.Hurwitz.poly.13},
because $v$ is the number of strong sign changes, so
$v(a_n,a_{n-2},\ldots,1)=v(a_{n-2},a_{n-4},\ldots,1)$ if $a_n=0$.
Now Theorem~\ref{Th.generalized.Lienard.Chipart} implies
\begin{equation*}\label{Theorem.Lienard.Chipart.Generalized.proof.2}
m=v(a_{2l},a_{2l-2},\ldots,a_0)=v(a_{2l},a_{2l-2},\ldots,a_0)+v(a_0,1)=v(a_{2l},a_{2l-2},\ldots,a_0,1),
\end{equation*}
which is
exactly~\eqref{Theorem.Lienard.Chipart.Generalized.proof.1}.

\vspace{3mm}

Let us now consider the function
\begin{equation*}\label{Theorem.Lienard.Chipart.Generalized.proof.4.5}
\Psi(u)=-\dfrac1{\Phi(u)}=-\dfrac{p_0(u)}{p_1(u)}=-\dfrac{a_0u^{l}+a_2u^{l-1}+\ldots+a_{2l-2}u+a_{2l}}{a_1u^{l-1}+a_3u^{l-2}+\ldots+a_{2l-3}u+a_{2l-1}}.
\end{equation*}
This function is also an \textit{R}-function by
Theorem~\ref{Theorem.properties.R-functions}, since $\Phi$ is an
\textit{R}-function with exactly $k$ nonnegative poles, say
$0\leqslant\omega_1<\ldots,\omega_k$. Note that $\Phi(u)>0$ for
$u\in(\omega_k,+\infty)$, so it can have at least $k-1$ and at
most $k$ positive zeroes. There can be only four possibilities:
\begin{itemize}
\item[I.] $\Phi(0)=\dfrac{a_{2l-1}}{a_{2l}}<0$. Then
$v(a_{2l},a_{2l-1})=1$, and $\Phi$ has no zeroes on the interval
$(0,\omega_1)$, since it is decreasing on this interval and
$\Phi(u)\to-\infty$ as $u\nearrow\omega_1$. Therefore, $\Phi$ has
exactly $k-1$ positive zeroes, so $\Psi$ has exactly $k-1$
positive poles. By Theorem~\ref{Th.generalized.Lienard.Chipart}
and Remark~\ref{Th.generalized.Lienard.Chipart.remark} applied to
the function $\Psi$, we have
\begin{equation*}\label{Theorem.Lienard.Chipart.Generalized.proof.5}
\begin{array}{c}
k-1=v(a_{2l-1},a_{2l-3},\ldots,a_1)=-1+v(a_{2l},a_{2l-1})+v(a_{2l-1},a_{2l-3},\ldots,a_1,1)+v(a_1,1)=\\
 \\
=-1+v(a_{2l},a_{2l-1},a_{2l-3},\ldots,a_1,1),
\end{array}
\end{equation*}
which is exactly
\begin{equation}\label{Theorem.Lienard.Chipart.Generalized.proof.6}
k=v(a_n,a_{n-1},a_{n-3},\ldots,1).
\end{equation}
Here we used the inequality $\Delta_1(p)=a_1>0$, which follows
from the inequalities~\eqref{general.Hurwitz.main.det.noneq} for
$n=2l$.
\item[II.]  $\Phi(0)=\dfrac{a_{2l-1}}{a_{2l}}>0$. Then
$v(a_{2l},a_{2l-1})=0$, and $\Phi$ has exactly $k$ positive
zeroes, so $\Psi$ has exactly $k$ positive poles. Thus, from
Theorem~\ref{Th.generalized.Lienard.Chipart} and
Remark~\ref{Th.generalized.Lienard.Chipart.remark} applied to the
function $\Psi$, we obtain
\begin{equation*}\label{Theorem.Lienard.Chipart.Generalized.proof.7}
\begin{array}{c}
k=v(a_{2l-1},a_{2l-3},\ldots,a_1)=v(a_{2l},a_{2l-1})+v(a_{2l-1},a_{2l-3},\ldots,a_1)+v(a_1,1)=\\
 \\
=v(a_{2l},a_{2l-1},a_{2l-3},\ldots,a_1,1),
\end{array}
\end{equation*}
which is
exactly~\eqref{Theorem.Lienard.Chipart.Generalized.proof.6}.
\item[III.]  $\Phi(0)=0$. Then
$\Psi$ has a pole at zero, that is, $a_{2l-1}=0$. Since $\Psi$ is
an \textit{R}-function, it has positive residues at each its pole,
so $\lim\limits_{u\to0}u\Psi(u)=-\dfrac{a_{2l}}{a_{2l-3}}>0$.
Thus, $v(a_{2l},a_{2l-3})=v(a_{2l},a_{2l-1},a_{2l-3})=1$. As in
Case~I, $\Phi$ has exactly $k-1$ positive zeroes, so the function
$\Psi$ has exactly $k-1$ positive poles. Therefore, by
Theorem~\ref{Th.generalized.Lienard.Chipart} and
Remark~\ref{Th.generalized.Lienard.Chipart.remark} applied to the
function $\Psi$, we have
\begin{equation*}\label{Theorem.Lienard.Chipart.Generalized.proof.8}
\begin{array}{c}
k-1=v(a_{2l-3},a_{2l-5},\ldots,a_1)=-1+v(a_{2l},a_{2l-1},a_{2l-3})+v(a_{2l-3},a_{2l-5},\ldots,a_1)+v(a_1,1)=\\
 \\
=-1+v(a_{2l},a_{2l-1},a_{2l-3},\ldots,a_1,1),
\end{array}
\end{equation*}
which is
precisely~\eqref{Theorem.Lienard.Chipart.Generalized.proof.6}.
\item[IV.]  $\Psi(0)=0$. Then $p(0)=a_n=a_{2l}=0$, and the function $\Phi$
has a pole at zero, so it has exactly $k-1$ positive zeroes. Thus,
the function $\Psi$ has exactly $k-1$ positive poles, and by
Theorem~\ref{Th.generalized.Lienard.Chipart} and
Remark~\ref{Th.generalized.Lienard.Chipart.remark} applied to the
function $\Psi$, we have
\begin{equation*}\label{Theorem.Lienard.Chipart.Generalized.proof.9}
\begin{array}{c}
k-1=v(a_{2l-1},a_{2l-3},\ldots,a_1)=v(a_{2l-1},a_{2l-3},\ldots,a_1,1),
\end{array}
\end{equation*}
therefore,
\begin{equation}\label{Theorem.Lienard.Chipart.Generalized.proof.10}
k=v(a_{n-1},a_{n-3},\ldots,1)+1.
\end{equation}
\end{itemize}

\vspace{4mm}

Let now $n=2l+1$. By Theorem~\ref{Theorem.main.general.Hurwitz},
$p$ is a generalized Hurwitz polynomial if and only if $\Phi$ is
an \textit{R}-function with exactly $l$ zeroes and $l$ or $l-1$
poles. If $k$ is order of the polynomial $p$, then $\Phi$ has
exactly $k$ nonnegative poles for $a_1>0$ and exactly $k-1$
nonnegative poles for $a_1\leqslant0$.

At first, we establish the
formula~\eqref{Theorem.Lienard.Chipart.Generalized.proof.1}, where
$m$ is the number of positive poles of the function $\Phi$. As we
mentioned above, this formula is equivalent to the first
equalities
in~\eqref{general.Hurwitz.poly.12}--\eqref{general.Hurwitz.poly.13}.

If $a_1>0$, then from Theorem~\ref{Th.generalized.Lienard.Chipart}
applied to the function $\Phi$ we obtain
\begin{equation*}\label{Theorem.Lienard.Chipart.Generalized.proof.11}
m=v(a_{2l+1},a_{2l-1},\ldots,a_1)=v(a_{2l+1},a_{2l-1},\ldots,a_1,1).
\end{equation*}
This is
exactly~\eqref{Theorem.Lienard.Chipart.Generalized.proof.1}.

Let $a_1\leqslant0$. As we established above, if $a_1=0$, then
$a_3<0$. Therefore, for $a_1\leqslant0$, we always have
$v(a_3,a_1)=v(a_3,a_1,1)-1$. Thus,
Theorem~\ref{Th.generalized.Lienard.Chipart} implies
\begin{equation*}\label{Theorem.Lienard.Chipart.Generalized.proof.12}
m=v(a_{2l+1},a_{2l-1},\ldots,a_1)=v(a_{2l+1},a_{2l-1},\ldots,a_1,1)-1,
\end{equation*}
which is
again~\eqref{Theorem.Lienard.Chipart.Generalized.proof.1}, since
for $a_1\leqslant0$, $m=k-1$ if $p(0)\neq0$, and $m=k-2$ if
$p(0)=0$.

\vspace{3mm}

As above, consider the function
\begin{equation*}\label{Theorem.Lienard.Chipart.Generalized.proof.12.5}
\Psi(u)=-\dfrac1{\Phi(u)}=-\dfrac{p_0(u)}{p_1(u)}=-\dfrac{a_1u^{l}+a_3u^{l-1}+\ldots+a_{2l-1}u+a_{2l+1}}{a_0u^{l}+a_2u^{l-1}+\ldots+a_{2l-2}u+a_{2l}}.
\end{equation*}

This function is also an \textit{R}-function by
Theorem~\ref{Theorem.properties.R-functions}. As above, there can
be only four possibilities:
\begin{itemize}
\item[I.] $\Phi(0)=\dfrac{a_{2l}}{a_{2l+1}}<0$, so
$v(a_{2l+1},a_{2l})=1$. If $a_1>0$, then by
Theorem~\ref{Theorem.main.general.Hurwitz}, $\Phi$ has exactly $k$
positive poles, say $0<\omega_1<\ldots<\omega_k$. Since
$\Phi(u)\to\dfrac{a_0}{a_1}>0$ as $u\to+\infty$ and $\Phi$ is
decreasing between its poles, the function $\Phi$ has no zeroes on
the interval $(\omega_k,+\infty)$. Also $\Phi$ has no zeroes on
the interval $(0,\omega_1)$, since $\Phi(0)<0$ by assumption.
Thus, $\Phi$ has exactly $k-1$ positive zeroes, so $\Psi$ has
exactly $k-1$ positive poles.

If $a_1\leqslant0$, then by
Theorem~\ref{Theorem.main.general.Hurwitz}, the function $\Phi$
has exactly $k-1$ positive poles, say
$0<\omega_1<\ldots<\omega_{k-1}$. Since $\Phi(u)$ is decreasing
between its poles, it is negative for sufficiently large
positive~$u$. Moreover, $\Phi(u)\nearrow+\infty$ as
$u\searrow\omega_{k-1}$. Thus, $\Phi$ has a zero on the interval
$(\omega_{k-1},+\infty)$, but it has no zeroes $(0,\omega_1)$,
since $\Phi(0)<0$ by assumption. Consequently, $\Phi$ has exactly
$k-1$ positive zeroes, so $\Psi$ has exactly $k-1$ positive poles.

Thus, we showed that the function $\Psi$ has exactly $k-1$
positive poles regardless of the sign of $a_1$. Now from
Theorem~\ref{Th.generalized.Lienard.Chipart} and
Remark~\ref{Th.generalized.Lienard.Chipart.remark} applied to the
function $\Psi$, we obtain
\begin{equation*}\label{Theorem.Lienard.Chipart.Generalized.proof.13}
\begin{array}{c}
k-1=v(a_{2l},a_{2l-2},\ldots,a_0)=-1+v(a_{2l+1},a_{2l})+v(a_{2l},a_{2l-2},\ldots,a_0,1)=\\
 \\
=-1+v(a_{2l+1},a_{2l},a_{2l-2},\ldots,a_0,1),
\end{array}
\end{equation*}
which is
exactly~\eqref{Theorem.Lienard.Chipart.Generalized.proof.6}.
\item[II.]  $\Phi(0)=\dfrac{a_{2l}}{a_{2l+1}}>0$, so
$v(a_{2l+1},a_{2l})=0$. In the same way as above, one can show
that the function $\Psi$ has exactly $k$ positive poles regardless
of the sign of $a_1$. Thus, from
Theorem~\ref{Th.generalized.Lienard.Chipart} and
Remark~\ref{Th.generalized.Lienard.Chipart.remark} applied to the
function $\Psi$, we obtain
\begin{equation*}\label{Theorem.Lienard.Chipart.Generalized.proof.15}
\begin{array}{c}
k=v(a_{2l},a_{2l-2},\ldots,a_0)=v(a_{2l+1},a_{2l})+v(a_{2l},a_{2l-2},\ldots,a_0,1)=\\
 \\
=v(a_{2l+1},a_{2l},a_{2l-2},\ldots,a_0,1).
\end{array}
\end{equation*}
So the formula~\eqref{Theorem.Lienard.Chipart.Generalized.proof.6}
is also valid in this case.
\item[III.]  $\Phi(0)=0$. Then
$\Psi$ has a pole at zero, that is, $a_{2l}=0$. Since $\Psi$ is an
\textit{R}-function, it has positive residues at each its pole, so
$\lim\limits_{u\to0}u\Psi(u)=-\dfrac{a_{2l+1}}{a_{2l-2}}>0$. Thus,
$v(a_{2l+1},a_{2l-2})=v(a_{2l+1},a_{2l},a_{2l-2})=1$. As above,
one can show that $\Psi$ has exactly $k-1$ positive poles
regardless of the sign of $a_1$. Therefore, by
Theorem~\ref{Th.generalized.Lienard.Chipart} and
Remark~\ref{Th.generalized.Lienard.Chipart.remark} applied to the
function $\Psi$, we have
\begin{equation*}\label{Theorem.Lienard.Chipart.Generalized.proof.16}
\begin{array}{c}
k-1=v(a_{2l-2},a_{2l-4},\ldots,a_0)=-1+v(a_{2l+1},a_{2l},a_{2l-2})+v(a_{2l-2},a_{2l-4},\ldots,a_0,1)=\\
 \\
=-1+v(a_{2l+1},a_{2l},a_{2l-2},\ldots,a_0,1),
\end{array}
\end{equation*}
which is
precisely~\eqref{Theorem.Lienard.Chipart.Generalized.proof.6}.
\item[IV.]  $\Psi(0)=0$. Then $p(0)=a_n=a_{2l+1}=0$, and the function $\Phi$
has a pole at zero. In the same way as above, one can prove that
the function $\Psi$ has exactly $k-1$ positive poles regardless of
the sign of~$a_1$, so Theorem~\ref{Th.generalized.Lienard.Chipart}
and Remark~\ref{Th.generalized.Lienard.Chipart.remark} applied to
the function $\Psi$ imply
\begin{equation*}\label{Theorem.Lienard.Chipart.Generalized.proof.17}
\begin{array}{c}
k-1=v(a_{2l},a_{2l-2},\ldots,a_0),
\end{array}
\end{equation*}
which is equivalent
to~\eqref{Theorem.Lienard.Chipart.Generalized.proof.10}, since
$a_0>0$.
\end{itemize}
\end{proof}

This theorem implies a necessary condition for polynomials to be
generalized Hurwitz. This condition plays a role of Stodola's
necessary condition for Hurwitz stable polynomials
(Theorem~\ref{Th.Stodola.necessary.condition.Hurwitz}).
\begin{corol}\label{Corol.Stodola.generalized.Hurwitz}
If the polynomial
\begin{equation}\label{Corol.Stodola.generalized.Hurwitz.poly}
p(z)=a_0z^n+a_1z^{n-1}+\dots+a_n,\qquad
a_1,\dots,a_n\in\mathbb{R},\ a_0>0,
\end{equation}
is generalized Hurwitz (of type~I) of order $k$, then
\begin{equation*}\label{Corol.Stodola.necessary.condition.generalized.Hurwitz}
k=v(a_n,a_{n-2},\ldots,1)=v(a_n,a_{n-1},a_{n-3},\ldots,1)\qquad\text{if}\quad
p(0)\neq0,
\end{equation*}
or
\begin{equation*}\label{Corol.Stodola.necessary.condition.generalized.Hurwitz.2}
k-1=v(a_n,a_{n-2},\ldots,1)=v(a_n,a_{n-1},a_{n-3},\ldots,1)\qquad\text{if}\quad
p(0)=0,
\end{equation*}
\end{corol}

\vspace{2mm}

This corollary implies
Theorem~\ref{Th.Stodola.necessary.condition.Hurwitz} (Stodola's
theorem) for $k=0$, and it implies
Theorem~\ref{Th.Stodola.necessary.condition.self-interlacing} for
$k=\left[\dfrac{n+1}2\right]$ and $a_n\neq0$.

Note that G.\,F.\,Korsakov~\cite{Korsakov} made an attempt to
prove Theorem~\ref{Theorem.Lienard.Chipart.Generalized}. However,
his methods did not allow him to prove the whole theorem. He
proved only sufficiency, that is, he proved that if the
inequalities~\eqref{general.Hurwitz.main.det.noneq} hold for the
polynomial~\eqref{Corol.Stodola.generalized.Hurwitz.poly}, then
$p$ is generalized Hurwitz, and the number $k$ of its positive
zeroes can be found by the formul\ae
\begin{itemize}
\item for $n=2l$,
\begin{equation}\label{Korsakov.1}
k=v(a_n,a_{n-2},\ldots,a_{2},a_0);
\end{equation}
\item for $n=2l+1$,
\begin{equation}\label{Korsakov.2}
k=v(a_n,a_{n-2},\ldots,a_3,a_{1},a_0).
\end{equation}
\end{itemize}
Since in~\cite{Korsakov}, it was assumed that $a_n\neq0$, the
number $k$ in~\eqref{Korsakov.1}--\eqref{Korsakov.2} is exactly
order of the polynomial~$p$. Nevertheless, the generalized Hurwitz
polynomials with a root at zero are also partially described
in~\cite{Korsakov}. In fact, there was proved that if the
following inequalities hold for the polynomial $p$ defined
in~\eqref{Corol.Stodola.generalized.Hurwitz.poly}
\begin{equation}\label{Korsakov.3}
\Delta_{n}(p)>0,\ \Delta_{n-2}(p)>0,\ \Delta_{n-4}(p)>0,\ldots,
\end{equation}
then $p$ has exactly $k$ (defined
by~\eqref{Korsakov.1}--\eqref{Korsakov.2}) zeroes on the positive
real half-axis, and other its zeroes are located in the open left
half-plane. Generally speaking, the polynomial $p$ satisfied the
inequalities~\eqref{Korsakov.3} is not a generalized Hurwitz
polynomial. However, if $p$ satisfies~\eqref{Korsakov.3}, then the
polynomial $q(z)=zp(z)$ of degree $n+1$ is a generalized Hurwitz
polynomial of order $k+1$. Indeed, it is easy to see that the
following equalities are true:
\begin{equation*}\label{Korsakov.4}
\Delta_i(q)=\Delta_i(p),\qquad i=1,\ldots,n.
\end{equation*}
So the inequalities~\eqref{Korsakov.3} for the polynomial $p$ of
degree $n$ are exactly the
inequalities~\eqref{general.Hurwitz.main.det.noneq} hold  for the
polynomial $q$ of degree $n+1$, so $q$ is a generalized Hurwitz
polynomial having $k+1$ nonnegative zeroes.

\vspace{3mm}

At last, we establish a relation between generalized Hurwitz
polynomials and continued fractions of Stieltjes type (if any).
\begin{theorem}\label{Theorem.general.Hurwitz.Stieltjes.cont.frac.criteria}
If the polynomial $p$ of degree $n$ defined
in~\eqref{Corol.Stodola.generalized.Hurwitz.poly} is generalized
Hurwitz \emph{(}of type~I\emph{)} of order~$k$ and 
 $\Delta_{n-2j}(p)\neq0$, $i=1,\ldots,l$, where $\Delta_0(p)\equiv1$,
then the function~$\Phi$ has the following Stieltjes type
continued fraction expansion:
\begin{equation}\label{Stieltjes.fraction.for.general.Hurwitz}
\Phi(u)=\dfrac{p_1(u)}{p_0(u)}=c_0+\dfrac1{c_1u+\cfrac1{c_2+\cfrac1{c_{3}u+\cfrac1{\ddots+\cfrac1{T}}}}},\qquad
T=
\begin{cases}
&c_{2l}\qquad\text{if}\;\;\; p(0)\neq0,\\
&c_{2l-1}u\quad\text{if}\;\;\; p(0)=0,
\end{cases}
\end{equation}
where $c_0=0$ if $n$ is even, and $c_0\neq0$ if $n$ is odd, and
\begin{equation}\label{Stieltjes.fraction.for.general.Hurwitz.2}
c_{2j-1}>0,\qquad j=1,\ldots,l.
\end{equation}
Moreover, if $p(0)\neq0$, then the number of negative coefficients
$c_{2j}$, $j=0,1,\ldots,l$, equals $k$, and if $p(0)=0$, then the
number of negative coefficients $c_{2j}$, $j=0,1,\ldots,l-1$,
equals $k-1$.

\vspace{2mm}

Conversely, if for given polynomial $p$ of degree $n$, its
associated function $\Phi$ has a Stieltjes continued
fraction~\eqref{Stieltjes.fraction.for.general.Hurwitz}--\eqref{Stieltjes.fraction.for.general.Hurwitz.2},
then $p$ is a generalized Hurwitz polynomial \emph{(}with
$\Delta_{n-2j}(p)\neq0$, $i=1,\ldots,l$\emph{)}, where~$l$ defined
in~\eqref{floor.poly.degree}. Its order equals the number of
negative coefficients $c_{2j}$, $j=0,1,\ldots,l$, if $p(0)\neq0$,
or the number of negative coefficients $c_{2j}$,
$j=0,1,\ldots,l-1$, plus one if $p(0)=0$.
\end{theorem}
\begin{proof}
The theorem follows from
Theorem~\ref{Theorem.general.Hurwitz.criterion} and from the
formul\ae~\eqref{coeff.Stieltjes.fraction.main.formula.quasi-stability.even.degree}--\eqref{coeff.Stieltjes.fraction.main.formula.quasi-stability.odd.degree}î
\end{proof}

Theorem~\ref{Theorem.general.Hurwitz.Stieltjes.cont.frac.criteria}
does not cover the case of $a_1=0$. The following theorem fills
this gap.
\begin{theorem}\label{Theorem.general.Hurwitz.Stieltjes.cont.frac.criteria.2}
If the polynomial $p$ of degree $n=2l+1$ defined
in~\eqref{Corol.Stodola.generalized.Hurwitz.poly} with $a_1=0$ is
generalized Hurwitz \emph{(}of type~I\emph{)} of order~$k$ and
 $\Delta_{2j+1}(p)\neq0$, $i=1,\ldots,l-1$,
then the function~$\Phi$ has the following Stieltjes continued
fraction expansion:
\begin{equation}\label{Stieltjes.fraction.for.general.Hurwitz.00}
\Phi(u)=\dfrac{p_1(u)}{p_0(u)}=-c_{-1}u+c_0+\dfrac1{c_1u+\cfrac1{c_2+\cfrac1{c_{3}u+\cfrac1{\ddots+\cfrac1{T}}}}},\qquad
T=
\begin{cases}
&c_{2l-2}\qquad\text{if}\;\;\; p(0)\neq0,\\
&c_{2l-3}u\quad\text{if}\;\;\; p(0)=0,
\end{cases}
\end{equation}
where $c_0\in\mathbb{R}$, and
\begin{equation}\label{Stieltjes.fraction.for.general.Hurwitz.22}
c_{2j-1}>0,\qquad j=0,1,\ldots,l-1.
\end{equation}
Moreover, if $p(0)\neq0$, then the number of negative coefficients
$c_{2j}$, $j=1,\ldots,l-1$, equals $k-1$. But if $p(0)=0$, then
the number of negative coefficients $c_{2j}$, $j=0,1,\ldots,l-2$,
equals $k-2$.

\vspace{2mm}

Conversely, if for given polynomial $p$ of degree $n=2l+1$, its
associated function $\Phi$ has a Stieltjes continued
fraction~\eqref{Stieltjes.fraction.for.general.Hurwitz.00}--\eqref{Stieltjes.fraction.for.general.Hurwitz.22},
then $p$ is a generalized Hurwitz polynomial \emph{(}with
$\Delta_{2j+1}(p)\neq0$, $i=1,\ldots,l-1$ and
$\Delta_1(p)=a_1=0$\emph{)}. Its order equals the number of
negative coefficients $c_{2j}$, $j=1,\ldots,l-1$, plus one if
$p(0)\neq0$, or the number of negative coefficients $c_{2j}$,
$j=1,\ldots,l-2$, plus two if $p(0)=0$.
\end{theorem}
\begin{proof}
Let the polynomial $p$ of degree $n=2l+1$ with $a_1=0$ be
generalized Hurwitz (of type~I) of order~$k$ and
$\Delta_{2j+1}(p)\neq0$, $i=1,\ldots,l-1$. By
Theorem~\ref{Theorem.general.Hurwitz.poly.via.Phi.2} and by
formul\ae~\eqref{Formulae.Gurwitz.3}, we obtain $a_0a_3<0$ and
\begin{equation}\label{Theorem.general.Hurwitz.Stieltjes.cont.frac.criteria.2.proof.1}
\begin{cases}
& D_j(\Phi)>0,\qquad j=1,\ldots,l-1,\\
& \widehat{D}_j(\Phi)\neq0\qquad j=1,\ldots,l-2.
\end{cases}
\end{equation}
It is easy to see that for the function $F(u):=\Phi(u)+c_{-1}u$,
where $c_{-1}=-\dfrac{a_0}{a_3}$, the following equalities hold
\begin{equation}\label{Theorem.general.Hurwitz.Stieltjes.cont.frac.criteria.2.proof.2}
\begin{cases}
& D_j(F)=D_j(\Phi),\\
& \widehat{D}_j(F)=\widehat{D}_j(\Phi)\neq0,
\end{cases}
\qquad j=1,\ldots,l-1,
\end{equation}
so the function $F$ has a Stiltjes type continued
fraction~\eqref{Stieltjes.fraction.1} with $c_{2j-1}>0$,
$j=1,\ldots,l-1$, that follows
from~\eqref{Theorem.general.Hurwitz.Stieltjes.cont.frac.criteria.2.proof.1}--\eqref{Theorem.general.Hurwitz.Stieltjes.cont.frac.criteria.2.proof.2}
and~\eqref{odd.coeff.Stieltjes.fraction.main.formula}. Also
from~\eqref{general.Hurwitz.poly.3},~\eqref{Theorem.general.Hurwitz.Stieltjes.cont.frac.criteria.2.proof.1}--\eqref{Theorem.general.Hurwitz.Stieltjes.cont.frac.criteria.2.proof.2}
and~\eqref{even.coeff.Stieltjes.fraction.main.formula} we obtain
that the number of negative coefficients $c_{2j}$,
$j=1,\ldots,l-1$, equals $k-1$ if $p(0)\neq0$, or $k-2$ if
$p(0)=0$.

\vspace{2mm}

The converse assertion of the theorem follows from the
formul\ae~\eqref{Theorem.general.Hurwitz.Stieltjes.cont.frac.criteria.2.proof.1}--\eqref{Theorem.general.Hurwitz.Stieltjes.cont.frac.criteria.2.proof.2},
from Theorem~\ref{Th.Stieltjes.fraction.criterion} applied to the
function $F(u):=\Phi(u)-\dfrac{a_0}{a_3}u$, and from
Theorem~\ref{Theorem.general.Hurwitz.poly.via.Phi.2}.
\end{proof}

Theorem~\ref{Corol.differentiation.of.R-functions} implies the
following theorem.
\begin{theorem}\label{Theorem.general.Hurwitz.with.differentiation}
Let the polynomial $p$ of degree $n\geqslant2$ as
in~\eqref{Corol.Stodola.generalized.Hurwitz.poly} be generalized
Hurwitz polynomial. Then all polynomials
\begin{equation*}\label{Theorem.general.Hurwitz.with.differentiation.condition}
p_j(z)=\sum\limits_{i=0}^{n-2j}\left[\dfrac{n-i}2\right]\left(\left[\dfrac{n-i}2\right]-1\right)\cdots\left(\left[\dfrac{n-i}2\right]+j-1\right)a_iz^{n-2j-i},\quad
k=1,\ldots,\left[\dfrac{n}2\right]-1,
\end{equation*}
are also generalized Hurwitz.
\end{theorem}

\subsection{Application to bifurcation theory. Polynomials dependent on
parameters.}\label{subsection:application.bifurcation}

\hspace{4mm} Let us consider a system of ordinary differential
equations dependent on a real parameter, say $\alpha$:
\begin{equation}\label{diff.sys}
\begin{cases}
&\dot{x}(t)=F(x(t),\alpha)=A(\alpha)x(t)+o(x(t)),\\
&x(t_0)=x_0
\end{cases}
\qquad F(x),\,\,x\in\mathbb{R}^n,
\end{equation}
where the $n\times n$ matrix $A(\alpha)$ is the linear part of the
function $F(x,\alpha)$, and $F(0,\alpha)\equiv0$. Suppose that for
some value of the parameter $\alpha_0$ all the eigenvalues of the
matrix $A(\alpha_0)$ are located in the open left half-plane of
the complex plane. In this case, the stability theory
says~\cite{Lyapunoff,Marsden&McCracen,Iooss_Joseph} that the zero
solution of the system~\eqref{diff.sys} is (asymptotically)
stable, that is, all the solutions of the system~\eqref{diff.sys},
whose initial conditions $x_0$ are sufficiently small, tend to
zero as $t\to+\infty$.

When the parameter $\alpha$ arrives at its critical value, say
$\alpha_1^*$ the following two basic cases can occur, generically:
either a pair of complex conjugate eigenvalues (of multiplicity
one) of the matrix $A(\alpha_1^*)$ appear on the imaginary axis or
an eigenvalue (of multiplicity one) becomes $0$. In the former
case, for $\alpha=\alpha_1^*+\varepsilon$, where $\varepsilon>0$
is sufficiently small, there appears a small-amplitude limit cycle
branching from the zero solution (a fixed point) of the
system~\eqref{diff.sys}, but the zero solution loses its
stability\footnote{Here we assume that for $\alpha>\alpha_1^*$,
the matrix $A(\alpha)$ has an eigenvalue with a positive real
part.}. This type of bifurcation is called \text{Hopf
bifurcation}~\cite{Lyapunoff,Marsden&McCracen,Iooss_Joseph}. In
hydrodynamics~\cite{Yudovich_11_pr}, such type of bifurcation is
called an \textit{oscillatory loss of stability} or
\textit{oscillatory instability}. In the case, when the matrix
$A(\alpha_1^*)$ has all the eigenvalues in the open left
half-plane except one, which is equal to zero and is of
multiplicity one, new fixed points of the system~\eqref{diff.sys}
branch from the zero fixed point. In
hydrodynamics~\cite{Yudovich_11_pr}, such type of bifurcation is
called a \textit{monotone loss of stability} or \textit{monotone
instability}. To determine the type of bifurcation of a fixed
point of a given system is a very important problem in
hydrodynamics and mechanics
(see~\cite{Marsden&McCracen,Yudovich_11_pr} and references there).

Let us consider the characteristic polynomial
$p(z,\alpha)=\det(zE-A(\alpha))$ of the matrix $A(\alpha)$. By
assumption $p(z,\alpha_0)$ is Hurwitz stable. Recall that
$\alpha_1^*>\alpha_0$ is the critical value of the parameter
$\alpha$ such that the zero solution of the
system~\eqref{diff.sys} loses its stability. The
\textit{monotonicity principle} takes place if and only if the
polynomial $p(z,\alpha_1^*)$ is a generalized Hurwitz polynomial,
that is, the inequalities~\eqref{general.Hurwitz.main.det.noneq}
hold. This fact follows from
Theorem~\ref{Theorem.general.Hurwitz.criterion}. Moreover, the
polynomial $p(z,\alpha)$ is generalized Hurwitz of order $0$
(Hurwitz stable) for $\alpha<\alpha_1*$, and $p(z,\alpha_1^*)$ is
a generalized Hurwitz polynomial of order $1$.

Now we describe further changes of zero location of the polynomial
$p(z,\alpha)$ when $\alpha$ goes beyond~$\alpha_1^*$ until
$p(z,\alpha)$ becomes self-interlacing. Let
$\{\alpha_2^{*},\alpha_3^{*},\ldots,\alpha_m^{*}\}$, $2\leqslant
m\leqslant\left[\dfrac{n+1}2\right]$, be the next critical values
of the parameter $\alpha$ at each of which order of the
generalized Hurwitz polynomial $p(z,\alpha)$ increases. This means
that $p(z,\alpha_i^*)$ has a simple root at zero, and for
$\alpha>\alpha_i^*$ sufficiently  close to $\alpha_i^*$, this root
becomes a positive simple root of $p(z,\alpha)$ with the minimal
absolute value (among real roots of $p(z,\alpha)$). Therefore,
when $\alpha$ transfers $\alpha_m^*$, the polynomial $p(z,\alpha)$
has exactly $m$ positive simple zeroes for $\alpha>\alpha_m^*$
sufficiently close to $\alpha_m^*$. Denote the minimal positive
zero of $p(z,\alpha)$ by $\mu_m$. Obviously, if the difference
$\alpha-\alpha_m^*>0$ is sufficiently small, then the number
$\mu_m>0$ is close to zero, so $p(z,\alpha)$ has no zeroes on the
interval $(-\mu_m,0)$. Let $\alpha$ increase. One more simple
positive zero of the polynomial $p(z,\alpha)$ can
appear\footnote{It is possible only if
$m<\left[\dfrac{n+1}2\right]$.}, when two zeroes of a pair of
complex conjugate zeroes of $p(z,\alpha)$ (all its nonreal zeroes
lie in the open left half-plane) face each other on the interval
$(-\mu_m,0)$ to become a real zero of $p(z,\alpha)$ of
multiplicity two. Further, when $\alpha$ grows, this zero splits
into to real simple roots on the interval $(-\mu_m,0)$. One of
those roots goes toward zero and crosses it, when $\alpha$
transfers through the next critical value $\alpha_{m+1}^*$ of the
parameter $\alpha$. Thus, for sufficiently small difference
$\alpha-\alpha_{m+1}^*>0$, this zero becomes the smallest positive
zero $\mu_{m+1}$. The second zero of the mentioned pair of zeroes
remains on the interval $(-\mu_m,-\mu_{m+1})$, since it cannot
cross the value $-\mu_m$ and cannot go away from the real axis.

Let us show that there are no any other scenarios for positive
zeroes to appear if $p(z,\alpha)$ is generalized Hurwitz for all
$\alpha$. In fact, suppose that for $\alpha>\alpha_m^*$, some
negative zeroes of $p(z,\alpha)$ come into the interval
$(-\mu_m,0)$ from the real axis, that is, from the interval
$(-\infty,-\mu_m)$. But in this case, the polynomial $p(z,\alpha)$
will have zeroes $\pm\mu_m$ at some value of $\alpha$, so by the
Orlando's formula~\cite{Gantmakher,Holtz_Tyaglov}, we have
$\Delta_{n+1}(p)=0$. According to
Theorem~\ref{Theorem.general.Hurwitz.criterion}, this contradicts
the assumption that $p(z,\alpha)$ is generalized Hurwitz.

Thus, when $\alpha$ increases from $\alpha_0$, at each interval
$(-\mu_j,-\mu_{j+1})$ the polynomial $p(z,\alpha)$ has an odd
number of zeroes. One of those zeroes appears on
$(-\mu_j,-\mu_{j+1})$ as it was described above, while the other
zeroes appear from complex conjugate pairs. The maximal number of
positive zeroes of $p(z,\alpha)$
equals~$\left[\dfrac{n+1}2\right]$. In this case, $p(z,\alpha)$
becomes self-interlacing.

So the \textit{monotonicity principle} takes place for the
system~\eqref{diff.sys} if and only if for the first critical
value~$\alpha_1^*$ of the parameter $\alpha$ the characteristic
polynomial $p(z,\alpha_1^*)=\det(zE-A(\alpha_1^*))$ is generalized
Hurwitz, that is, the
inequalities~\eqref{general.Hurwitz.main.det.noneq} hold.

\section*{Acknowledgement}
\hspace{4mm} I am grateful to Yu.S.\,Barkovsky and O.\,Holtz for
helpful discussions and for Barkovsky's detailed suggestions for
the improvement of some results in the paper.

\addcontentsline{toc}{section}{Bibliography} 
\bibliographystyle{plain}
\bibliography{References}

\def\cprime{$'$} \def\cprime{$'$} \def\cprime{$'$} \def\cprime{$'$}
  \def\cprime{$'$} \def\cprime{$'$} \def\cprime{$'$}
\begin{thebibliography}{10}

\bibitem{Akhiezer_moment}
N.~I. Akhiezer.
\newblock {\em The classical moment problem and some related questions in
  analysis}.
\newblock Translated by N. Kemmer. Hafner Publishing Co., New York, 1965.

\bibitem{AkhiezerKrein}
N.~I. Akhiezer and M.~G. Krein.
\newblock {\em Some questions in the theory of moments}.
\newblock translated by W. Fleming and D. Prill. Translations of Mathematical
  Monographs, Vol. 2. American Mathematical Society, Providence, R.I., 1962.

\bibitem{Asner}
B.~A. Asner, Jr.
\newblock On the total nonnegativity of the {H}urwitz matrix.
\newblock {\em SIAM J. Appl. Math.}, 18:407--414, 1970.

\bibitem{Atkinson}
F.~V. Atkinson.
\newblock {\em Discrete and continuous boundary problems}.
\newblock Mathematics in Science and Engineering, Vol. 8. Academic Press, New
  York, 1964.

\bibitem{Barkovsky.pr}
Yu.~S. Barkovsky.
\newblock {\em Private communication}, 2004.

\bibitem{Barkovsky.2}
Yu.~S. Barkovsky.
\newblock Lectures on the {R}outh-{H}urwitz problem.
\newblock {\em arXiv:0802.1805}, 2008.

\bibitem{ChebotarevMeiman}
N.~G. {\v{C}}ebotarev and N.~N. Me{\u\i}man.
\newblock The {R}outh-{H}urwitz problem for polynomials and entire functions.
  {R}eal quasipolynomials with {$r=3$}, {$s=1$}.
\newblock {\em Trudy Mat. Inst. Steklov.}, 26:331, 1949.
\newblock Appendix by G. S. Barhin and A. N. Hovanski{\u\i}.

\bibitem{Fisk}
S.~Fisk.
\newblock Polynomials, roots, and interlacing.
\newblock {\em arXiv:0612833}, 2006.

\bibitem{Frobenius}
G.~Frobenius.
\newblock {\"U}ber das {T}r\"agheitsgesetz der quadratischen {F}ormen.
\newblock {\em Sitz.-Ber. Acad. Wiss. Phys.-Math. Klasse, Berlin}, pages
  241--256; 407--431, 1894.

\bibitem{KreinGantmaher}
F.~P. Gantmacher and M.~G. Krein.
\newblock {\em Oscillation matrices and kernels and small vibrations of
  mechanical systems}.
\newblock AMS Chelsea Publishing, Providence, RI, revised edition, 2002.
\newblock Translation based on the 1941 Russian original, Edited and with a
  preface by Alex Eremenko.

\bibitem{Gantmakher.1}
F.~R. Gantmacher.
\newblock {\em The theory of matrices. {V}ol. 1}.
\newblock Translated by K. A. Hirsch. Chelsea Publishing Co., New York, 1959.

\bibitem{Gantmakher}
F.~R. Gantmacher.
\newblock {\em The theory of matrices. {V}ol. 2}.
\newblock Translated by K. A. Hirsch. Chelsea Publishing Co., New York, 1959.

\bibitem{Hermite}
C.~Hermite.
\newblock On the number of roots of an algebraic equation contained between
  given limits.
\newblock {\em Internat. J. Control}, 26(2):183--195, 1977.
\newblock Translated from the French original by P. C. Parks, Routh centenary
  issue.

\bibitem{Holtz1}
O.~Holtz.
\newblock Hermite-{B}iehler, {R}outh-{H}urwitz, and total positivity.
\newblock {\em Linear Algebra Appl.}, 372:105--110, 2003.

\bibitem{H4}
O.~Holtz.
\newblock The inverse eigenvalue problem for symmetric anti-bidiagonal
  matrices.
\newblock {\em Linear Algebra Appl.}, 408:268--274, 2005.

\bibitem{Holtz_Tyaglov}
O.~Holtz and Tyaglov M.
\newblock Structured matrices, continued fractions, and root localization of
  polynomials.
\newblock {\em arXiv:0912.4703}, 2009.

\bibitem{Hurwitz}
A.~Hurwitz.
\newblock \"{U}ber die {B}edingungen, unter welchen eine {G}leichung nur
  {W}urzeln mit negativen reellen {T}heilen besitzt.
\newblock In {\em Stability theory ({A}scona, 1995)}, volume 121 of {\em
  Internat. Ser. Numer. Math.}, pages 239--249. Birkh\"auser, Basel, 1996.
\newblock Reprinted from Math. Ann. {{\bf{4}}6} (1895), 273--284 [JFM
  26.0119.03].

\bibitem{Iooss_Joseph}
G.~Iooss and D.~Joseph.
\newblock {\em Elementary stability and bifurcation theory}.
\newblock Undergraduate Texts in Mathematics. Springer-Verlag, New York, second
  edition, 1990.

\bibitem{KreinKatz}
I.~V. Kac, , and M.~G. Krein.
\newblock {\em On the spectral functions of the string}, volume 103 of {\em
  Amer. Math. Soc. Transl.}
\newblock American Mathematical Society, 1974.

\bibitem{Kemperman}
J.~H.~B. Kemperman.
\newblock A {H}urwitz matrix is totally positive.
\newblock {\em SIAM J. Math. Anal.}, 13(2):331--341, 1982.

\bibitem{Korsakov}
G.~F. Korsakov.
\newblock A generalization of a theorem of {L}i\'enard and {C}hipart.
\newblock {\em Mat. Zametki}, 22(1):13--21, 1977.

\bibitem{KreinNaimark}
M.~G. Kre{\u\i}n and M.~A. Na{\u\i}mark.
\newblock The method of symmetric and {H}ermitian forms in the theory of the
  separation of the roots of algebraic equations.
\newblock {\em Linear and Multilinear Algebra}, 10(4):265--308, 1981.
\newblock Translated from the Russian by O. Boshko and J. L. Howland.

\bibitem{Nudelman}
M.~G. Kre{\u\i}n and A.~A. Nudel{\cprime}man.
\newblock {\em The {M}arkov moment problem and extremal problems}.
\newblock American Mathematical Society, Providence, R.I., 1977.
\newblock Ideas and problems of P. L. {\v{C}}eby{\v{s}}ev and A. A. Markov and
  their further development, Translated from the Russian by D. Louvish,
  Translations of Mathematical Monographs, Vol. 50.

\bibitem{Kronecker}
L.~Kronecker.
\newblock Algebraische {R}eduction der {S}chaaren bilinearer {F}ormen.
\newblock S.-B. Akad. Berlin, 1890.

\bibitem{LienardChipart}
A.~Li\'enard and M.~Chipart.
\newblock Sur la signe de la partie r\'eelle des racines d'une \'equation
  alg\'ebraique.

\bibitem{Lyapunoff}
A.~M. Lyapunov.
\newblock {\em The general problem of the stability of motion}.
\newblock Taylor \& Francis Ltd., London, 1992.
\newblock Translated from Edouard Davaux's French translation (1907) of the
  1892 Russian original and edited by A. T. Fuller, With an introduction and
  preface by Fuller, a biography of Lyapunov by V. I. Smirnov, and a
  bibliography of Lyapunov's works compiled by J. F. Barrett, Lyapunov
  centenary issue, Reprint of Internat. J. Control {{\bf{5}}5} (1992), no. 3 [
  MR1154209 (93e:01035)], With a foreword by Ian Stewart.

\bibitem{Markov}
A.~A. Markov.
\newblock {\em Collected works (Russian)}.
\newblock Academy of Scines USSR, Moscow, 1948.

\bibitem{Marsden&McCracen}
J.~E. Marsden and M.~McCracken.
\newblock {\em The {H}opf bifurcation and its applications}.
\newblock Springer-Verlag, New York, 1976.
\newblock With contributions by P. Chernoff, G. Childs, S. Chow, J. R. Dorroh,
  J. Guckenheimer, L. Howard, N. Kopell, O. Lanford, J. Mallet-Paret, G. Oster,
  O. Ruiz, S. Schecter, D. Schmidt and S. Smale, Applied Mathematical Sciences,
  Vol. 19.

\bibitem{Pick}
G.~Pick.
\newblock {\"U}ber die {B}eschr{\"a}nkungen analytische {F}unktionen, welche
  durch vorgegebene {F}unktionswerte bewirkt werden.
\newblock {\em Math. Ann.}, 77:7--23, 1916.

\bibitem{Pivovarchik_2007}
V.~Pivovarchik.
\newblock On spectra of a certain class of quadratic operator pencils with
  one-dimensional linear part.
\newblock {\em Ukra\"\i n. Mat. Zh.}, 59(5):702--716, 2007.

\bibitem{Pivovarchik_gen_Hurw}
V.~Pivovarchik.
\newblock Symmetric {H}ermite-{B}iehler polynomials with defect.
\newblock In {\em Operator theory in inner product spaces}, volume 175 of {\em
  Oper. Theory Adv. Appl.}, pages 211--224. Birkh\"auser, Basel, 2007.

\bibitem{Mee_Piv_2001}
V.~Pivovarchik and C.~van~der Mee.
\newblock The inverse generalized {R}egge problem.
\newblock {\em Inverse Problems}, 17(6):1831--1845, 2001.

\bibitem{Pivovarchik_Worachek_1}
V.~Pivovarchik and H.~Woracek.
\newblock Shifted {H}ermite-{B}iehler functions and their applications.
\newblock {\em Integral Equations Operator Theory}, 57(1):101--126, 2007.

\bibitem{Pivovarchik_Worachek_2}
V.~Pivovarchik and H.~Woracek.
\newblock The square-transform of {H}ermite-{B}iehler functions. {A} geometric
  approach.
\newblock {\em Methods Funct. Anal. Topology}, 13(2):187--200, 2007.

\bibitem{Routh}
E.~J. Routh.
\newblock {\em A treatise on the stability of a given state of motion}.
\newblock London, 1877.

\bibitem{Sheil-Small}
T.~Sheil-Small.
\newblock {\em Complex polynomials}, volume~75 of {\em Cambridge Studies in
  Advanced Mathematics}.
\newblock Cambridge University Press, Cambridge, 2002.

\bibitem{Simon_1998}
B.~Simon.
\newblock The classical moment problem as a self-adjoint finite difference
  operator.
\newblock {\em Adv. Math.}, 137(1):82--203, 1998.

\bibitem{Stieltjes1}
T.~J. Stieltjes.
\newblock Recherches sur les fractions continues.
\newblock {\em Ann. Fac. Sci. Toulouse Math. (6)}, 4(3):J76--J122, 1995.
\newblock Reprint of Ann. Fac. Sci. Toulouse {{\bf{8}}} (1894), J76--J122.

\bibitem{Stieltjes2}
T.~J. Stieltjes.
\newblock Recherches sur les fractions continues.
\newblock {\em Ann. Fac. Sci. Toulouse Math. (6)}, 4(4):A5--A47, 1995.
\newblock Reprint of Ann. Fac. Sci. Toulouse {{\bf{9}}} (1895), A5--A47.

\bibitem{Mee_Piv_2005}
C.~van~der Mee and V.~Pivovarchik.
\newblock The {S}turm-{L}iouville inverse spectral problem with boundary
  conditions depending on the spectral parameter.
\newblock {\em Opuscula Math.}, 25(2):243--260, 2005.

\bibitem{Wall}
H.~S. Wall.
\newblock {\em Analytic {T}heory of {C}ontinued {F}ractions}.
\newblock D. Van Nostrand Company, Inc., New York, 1948.

\bibitem{Yudovich_11_pr}
V.~I. Yudovich.
\newblock Eleven great problems of mathematical hydrodynamics.
\newblock {\em Mosc. Math. J.}, 3(2):711--737, 746, 2003.
\newblock Dedicated to Vladimir I. Arnold on the occasion of his 65th birthday.

\end{thebibliography}

\end{document}